\input amstex
\documentstyle{amsppt}
\magnification=1200
\nopagenumbers
\hsize=470 true pt
\vsize=640 true pt

\document

\topmatter
\title  A standard zero free region for Rankin Selberg L-functions  \endtitle 
\author  Dorian Goldfeld, \; Xiaoqing Li
 \endauthor
 
 \abstract
A standard zero free region is obtained for Rankin Selberg L-functions 
$L(s, f\times \widetilde{f})$ where $f$ is an almost everywhere tempered Maass form on $GL(n)$ and $f$ is not necessarily self dual. The method is based on the theory of Eisenstein series generalizing a work of Sarnak.  \endabstract

\address
Dorian Goldfeld\hfill\break 
\phantom{24.}Department of Mathematics\hfill\break
\phantom{24.}Columbia University\hfill\break
\phantom{24.}Rm 422, MC4428\hfill\break
\phantom{24.}2990 Broadway\hfill\break
\phantom{24.}New York, NY 10027
\endaddress

\email
goldfeld\@columbia.edu\endemail

\address
Xiaoqing Li,\hfill\break 
\phantom{24.}Department of Mathematics\hfill\break
\phantom{24.}244 Mathematics Building\hfill\break
\phantom{24.}University at Buffalo\hfill\break
\phantom{24.}Buffalo, NY 14260-2900
\endaddress

\email
xl29\@buffalo.edu\endemail

\endtopmatter

\noindent
{\bf \S 1 Introduction}
\vskip 10pt
Let $\pi_1, \pi_2$ be cuspidal automorphic representations of $GL(n_i, \Bbb A_F)$ $(i = 1,2$) for some number field $F$ whose central characters are
trivial on $\Bbb R^+$  imbedded diagonally in the (archimedean) id\`eles. As pointed out in \cite{Moreno, 1985},  \cite{Sarnak, 2004}, de La Vall\'ee Poussin's proof of the prime number theorem can be applied to any Rankin-Selberg L-function $L(s, \pi_1\times \pi_2)$
 provided   one of the $\pi_i$ ($i = 1,2$) is self dual. In this case, the  zero free region for $L(\sigma+it, \pi_1\times \pi_2)$ takes the form
 $$\sigma \; > \; 1 - \frac{c}{\log\left(Q_{\pi_1} Q_{\pi_2}\cdot \big(|t|+2\big)   \right)},$$
 where $c > 0$ is an absolute constant depending only on $n_1, n_2, F$ and $Q_{\pi_1}, Q_{\pi_2}$   are the analytic conductors (see \cite{Iwaniec-Sarnak, 1999} for the definition of analytic conductor) of $\pi_1, \pi_2,$ respectively. A zero free region of the form
  $$\sigma \; > \; 1 - \frac{c}{\left(\log\big(Q_{\pi_1} Q_{\pi_2}\cdot \big(|t|+2\big)   \big)\right)^B},$$
 for some fixed $B > 0$ is called a standard zero free region.
 
Prior to this work, a standard zero free region for $L(s, \pi_1\times \pi_2)$ was not known in the case that both $\pi_1, \pi_2$ are not self dual. The best zero-free region (non self dual case)  known to date is due to \cite{Brumley, 2006}, and is of the form
 $$\sigma \; > \; 1 - \frac{c}{\big(Q_{\pi_1} Q_{\pi_2}\cdot \big(|t|+2\big)\big)^N},\tag{1.1}$$
 where $c > 0, N>1$ depend only on $n_1, n_2, F.$ 
 The bound (1.1) was generalized to L-functions of Langlands-Shahidi type in \cite{Gelbart-Lapid, 2006}. They invoked the method of \cite{Sarnak, 2004} which made use of Eisenstein series, the Maass--Selberg relations, and a sieving argument. Sarnak's work can be viewed as an effectuation of the non-vanishing results of \cite{Jacquet-Shalika, 1976} for the standard L-function.

  In this paper we also follow \cite{Sarnak, 2004}, obtaining a much stronger result than (1.1), i.e., a standard zero free region, for a restricted class of Rankin-Selberg L-functions.  Our main theorems are the following.
  
  \vskip 10pt\noindent
{\bf Theorem 1.2} {\it  Let $\pi$ be an irreducible cuspidal  unramified representation of $GL(n, \Bbb A_{\Bbb Q})$ for $n \ge 2$ which is tempered at all finite primes (except possibly a set of measure zero). Let $\widetilde \pi$ be the contragredient representation of $\pi.$ Then
a zero free region for $L\left(\sigma+it, \pi\times\widetilde \pi\right)$ is given by
$$\sigma \; > \; 1 \; - \; \frac{c}{\big(\log (|t|+2)\big)^5 },$$
where $c > 0$ is a fixed constant (independent of $t$) which depends at most on $\pi.$}

\vskip 10pt\noindent
{\bf Theorem 1.3} {\it  Let $\pi$ be an irreducible cuspidal  unramified representation of $GL(n, \Bbb A_{\Bbb Q})$ for $n \ge 2$ which is tempered at all finite primes (except possibly a set of Dirichlet density zero). Let $f$ be a Maass cusp form in the space of $\pi$  with dual  $\widetilde f \in \widetilde \pi.$ Then
 $$|L(1+it, \; f\times \widetilde f \,) |\; \gg\; \frac{1}{ \big(\log (|t|+2)\big)^3  }\tag{1.4}$$
for all $t\in \Bbb R$ and where the  constant implied by the $\gg$-symbol is effectively computable and depends at most on $f.$
}

\vskip 10pt
Theorem 1.2 follows from the mean value theorem together with theorem 1.3 and the fact that
$L'(\sigma+it, f\times \widetilde f\,) \ll_{\pi,c} \big(\log (|t|+2)\big)^2$ for $1 \, - \,c/\big(\log (|t|+2)\big)^5 \le \sigma \le 1.$ When $|t|$ is less than any fixed positive constant, the lower bound (1.4) was already proved in \cite{Brumley, 2012}. To simplify the exposition we shall assume $|t| \gg 1$  for the remainder of this paper, where the constant  implied by the $\gg$-symbol is sufficiently large and  effectively computable.
 
\vskip 6pt
Theorem 1.3 can be viewed as an effectuation of Shahidi's non-vanishing result \cite{Shahidi, 1981} and an improvement of (1.1).
Theorem 1.3 will have many applications since  the  relative trace formula for $GL(2n)$ (see \cite{Jacquet-Lai, 1985}, \cite{Lapid, 2006}) will have a spectral contribution which involves integrating the Rankin-Selberg L-function  $L(1+it, \pi\times\widetilde{\pi})^{-1}$ over $t\in\Bbb R$. It is crucial to have good lower bounds for $L(1+it, \pi\times \widetilde{\pi}).$
\vskip 8pt
The key step in the proof of theorem 1.3 is the introduction of the integral $\Cal I$ which is defined as  $|L(1+2int, \; f\times \widetilde f \,)|^2$ multiplied by a certain transform of a smoothly truncated Eisenstein series on $GL(2n, \Bbb A_\Bbb Q).$ The precise definition of $\Cal I$ is given in 11.1. By the use of the Maass-Selberg relation it is possible to obtain an almost sharp upper bound for $\Cal I.$ The upper bound is given in theorem 11.2.  Remarkably, it is also possible to obtain an almost sharp lower bound for $\Cal I$ using bounds for Whittaker functions and a sieving argument. The lower bound is given in theorem 12.1. The combination of the upper and lower bounds for $\Cal I$ immediately prove theorem 1.3.
\vskip 6pt
We believe that the methods of this paper can be vastly generalized to  global (not  necessarily unramified) irreducible  cuspidal automorphic representations of reductive groups over a number field, leading to standard zero free regions for all Rankin-Selberg L-functions.  It is very likely that the assumptions of temperedness at finite primes can also be dropped.

  \vskip 6pt
   To keep technicalities as simple as possible we work over $\Bbb A_{\Bbb Q}$ with the assumption that $\pi$ is globally unramified. In this case, it is only necessary to do archimedean computations. The assumption that $\pi$ is globally unramified can be removed. It is not hard to modify the proof we present if $\pi$ is ramified at $\infty.$ If, in addition, $\pi$ is ramified at a finite set of  primes then the entire proof will go through except that the $\gg$-constant in the lower bound for $|L(1+it, \; f\times \widetilde f \,) |$ will now depend  on the finite set of ramified primes. This is clear because the key integrals coming up in the proof will all be Eulerian and the ramified primes will only involve a finite Euler product. So bounds for these integrals will only change by a finite factor when doing the proof in the ramified case.

\vskip 20pt
\noindent
  {\bf \S 2 Basic Notation:} 
 \vskip 10pt 
    We shall be working  with
    $$G := GL(n,\Bbb R), \;\;\quad \Gamma := GL(n, \Bbb Z), \;\;\quad K := O(n, \Bbb R),$$
    for $n\ge 2.$ Every $z \in GL(n, \Bbb R)\big /(K\cdot\Bbb R^\times)$ can be given in Iwasawa form $$z = xy\tag{2.1}$$     $$x =  \pmatrix 1 & x_{1,2}& x_{1,3}& \cdots  & &
x_{1,n}\\
 & 1& x_{2,3} &\cdots & & x_{2,n}\\
& &\hskip 2pt \ddots & & & \vdots\\
& && & 1& x_{n-1,n}\\
& & & & &1\endpmatrix, \;\;\;\;\;\;\;\;
y =  
\pmatrix y_1y_2\cdots
y_{n-1} & & & \\
& \hskip -30pt y_1y_2\cdots y_{n-2} & & \\
& \ddots &  & \\
& & \hskip -5pt y_1 &\\
& & &  1\endpmatrix.$$
with $x_{i,j} \in \Bbb R$ ($1\le i < j\le n$) and $y_i > 0$ ($1 \le i \le n-1$).
\vskip 8pt
For $\nu = (\nu_1, \ldots, \nu_{n-1}) \in \Bbb C^{n-1}$ consider the $I_\nu$-function
$$I_\nu(z) := \prod_{i=1}^{n-1} \prod_{j=1}^{n-1}
y_i^{b_{i,j} \nu_j},\tag{2.2}$$ where
$$b_{i,j} = \cases i j & \text{if $i + j \le n,$}\\
(n-i)(n-j) & \text{if $i + j \ge n.$}\endcases$$

Associated to $\nu\in\Bbb C^{n-1}$ there are Langlands parameters $\alpha_1, \ldots,\alpha_n \in \Bbb C$ defined as follows. For each $1 \le i\le n-1$ define
$$B_i(\nu) := \sum_{j=1}^{n-1} b_{i,j} \nu_j.$$
Then the Langlands parameters associated to $(n, \nu)$  are given by
$$\alpha_i := \cases B_{n-1}(\nu) +\frac{1-n}{2}, & \text{if} \; i = 1,\\
B_{n-i}(\nu) - B_{n-i+1}(\nu) + \frac{2i-n-1}{2}, & \text{if}\; 1 < i < n,\\
-B_1(\nu) + \frac{n-1}{2}, &\text{if} \; i = n.\endcases \tag{2.3}$$

Fix a partition
$n = n_1+n_2 + \cdots +n_r$ with $1 \le n_1,n_2, \ldots, n_r < n$. We define the standard parabolic subgroup
$$\Cal P := P_{n_1,n_2, \ldots,n_r} := \left\{\pmatrix GL(n_1) & * & \cdots &*\\
0 & GL(n_2) & \cdots & *\\
\vdots & \vdots & \ddots & \vdots \\
0 & 0 &\cdots & GL(n_r)\endpmatrix\right\}\tag{2.4}$$
with nilpotent radical $$N^{\Cal P} := \left\{\pmatrix I_{n_1} & * & \cdots &*\\
0 & I_{n_2} & \cdots & *\\
\vdots & \vdots & \ddots & \vdots \\
0 & 0 &\cdots & I_{n_r}\endpmatrix\right\},\tag{2.5}$$
and standard Levi
$$M^{\Cal P} := \left\{\pmatrix GL(n_1) & 0 & \cdots &0\\
0 & GL(n_2) & \cdots & 0\\
\vdots & \vdots & \ddots & \vdots \\
0 & 0 &\cdots & GL(n_r)\endpmatrix\right\}.\tag{2.6}$$

 Define $$\frak h^n := G/(K\cdot\Bbb R^\times).$$For a function $F: G/ KR^\times \to \Bbb C$ and $\gamma \in \Gamma$, we define the slash operator by
$$F(z) \;\big |_\gamma := F(\gamma z), \qquad (z \in G/ KR^\times).\tag{2.7}$$

For two functions  $F_1, F_2  \in \Cal L^2\big(SL(n,\Bbb Z)\backslash\frak h^n\big)$, we define the inner product
$$\langle F_1, F_2\rangle := \int\limits_{SL(n, \Bbb Z)\backslash \frak h^n} F_1(z) \overline{F_2(z)} \; d^* z,\tag{2.8}$$
where $d^* z := d^\times x \;d^\times y = \prod\limits_{1 \le i < j \le n} dx_{i,j} \prod\limits_{k=1}^{n-1} y_k^{-k(n-k)-1} dy_k$ is the  left invariant measure.

\pagebreak

\vskip 20pt
\noindent
{\bf \S 3 The maximal parabolic Eisenstein series $\pmb E_{\pmb P_{\pmb 2\pmb n-\pmb1, \pmb1}}$:}

\vskip 10pt Let $z \in GL(2n, \Bbb R)\big /(K\cdot\Bbb R^\times)$ and $s \in   \Bbb C$ with $\Re(s) > 1$. Let $P_{2n-1,1}$ denote the maximal parabolic subgroup as in
(2.4). We may then define as in proposition 10.7.5 of \cite{Goldfeld 2006} the maximal parabolic Eisenstein series
$$E_{P_{2n-1,1}}(z, s) := \sum_{\gamma\, \in \, (P_{2n-1,1}(\Bbb Z)\cap SL(2n,\Bbb Z))\backslash SL(2n, \Bbb Z)} \text{Det}(\gamma z)^s,$$
and the completed maximal parabolic Eisenstein series
$$E^*_{P_{2n-1,1}}(z,s) = \pi^{-ns} \Gamma(ns) \zeta(2ns) E_{P_{2n-1,1}}(z,s) = E^*_{P_{2n-1,1}}(\phantom{.}^tz^{-1},1-s)$$
where $\phantom{.}^tz$ denotes the transpose of the matrix $z$.

 For $u\in \Bbb R$ with $u > 0$, and $z = xy$ as in (2.1),  define the theta function
 $$\theta_z(u) := \sum_{(a_1, a_2, \ldots, a_{2n})\in \Bbb Z^{2n}} e^{-\pi\left(b_1^2 + b_2^2 + \cdots + b_{2n}^2   \right)\cdot u},\tag{3.1}$$
 where
 $$\align
 & b_1 = a_1Y_1,\tag{3.2}\\
 & b_2 = \left(a_1 x_{1,2} + a_2\right) Y_2\\
 & \phantom{xxi}\vdots\\
 & b_k = \left(a_1 x_{1,k} + a_2 x_{2,k} + \cdots + a_{k-1} x_{k-1,k} + a_{k}   \right)Y_k\\
& \phantom{xxi}\vdots\\
& b_{2n} = \left(a_1 x_{1,2n} + a_2 x_{2,2n} + \cdots + a_{2n}\right) Y_{2n}
  \endalign$$
  and
  $$Y_k := y_1y_2 \cdots y_{2n-k} \left[y_1^{2n-1} y_2^{2n-2} \cdots y_{2n-1}  \right]^{-\frac{1}{2n}}, \qquad (\text{for} \; 1 \le k \le 2n).$$
  In the formula above for $Y_k$ we define $y_0 = 1.$
  
    It then follows as in the proof of  proposition 10.7.5 \cite{Goldfeld 2006} that
    $$\align E^*_{P_{2n-1,1}}(z,s) &= \int\limits_1^\infty \left[\theta_z(u) - 1   \right] u^{ns} \; \frac{du}{u} \; + \; \int\limits_1^\infty \left[ \theta_{w^tz^{-1}w}(u) - 1  \right] \, u^{n(1-s)} \; \frac{du}{u}\tag{3.3}\\
    & \hskip 180pt-\frac{1}{n}\left(\frac{1}{1-s} \; + \; \frac{1}{s}  \right).\endalign$$
    where $w$ is the long element in the Weyl group.
  We shall use (3.3) to prove the following proposition.
  
  \vskip 10pt\noindent
  {\bf Proposition 3.4} {\it Let  $z \in GL(2n, \Bbb R)\big /(K\cdot\Bbb R^\times)$  as in (2.1). Let $y_0=1$ and assume that $y_i \ge 1 \;(1 \le i \le 2n-1)$. 
 Then for $w \in \Bbb C$ with $\Re(w) = \frac12$, we have
 $$\align & \big|E_{P_{2n-1,1}}(z,w)\big| \,\ll \,\frac{e^{\frac{\pi}{2}n|w|}}{|w|^{\frac{n-1}{2}}}\sum_{1\le k \le 2n}\Bigg[ \Big(y_1 y_2^2 \; \cdots \; y_{2n-k}^{2n-k}    \Big)^{\frac12}\Big(y_{2n-k+1}^{k-1} y_{2n-k+2}^{k-2} \;\;\cdots\;\; y_{2n-1}   \Big)^{ \frac12  }\\
 &\hskip 36pt + \Big(y_1 y_2^2 \; \cdots \; y_{2n-k}^{2n-k}    \Big)^{\frac{k-1}{2n}} \Big(y_{2n-k+1}^{k-1} y_{2n-k+2}^{k-2} \;\;\cdots\;\; y_{2n-1}   \Big)^{ \frac{2n-k+1}{2n}  }\Bigg] \cdot \log(y_1\cdots y_{2n-1}),\endalign$$
 where the implied constant only depends on $n$. }

 \vskip 10pt
 {\bf Proof of Proposition 3.4:} The factor $\frac{e^{\frac{\pi}{2}n|w|}}{|w|^{\frac{n-1}{2}}}$ comes from Stirling's asymptotic formula for $|\Gamma(nw)|^{-1},$ which arises when passing from $E^*_{P_{2n-1,1}}(z,w)$ to $E_{P_{2n-1,1}}(z,w)$.
 
  It follows from (3.3) that for  $w \in \Bbb C$ with $\Re(w) = \frac12$, to bound $E_{P_{2n-1,1}}(z,w)$, it is enough to bound the integral
 $$\int\limits_1^\infty \big|\theta_z(u) - 1   \big|\, u^{n/2} \; \frac{du}{u}.$$
 The other integral can be bounded similarly.
 
 First of all, we note that by (3.1), (3.2), we have
 $$\theta_z(u) - 1 \; = \; \sum_{k=1}^{2n} \underset a_i=0\; (1\le i\le k-1) \to{\sum_{a_k\ne 0}} \sum_{a_{k+1}, \ldots,a_{2n}\in\Bbb Z} e^{-\pi u\left[ a_k^2Y_k^2     \;+\; b_{k+1}^2 \; + \; \cdots \; + \; b_{2n}^2\right] }.\tag{3.5}$$
 
  Observe that for $t > 0$ and $\alpha\in\Bbb R$, we may compute (using the Poisson summation formula) that
 $$\sum_{n\in\Bbb Z} e^{-(n+\alpha)^2 t}  = \sum_{n\in\Bbb Z}\; \int\limits_{-\infty}^{\infty} e^{-(x+\alpha)^2 t } e^{2\pi i n x} \, dx \;  \ll \;1 + \frac{1}{t^\frac12}.\tag{3.6}$$ 
 It immediately follows from (3.2), (3.5), (3.6) that
 $$\big| \theta_z(u) - 1\big| \ll \sum_{k=1}^{2n-1} \sum_{a_k\ne 0} e^{-\pi u a_k^2Y_k^2} \prod_{i=k+1}^{2n} \left( 1 + \frac{1}{Y_i\cdot u^\frac12}\right) \; + \; \sum_{a_{2n}\ne 0} e^{-\pi u \, a_{2n}^2 Y_{2n}^2}.$$
 
 Consequently
  $$\align \int\limits_1^\infty \big|\theta_z(u) - 1   \big|\, u^{n/2} \; \frac{du}{u} \; & \ll  \; \sum_{k=1}^{2n-1} \,\int\limits_{Y_k^2}^\infty  \sum_{a_k\ne 0} e^{-\pi u a_k^2} \prod_{i=k+1}^{2n} \left( 1 + \frac{1}{Y_i Y_k^{-1}\cdot u^\frac12}\right) \,Y_k^{-n} \,u^{n/2} \; \frac{du}{u}
  \\
  &
  \\
  & +  \int\limits_{Y_{2n}^2}^\infty  \sum_{a_{2n}\ne 0} e^{-\pi u a_{2n}^2} \,Y_{2n}^{-n} \, u^{n/2} \; \frac{du}{u}.\endalign$$
 
 We now define
 $$\Cal I_k(z) := \cases \int\limits_{Y_k^2}^\infty  \sum\limits_{a_k\ne 0} e^{-\pi u a_k^2} \prod\limits_{i=k+1}^{2n} \left( 1 + \frac{1}{Y_i Y_k^{-1}\cdot u^\frac12}\right) \,Y_k^{-n} \,u^{n/2} \; \frac{du}{u}, & \text{if} \; 1 \le k \le 2n-1,\\
 \int\limits_{Y_{2n}^2}^\infty \;  \sum\limits_{a_{2n}\ne 0} e^{-\pi u a_{2n}^2} \,Y_{2n}^{-n} \, u^{n/2} \; \frac{du}{u}, & \text{if} \; k = 2n.\endcases$$
 so that
$$\int\limits_1^\infty \big|\theta_z(u) - 1   \big|\, u^{n/2} \; \frac{du}{u} \; \ll  \; \sum_{k=1}^{2n} \Cal I_k(z).$$
 The bound for $\Cal I_k(z)$ will depend on whether $Y_k \ge 1$ or $Y_k < 1.$ Clearly
 $$Y_1 \ge Y_2 \ge \cdots \ge Y_{2n}.$$

 \vskip 10pt
 $\underline{\text{\bf Case 1: $Y_k \ge 1$ and $1\le k \le 2n-1$}}$
 \vskip 6pt
  In this case we have
 $$\align  \Cal I_k(z) & \ll \int\limits_{1}^\infty  \sum_{a_k\ne 0} e^{-\pi u a_k^2}   \left[1 + Y_{k+1}^{-1} Y_k\right]\; \cdots \;
 \left[1 + Y_{2n}^{-1} Y_k \right] \, Y_k^{-n}\, u^{n/2}  \; \frac{du}{u}\tag{3.6}
 \\
 &
 \ll \;Y_{k+1}^{-1} \cdots Y_{2n}^{-1} \, Y_k^{n-k} .
 \endalign$$
 since $Y_j^{-1} Y_k \ge 1$ for $j > k.$

  \vskip 20pt
 $\underline{\text{\bf Case 2: $Y_k < 1$ and $1 \le k \le 2n-1$}}$
 \vskip 6pt
 
   We write
   $$\align \Cal I_k(z) &= \left( \int\limits_{Y_k^2}^1  \;\; + \;\; \int\limits_1^\infty\right)  \sum_{a_k\ne 0} e^{-\pi u a_k^2} \prod_{i=k+1}^{2n} \left( 1 + \frac{1}{Y_i Y_k^{-1}\cdot u^\frac12}\right) \,Y_k^{-n} \,u^{n/2} \; \frac{du}{u}\\
   &:= \; I_1 + I_2. \endalign$$
 
 Note that for $u \le 1$, we have
 $$\sum_{a_k\ne 0} e^{-\pi a_k^2 u} \ll u^{-\frac12}.$$ 
 It follows that 
 $$\align  I_1 & = \int\limits_{Y_k^2}^1 \sum_{a_k\ne 0} e^{-\pi u a_k^2} \prod_{i=k+1}^{2n} \left( 1 + \frac{1}{Y_i Y_k^{-1}\cdot u^\frac12}\right) \,Y_k^{-n} \,u^{n/2} \; \frac{du}{u} \;
 \tag{3.7}\\
 &
  \ll \;\int\limits_{Y_k^2}^1 u^{-\frac12} \Big( Y_{k+1}^{-1} Y_k u^{-\frac12} \Big)\; \cdots  \; \Big( Y_{2n}^{-1} Y_k u^{-\frac12} \Big)\, Y_k^{-n} \, u^{n/2} \;  \frac{du}{u}
 \\
 & 
  \ll \;\int\limits_{Y_k^2}^1  \Big( Y_{k+1}^{-1} Y_{k+2}^{-1} \; \cdots \; Y_{2n}^{-1} \Big)  \, Y_k^{n-k} \, u^{-\frac32 - \frac{n}{2} + \frac{k}{2}} \;  du
  \\
  &
 \ll \cases Y_{k+1}^{-1} \cdots Y_{2n}^{-1} \, Y_k^{n-k} \cdot \log(y_1 \cdots y_{2n-1}), & \text{if} \; k = n+1,\\
 Y_{k+1}^{-1} \cdots Y_{2n}^{-1} \, Y_k^{-1}, & \text{if} \; k \ne n+1.\endcases
     \endalign$$
     
     We also have
     $$\align I_2 & = \int\limits_1^\infty \sum_{a_k \ne 0} e^{-\pi a_k^2 u} \left[1 + Y_{k+1}^{-1} Y_k\, u^{-\frac12}\right] \cdots 
         \left[1 + Y_{2n}^{-1} Y_k\, u^{-\frac12}\right] \cdot u^{\frac{n}{2}} \, Y_k^{-n} \;\frac{du}{u}\tag{3.8}\\
         &
         \\
         & \ll Y_{k+1}^{-1} \cdots Y_{2n}^{-1} \, Y_k^{n-k}.\endalign$$
         
         \vskip 20pt
 $\underline{\text{\bf Case 3: $k = 2n$}}$

 \vskip 6pt
Clearly $Y_{2n} \le 1.$  We compute
 $$\align\Cal I_{2n}(z)  & = \int\limits_{Y_{2n}^2}^\infty \sum_{a_{2n}\ne 0} e^{-\pi a_{2n}^2  u} \, Y_{2n}^{-n} \,u^{n/2} \; \frac{du}{u}
 \\
 & = \int\limits_{Y_{2n}^2}^1 \;\; + \;\; \int\limits_1^\infty
 \endalign$$
 where
 $$\align \int_{Y_{2n}^2}^1 & \ll \;\int_{Y_{2n}^2}^1 u^{-\frac12+\frac{n}{2}-1} \; du \cdot Y_{2n}^{-n}
\; \ll \; \left(1 + Y_{2n}^{n-1}\right) \cdot Y_{2n}^{-1} \; \ll \;Y_{2n}^{-n}
 \\
 \endalign$$
 and
 $$\int\limits_1^\infty \ll \;Y_{2n}^{-n}.$$
 Hence
 $$\Cal I_{2n}(z)\; \ll \;Y_{2n}^{-n}.\tag{3.9}$$
 
  To complete the proof of Proposition 3.4, we make use of the following two identities.
 $$Y_{k+1}^{-1}\cdots Y_{2n}^{-1} \, Y_k^{n-k} = \left(y_1 y_2^2 \cdots y_{2n-k}^{2n-k}   \right)^\frac12 \left(y_{2n-k+1}^{k-1} \, y_{2n-k+2}^{k-2} \cdots y_{2n-1}\right)^\frac12. \tag{3.10}$$ 
 $$Y_{k+1}^{-1}\cdots Y_{2n}^{-1} \, Y_k^{-1} = \left(y_1 y_2^2 \cdots y_{2n-k}^{2n-k}   \right)^\frac{k-1}{2n} \left(y_{2n-k+1}^{k-1} \, y_{2n-k+2}^{k-2} \cdots y_{2n-1}\right)^{ \frac{2n-k+1}{2n} }. \tag{3.11}$$

  The proof of Proposition 3.4 then follows form (3.6), (3.7), (3.8), (3.9) in conjunction with (3.10), (3.11). 
 \qed
 
 \vskip 20pt
 \noindent
  {\bf \S 4 Eisenstein series associated to the $P_{n,n}$ parabolic on $GL(2n)$:} 
 \vskip 10pt 

 We  fix an even Maass form (normalized with first Fourier coefficient = 1)
$$f:GL(n, \Bbb Z)\big\backslash GL(n, \Bbb R)\big/(O(n, \Bbb R)\cdot \Bbb R^\times) \; \to \; \Bbb C,\tag{4.1}$$ of type $\nu = (\nu_1, \ldots, \nu_{n-1}) \in \Bbb C$
as in definition 5.1.3 \cite{Goldfeld, 2006}. Associated to $f$ there are Langlands parameters $\alpha_1, \ldots,\alpha_n \in \Bbb C^n$ as in (2.3).

Also associated to $f$ we have the L-function 
$$L(s, f) = \prod_p \prod_{i=1}^n \left(1 -\frac{\alpha_{p,i}}{p^s}   \right)^{-1}, \qquad \left(\alpha_{p,i}\in\Bbb C\right),$$ 
and the Rankin-Selberg L-function 
$$L(s, f\times \widetilde f\,) = \prod_p\prod_{i=1}^n\prod_{j=1}^n \left(1 -\frac{\alpha_{p,i} \overline{\alpha_{p,j}}}{p^s}   \right)^{-1}$$
 which satisfy the functional equations
$$\align \Lambda(s, f) & = \pi^{-\frac{ns}{2}} \prod_{i=1}^n \Gamma\left(  \frac{s+\alpha_i}{2} \right) L(s, f) = \epsilon(f) \cdot \Lambda(1-s, \widetilde f),\tag{4.2}\\
\Lambda\left(s, \,f\times \widetilde f \,\right) & = \pi^{-\frac{n^2 s}{2}} \prod_{i=1}^{n} 
\prod_{j=1}^{n}\Gamma\left(  \frac{s+\alpha_i + \overline{\alpha_j}}{2} \right) L(s, f\times \widetilde f\, ) = \epsilon(f\times \widetilde f\,) \cdot\Lambda(1-s, \widetilde f\times  f),
\endalign$$
with root numbers  $\epsilon(f), \,\epsilon(f\times \widetilde f\,)$ of absolute value one, and where $\widetilde f$ denotes the dual Maass form,
i.e., $\widetilde f(z) =  f(w \phantom{i}^t(z^{-1}))$ with $w$ the long element of the Weyl group.
\vskip 10pt
\noindent
{\bf Remark 4.3:} {\it The completed Rankin-Selberg L-function given by $\Lambda(s, f\times f)$ is entire if $f$ is not self dual (see \cite{Jacquet, 1972} and \cite{Moeglin-Waldspurger, 1989} for a proof).}

\vskip 5pt
Let $s \in \Bbb C$ with $\Re(s) \gg 1.$ For the parabolic
$\Cal P  = P_{n,n}$, with nilpotent radical $N^{\Cal P}$ and standard Levi $M^{\Cal P}$, define
$$\phi_s(\frak n \, \frak  m \, k) :=  \left |\frac{ \text{det}( \frak m_1)}{\text{det}(\frak m_2)}\right |^{n s} f(\frak m_1)  f(\frak m_2),\tag{4.4}$$
where  $$\frak n \in N^{\Cal P}, \qquad \frak m = \pmatrix \frak m_1 &\\& \frak m_2\endpmatrix \in M^{\Cal P}, \qquad k \in K.$$
Let $z \in GL(2n, \Bbb R)/K R^\times$. By the Iwasawa decomposition, $z$ lies in the minimal parabolic, and hence lies in every standard parabolic, so $z$ has a decomposition of the form $z = \frak n \frak m  k$ as above.  We may then define the Eisenstein series 
$$E(z,f; s) := \sum_{\gamma \,\in\, (P_{n,n}(\Bbb Z)\cap\Gamma)\backslash \Gamma} \phi_s(\gamma z).\tag{4.5}$$
Let
 $$C_{\Cal P}\, E(z,f;s) := \phi_s(z) + \frac{\Lambda(2ns - n, f\times \widetilde f\,)}{\Lambda(1+2ns-n, f\times \widetilde f\,)} \left| \frac{ \text{det}(\frak m_1)  }{ \text{det}(\frak m_2)  }  \right|^{n(1-s)}\hskip -2pt f(\frak m_1)\, f(\frak m_2)\tag{4.6}$$
 denote the constant term of $E$ along the parabolic $\Cal P$ (see proposition 8.5). It is convenient to introduce the height function
  $$h(z) := \left | \frac{\text{det}(\frak m_1)}{\text{det}(\frak m_2)}  \right |.\tag{4.7}$$
Define $$\widehat{E}_{A}(z, f; s) := E(z,f; s) \; - \;  \underset
 h(\gamma z) \ge A\to {\sum_{\gamma \, \in\, (P_{n,n}(\Bbb Z)\cap\Gamma)\backslash\Gamma}}C_{\Cal P} \, E(z,f; s) \; \bigg |_\gamma\tag{4.8}$$
 to be the truncated Eisenstein series in the sense of Arthur \cite{Arthur, 1980}. It appears difficult to compute a useful Fourier series expansion of Arthur's truncation, so we introduce a modified version.

 \vskip 10pt\noindent
 {\bf Definition 4.9 (Smoothed Arthur truncation of $E$):} Let $A \ge 1.$ We define
 $$\align & \widehat E_A^*(z, f;  s) := E(z,f; s) -A^{\frac{n}{2}} E(z, f; s-1/2) + \frac{\Lambda(2ns - 2n, f\times \widetilde f\,)}{\Lambda(1+2ns-2n, f\times \widetilde f\, )} \, E(z,  f; 2-s)\\
 & \hskip 45pt - \underset
 h(\gamma z) \ge A\to {\sum_{\gamma \in P_{n,n}(\Bbb Z)\backslash\Gamma}} h(z)^{n s} \left(1 - \frac{A^{\frac{n}{2}}}{h(z)^{\frac{n}{2}}}  \right) f(\frak m_1) f(\frak m_2) \; \bigg |_\gamma
 \\
 &
 \hskip 30pt
 - \frac{\Lambda(2ns-2n,\; f\times \widetilde f\, )}{\Lambda(1+2ns-2n,\; f\times \widetilde f\,)} \,\underset
 h(\gamma z) \ge A\to {\sum_{\gamma \,\in\, (P_{n,n}(\Bbb Z)\cap\Gamma)\backslash\Gamma}}  h(z)^{n(2-s)} \left(1 - \frac{A^{\frac{n}{2}}}{h(z)^{\frac{n}{2}}}  \right) f(\frak m_1) f(\frak m_2) \; \bigg |_\gamma\endalign$$
 to be the smoothed Arthur truncation of $E(z,f;s)$.

 \vskip 10pt\noindent
 {\bf Proposition 4.10} {\it Let $h(z)$ be given by (4.7). Then for any $\gamma \in P_{n,n}(\Bbb Z)\backslash SL(2n, \Bbb Z)$ and $y_i \ge 1,\; (i=1,2,\ldots ,2n-1)$, we have
 $h(\gamma z)  \le h(z).$}

 \vskip 10pt
 The proof  of proposition 4.10 relies on the following
 lemma.
 
 \vskip 10pt\noindent
 {\bf Lemma 4.11} {\it Suppose $M$ is an $m \times m$ upper triangular matrix. Let $k\le m.$ Define
 $$\align & I := (i_1, \ldots,i_k), \; \big(\text{with} \;1\le  i_1 < i_2 < \cdots < i_k \le m\big), \\
 &  J := (j_1, \ldots,j_k), \; \big(\text{with} \; 1\le  j_1 < j_2 < \cdots < j_k \le m\big),\endalign$$
 and let $[M]_{I,J}$ denote the determinant of the $k\times k$  sub matrix of $M$ that corresponds to the rows with index $I$ and the columns with index $J$. If $i_\ell \; > \; j_\ell, \; (\text{for some} \; 1 \le \ell\le k),$
 then $[M]_{I,J} = 0.$} 
 \vskip 10pt
 {\bf Proof of lemma 4.11:} Since $M = \big(m_{i,j}\big)_{\underset 1\le j\le m\to{\scriptscriptstyle 1\le i\le m}}$ is upper triangular we have 
 $m_{i,j}= 0$ if $i > j.$  Now $[M]_{I,J}$ is the minor of the matrix $M$ determined by $I,J$ and has the form
  $$[M]_{I,J} = \text{Det}\pmatrix
   & \overset \scriptscriptstyle j_1\to{\scriptstyle\downarrow} &  \overset \scriptscriptstyle j_2\to{\scriptstyle\downarrow} & \cdots & \overset \scriptscriptstyle j_\ell\to{\scriptstyle\downarrow} &  \overset \scriptscriptstyle j_r\to{\scriptstyle\downarrow} & \cdots & \overset \scriptscriptstyle j_k\to{\scriptstyle\downarrow}
   \\
\scriptscriptstyle i_1\, \rightarrow  & * & * & \cdots &* & *&  \cdots & *
\\
  & \vdots & \vdots & \cdots &\vdots & \vdots&  \cdots & \vdots
  \\
\scriptscriptstyle i_{\ell-1} \, \rightarrow  & * & * & \cdots &* & *&  \cdots & *  
\\
\scriptscriptstyle i_{\ell} \, \rightarrow  & 0 & 0 & \cdots &0 & *&  \cdots & * 
\\
  & \vdots & \vdots & \cdots &\vdots & \vdots&  \cdots & \vdots
  \\
  \scriptscriptstyle i_{k} \, \rightarrow  & 0 & 0 & \cdots &0 & *&  \cdots & * 
 \endpmatrix.$$
Laplace's theorem says that if we select any $\ell$ columns/ rows of a matrix $M$ and form all possible $\ell$-columned/rowed minors from these $\ell$ columns/rows, multiply each of these minors by its cofactor, and then add the results, we obtain $\text{Det}(M)$.
 Expand $[M]_{I,J}$ along the first $\ell$ columns. Since at least one row  in the $\ell\times \ell$ minors is zero, we immediately see that $[M]_{I,J} = 0.$
 \qed

\vskip 10pt
{\bf Proof of proposition 4.10:} Let
$$\mu_1 = \cdots = \mu_{n-1} = \mu_{n+1} = \cdots = \mu_{2n-1} = 0,\qquad \mu_n = \frac{1}{n}.$$
With this choice of $\mu = (0, \ldots, 0, \frac{1}{n},0,\ldots ,0)$ we see that (2.2) is given by
$$I_\mu(z) = \prod_{i=1}^{2n-1} y_i^{b_{i,n}\cdot \frac{1}{n}} = h(z).$$

Further, by  lemma 5.7.2 \cite{Goldfeld, 2006}, it follows that for $\gamma \in GL(2n,\Bbb Z)$ we have
$$\align  I_\mu(\gamma z) & = \prod_{i=0}^{2n-2} \big |\big | e_{2n-i} \,\gamma \, z \;\wedge \;\;\cdots \;\; \wedge \; e_{2n-1} \, \gamma \, z \; \wedge\; e_{2n} \; \gamma \, z\big | \big|^{-2n \mu_{2n-i-1}}
\\
&
\hskip 60pt
\cdot\prod_{i=0}^{2n-2} \big |\big | e_{2n-i} \, z \;\wedge \;\;\cdots \;\; \wedge \; e_{2n-1}  \, z \; \wedge\; e_{2n}  \, z\big | \big|^{2n \mu_{2n-i-1}} \;I_\mu(z)
\\
&
\\
&
= \big| \big| e_{n+1} \,\gamma\, z \; \wedge \cdots \; \wedge \; e_{2n}\,\gamma \, z\big|\big|^{-2} 
\cdot
 \big| \big| e_{n+1} \; z \; \wedge \cdots \; \wedge \; e_{2n} \, z\big|\big|^{2} \; I_\mu(z).
\endalign $$
Here $\big| \big| e_{n+1} \; z \; \wedge \cdots \; \wedge \; e_{2n} \, z\big|\big|^{2}$ is the sum of the squares of all the $n\times n$ minors of the last $n$ rows of $z$.

Fix $A,B \in GL(2n,\Bbb R),$ an integer $1\le k \le 2n$, and
$$\align & I := (i_1, \ldots,i_k), \; \big(\text{with} \;1\le  i_1 < i_2 < \cdots < i_k \le 2n\big), \\
 &  J := (j_1, \ldots,j_k), \; \big(\text{with} \; 1\le  j_1 < j_2 < \cdots < j_k \le 2n\big).\endalign$$  With these choices we have  the Cauchy-Binet formula
$$[A B]_{I,J} = \sum_K \, [A]_{I,K} \; [B]_{K,J},\tag{4.12}$$
where the sum goes over all $K = (\kappa_1, \ldots,\kappa_k), \; \big(\text{with} \;1\le  \kappa_1 <\kappa_2 < \cdots < \kappa_k \le 2n\big).$ 

\vskip 8pt
We will apply (4.12) with the choices $k = n$ and
$$B = z = \pmatrix  y_1\cdots y_{2n-1} & y_1\cdots y_{2n-2} x_{1,2} & \cdots & y_1 x_{1,2n-1} & x_{1,2n}\\
 & y_1\cdots y_{2n-2} & \cdots & y_1 x_{2,2n-1} & x_{2,2n}\\
  & & \ddots & \vdots & \vdots\\
  & &  & y_1 & x_{2n-1,2n}\\
  & & & 0 &1  \endpmatrix,$$ 
  $$A = \gamma = \pmatrix \alpha_{1,1} & \alpha_{1,2} & \cdots & \alpha_{1,2n}\\
  \vdots & \vdots & \cdots & \vdots\\
  \alpha_{n+1,1} & \alpha_{n+1,2} & \cdots & \alpha_{n+1,2n}\\
   \vdots & \vdots & \cdots & \vdots\\
   \alpha_{2n,1} & \alpha_{2n,2} & \cdots & \alpha_{2n,2n}
    \endpmatrix.$$
 
 Let $C := \gamma z.$ Then by (4.12) we see that
 $$[C]_{(n+1,n+2, \ldots, 2n), \; (1,2,\ldots,n)} = \sum_K\; [\gamma]_{(n+1,n+2, \ldots, 2n),\, K} \cdot [z]_{K, \, (1,2,\ldots,n)}.$$
 By lemma 4.11, we have $[z]_{K,\,(1,2,\ldots,n)} = 0$ unless $K = (1,2,\ldots, n).$ It follows that
 $$[C]_{(n+1,n+2, \ldots, 2n), \; (1,2,\ldots,n)} = [\gamma]_{(n+1,n+2, \ldots, 2n),\, (1,2,\ldots,n)} \cdot [z]_{(1,2,\ldots,n), \, (1,2,\ldots,n)}$$
Again, by lemma 4.11, we see that
 $$\align & [C]_{(n+1,n+2, \ldots, 2n), \; (1,2,\ldots,n-1,n+1)} =  \sum_K\; [\gamma]_{(n+1,n+2, \ldots, 2n),\, K} \cdot [z]_{K, \, (1,2,\ldots,n-1,n+1)}\\
 &\hskip 60pt = [\gamma]_{(n+1,n+2, \ldots,2n),\,(1,2,\ldots,n)}\cdot [z]_{(1,2,\ldots,n),\,(1,2,\ldots,n-1,n+1)}\\
 & \hskip 110pt+  [\gamma]_{(n+1,n+2, \ldots,2n),\,(1,2,\ldots,n-1,n+1)}\cdot [z]_{(1,2,\ldots,n),\,(1,2,\ldots n-1,,n+1)}.\endalign$$
 Continuing in the same manner it follows that
  $$\align & [C]_{(n+1,n+2, \ldots, 2n), \; (1,2,\ldots,n-2,n,n+1)} \\
  &\\
  & \hskip 30pt=  [\gamma]_{(n+1,n+2, \ldots,2n),\,(1,2,\ldots,n)}\cdot [z]_{(1,2,\ldots ,n),\,(1,2,\ldots,n-2,n,n+1)}\\
  & \hskip 57pt +\\
  &
 \hskip 67pt \ddots\\
  &
  \hskip 84pt
  + 
   [\gamma]_{(n+1,n+2, \ldots,2n),\,(1,2,\ldots ,n-1,n+1)}\cdot [z]_{(1,2,\ldots ,n-1,n+1),\,(1,2,\ldots n,n+1)}
   \\
   &
   \\
   &
  \hskip 108pt +    [\gamma]_{(n+1,n+2, \ldots,2n),\,(1,2,\ldots ,n,n+1)}\cdot [z]_{(1,2,\ldots ,n,n+1),\,(1,2,\ldots n,n+1)},
  \\
  &
  \\
  & \hskip 128pt \vdots
  \endalign$$

 $$\align & [C]_{(n+1,n+2, \ldots, 2n), \; (n+1,\ldots, 2n)} \\
  &\\
  & \hskip 30pt=  [\gamma]_{(n+1,n+2, \ldots,2n),\,(1,2, \ldots,n)}\cdot [z]_{(1,2,\ldots ,n),\,(n+1,2,\ldots,2n)}\\
  &
  \\&
  \hskip 84pt
  + 
   [\gamma]_{(n+1,n+2, \ldots,2n),\,(1,2,\ldots ,n-1,n+1)}\cdot [z]_{(1,2,\ldots n-1, n+1),\,(n+1,\ldots 2n)}
   \\
    & \hskip 105pt +\\
  &
 \hskip 115pt \ddots\\
&  \hskip 135pt +    [\gamma]_{(n+1,n+2, \ldots,2n),\,(n+1,\ldots ,2n)}\cdot [z]_{(n+1, \ldots ,2n),\,(n+1, \ldots 2n)}.
    \endalign$$
    
    For $$1\le i_1 < i_2 < \cdots < i_n \le 2n-1,$$
we have $$\align & [z]_{(i_1,i_2,\ldots,i_n),(i_1,i_2, \ldots i_n)} = y_1^{n} y_2^{n} \cdots y_{2n-i_n-1}^n\cdot y_{2n-i_n}^n\\
& \hskip 120pt\cdot y_{2n-i_n+1}^{n-1} \cdots y_{2n-i_{n-1}}^{n-1}
\\
& \hskip 120pt\cdot y_{2n-i_{n-1} +1}^{n-2} \cdots y_{2n-i_{n-2}}^{n-2}
\\
&
\hskip 160pt\vdots\\
& \hskip 120pt\cdot y_{2n-i_{2} +1} \cdots y_{2n-i_{1}}.
\endalign$$
Here, we define $y_0:= 1$ and
$$\align  [z]_{(n+1,n+2,\ldots,2n),\,(n+1,n+2,\ldots,2n)} &  = y_1^{n-1} y_2^{n-2} \cdots y_{n-2}^2 y_{n-1}
\\
& = ||e_{n+1} z \wedge \cdots \wedge e_{2n-1} z \wedge e_{2n} z||,\endalign$$
where the exponents of the powers of $y_j$ $(1\le j\le n-1)$ in the expression $$y_1^{n-1} y_2^{n-2} \cdots y_{n-2}^2 y_{n-1}\tag{4.13}$$ are less than or equal to the exponents of the powers of $y_j$ $(1\le j\le n-1)$  that occur in $[z]_{(i_1,i_2, \ldots,i_n),\,(i_1,i_2, \ldots. i_n)}$. This is because the exponent of $y_j$
in (4.13) is equal to $n-j$ while the exponent of $y_j$ in $[z]_{(i_1,i_2, \ldots,i_n),\,(i_1,i_2, \ldots. i_n)}$
is bigger than $n-j$ because $i_{n-j} \le 2n-j.$ Hence
$$ [z]_{(n+1,n+2,\ldots,2n),\,(n+1,n+2,\ldots,2n)} \;\le\;  [z]_{(i_1,i_2,\ldots,i_n),(i_1,i_2, \ldots i_n)},\tag{4.14}$$
for any $1\le i_1 < i_2 < \cdots < i_n \le 2n-1.$
\vskip 8pt 
Recall that
$$h(\gamma z) = I_\nu(\gamma z) =  ||e_{n+1}\gamma z \wedge \cdots \wedge e_{2n-1}\gamma z \wedge e_{2n}\gamma z||^{-2} \cdot || e_{n+1} z\wedge \cdots \wedge e_{2n-1} z \wedge e_{2n} z||^2 h(z).$$
To prove $h(\gamma z) \le h(z)$ it is enough to show that 
$$ ||e_{n+1}\gamma z \wedge \cdots \wedge e_{2n-1}\gamma z \wedge e_{2n}\gamma z|| \ge || e_{n+1} z\wedge \cdots \wedge e_{2n-1} z \wedge e_{2n} z||.\tag{4.15}$$
Here $$||e_{n+1} z\wedge \cdots \wedge e_{2n-1} z\wedge e_{2n} z||^2 =\left(y_1^{n-1} y_2^{n-2}\cdots y_{n-2}^2y_{n-1}\right)^2$$
and
$$\align & ||e_{n+1}\gamma z \wedge \cdots \wedge e_{2n-1}\gamma z \wedge e_{2n}\gamma z||^2 = [C]^2_{(n+1,n+2, \ldots,2n), \;(1,2,\ldots,n)}\tag{4.16}\\
&\hskip 30pt + \; [C]^2_{(n+1,n+2,\ldots,2n),\;(1,2,\ldots,n-1,n+1)} \; + \;\;\cdots\;\; +\; \; [C]^2_{(n+1,n+2,\ldots,2n), \;(n+1,\ldots,2n)}.\endalign$$
Now $$[C]_{(n+1,n+2,\ldots,2n),
 (1,2,\ldots,n)} = [\gamma]_{(n+1,n+2,\ldots,2n),
 (1,2,\ldots,n)} \cdot[z]_{(1,2,\ldots,n),
 (1,2,\ldots,n)}.$$
If $ [\gamma]_{(n+1,n+2,\ldots,2n),
 (1,2,\ldots,n)} \ne 0$ then $[\gamma]^2_{(n+1,n+2,\ldots,2n),
 (1,2,\ldots,n)} \ge 1$ since it is an integer. It immediately follows from (4.14) that $$\align [C]^2_{(n+1,n+2, \ldots,2n), \;(1,2,\ldots,n)} & \ge [z]^2_{(n+1,n+2,\ldots,2n),\,(n+1,n+2,\ldots,2n)}\\
 & =|| e_{n+1} z\wedge \cdots \wedge e_{2n-1} z \wedge e_{2n} z||^2.\endalign$$
 which proves (4.15), and, hence, also proves proposition 4.10. If $$ [\gamma]_{(n+1,n+2,\ldots,2n),
 (1,2,\ldots,n)} = 0,$$ then we use the next term, (i.e. $[C]^2_{(n+1,n+2,\ldots,2n),\;(1,2,\ldots,n-1,n+1)}$)  in equation (4.16) and repeat the previous argument. This process can be further continued and we eventually prove proposition 4.10.
\qed

\vskip 8pt
 For $s\in \Bbb C$ with $\Re(s)\gg 1$, let $\widehat{E}_{A}(z, f; s)$ denote Arthur's truncated Eisenstein series attached to the Maass cusp form $f$ defined in (4.8). Define
 $$c_s(f) := \frac{\Lambda(2ns - n, \; f\times \widetilde f\,)}{ \Lambda(1+2ns - n, \; f\times \widetilde f\,)}.\tag {4.18}$$

\vskip 10pt\noindent
{\bf Corollary 4.19} {\it Let $t\in\Bbb R$ with $t > 1.$  Let $\beta = t^{n^{10}}.$ Then for $1 \le y_i\le \left(t^{1+\epsilon}\right)^{\frac{2n(2n-1)}{2}}$ with $1\le i \le2n-1$ and $\beta \le A \le 2\beta$, we have
$$\widehat{E}_A^*(z, f; s) \; =\; \widehat{E}_{A}(z, f; s) \; - \; A^{\frac{n}{2}}\widehat{E}_{A}(z, f; s-1/2)\; + \; c_{s-\frac12}(f) \widehat{E}_{A}(z, f; 2-s).$$
}

{\bf Proof:} Recall from (4.6), (4.7), (4.8) that
$$\align
& \widehat{E}_A^*(z, f; s) = E(z,f; s) -A^{\frac{n}{2}} E(z, f; s-1/2) + \frac{\Lambda(2ns - 2n, f\times \widetilde f\,)}{\Lambda(1+2ns-2n, f\times \widetilde f\,)} \, E(z, f;  2-s)\\
& \hskip 40pt -  \underset
 h(\gamma z) \ge A\to {\sum_{\gamma \in P_{n,n}(\Bbb Z)\backslash\Gamma}} h(z)^{n s} \left(1 - \frac{A^{\frac{n}{2}}}{h(z)^{\frac{n}{2}}}  \right) f(\frak m_1)  f( \frak m_2) \; \bigg |_\gamma 
 \\
 &
 \hskip 20pt
 - \frac{\Lambda(2ns-2n,\; f\times \widetilde f\,)}{\Lambda(1+2ns-2n,\; f\times \widetilde f\,)} \,\underset
 h(\gamma z) \ge A\to {\sum_{\gamma \in P_{n,n}(\Bbb Z)\backslash\Gamma}}  h(z)^{n(2-s)} \left(1 - \frac{A^{\frac{n}{2}}}{h(z)^{\frac{n}{2}}}  \right)  f(\frak m_1) f(\frak m_2) \; \bigg |_\gamma.
\endalign$$
By the assumptions on the $y_i, (1\le i\le 2n-1)$ and proposition 4.10, we see that
$$h(\gamma z) \;\le\; h(z) = \prod_{i=1}^{2n-1} y_i^{b_{i,n}\cdot \frac{1}{n}} = y_1y_2^2y_3^3\cdots y_n^n \cdot y_{n+1}^{n-1} y_{n+2}^{n-2}  \cdots y_{2n-1} \le \left( t^{(2+\epsilon)n^2}  \right)^{n^2} < A.$$
Since $c_{s-\frac12}(f) = \frac{\Lambda(2ns - 2n, f\times \widetilde f\,)}{\Lambda(1+2ns-2n, f\times \widetilde f\,)}$, it immediately follows that
$$\ \widehat{E}_A^*(z, f;  s) = E(z, f; s) -A^{\frac{n}{2}} E(z, f; s-1/2) + c_{s-\frac12}(f) \, E(z,  f; 2-s).$$
A similar argument can be applied to $\widehat{E}_A(z, f; s),$ using (4.8).

\qed

 \vskip 20pt
 \noindent
  {\bf \S 5 Upper Bounds for Rankin-Selberg L-functions on the Line $\Re(s) = 1$:} 
 \vskip 10pt 
 
 The main result of this section is the following lemma.

  \vskip 10pt\noindent
 {\bf Lemma 5.1} {\it Let $t\in\Bbb R$ and $f$ be the Maass form on $GL(n,\Bbb Z)$ defined in (4.1). Then  for $|t| \ge 1,$ we have the bounds
\vskip-10pt  $$\align & L(1+2int, \; f\times \widetilde f\,) \; \ll_{n,f}  \log |t|, \\
 & L'(1+2int, \; f\times \widetilde f\,) \;  \ll_{n,f} (\log|t|)^2.,\endalign$$
}
 
 \vskip 5pt
 {\bf Proof of lemma 5.1:} It follows from theorem 5.3 in \cite{Iwaniec-Kowalski, 2004} that
 
 $$\align L(1+2int,\; f\times \widetilde f\,) \; &= \; \sum_{m=1}^\infty \frac{\lambda_{f\times \widetilde f\,}(m)}{m^{1+2int}} \, V_{1+2int}\left( \frac{m}{X}  \right) \, + \, \epsilon(f\times \widetilde f\,)\sum_{m=1}^\infty \frac{\overline{\lambda}_{f\times \widetilde f\,}(m)}{m^{-2int}} \, V_{-2int}(mX)\tag{5.2}\\
 & \hskip 71pt+ \left(\underset u = -2int\to{\text{Res}} \; + \; \underset u = -1-2int\to{\text{Res} } \right) \frac{\Lambda(s+u, \, f \times \widetilde f\,)}{\gamma(s, \, f\times \widetilde f\,)} \frac{G(u)}{u}\, X^u,
 \endalign$$
 where
 $X = \left(|t|+1\right)^{\frac{n^2}{2}}$. By proposition 5.4 in  \cite{Iwaniec-Kowalski, 2004}, we have
 $$V_{-2int}, \; V_{1+2int}(y) \; \ll \; \left( 1 + \frac{y}{\sqrt{q_\infty}}  \right)^{-100},$$
 with $q_\infty = \left(|t|+1\right)^{n^2} = X^2.$
 
 Consequently, $L(1+2int, \; f\times \widetilde f\,)$ is essentially approximated by the first piece in the approximate functional equation (5.2). Since
  $$\sum_{m\le N} \lambda_{f\times \widetilde f\,}(m) \ll_f N,$$
the bound for $L(1+2int, \; f\times \widetilde f,)$ follows immediately.

Next, we consider the derivative.  Let $G(u) := \left( \cos\left( \frac{\pi u}{400}  \right)  \right)^{-4n^3},$
 as in \cite{Iwaniec-Kowalski, 2004, p. 99}, and define
 $$I(X, f\times \widetilde f, \, s) := \frac{1}{2\pi i} \int\limits_{3-i\infty}^{3+i\infty} X^u \Lambda(s+u, f\times \widetilde f\,) \, G(u) \; \frac{du}{u^2},$$
 where $s = 1 + 2int.$ If we shift the line of integration to $\Re(u) = -3,$ we pick up a double pole at $u = 0$ and simple poles at $u = 1-s, \; -s.$ Note that
 $$G(0) = 1, \qquad G'(0) = 0.$$
Set $\Lambda(s, f\times \widetilde f\, ) = \gamma(s, f\times \widetilde f\,) L(s, f\times \widetilde f\,).$ It follows that
$$\align & \underset u=0\to{\text{Res}} \; X^u \Lambda(s+u, f\times \widetilde f\,) \, \frac{G(u)}{u^2}  = \Big[X^u \Lambda(s+u, f\times \widetilde f\,) \,G(u)\Big]_{u=o}'\\
& \hskip 43pt = \log X \,\Lambda(s, f\times \widetilde f\,) \; + \; \gamma'(s, f\times \widetilde f\,) L(s, f\times \widetilde f\,) \; + \; \gamma(s, f\times \widetilde f\,) L'(s, f\times \widetilde f\,).\endalign$$

Applying the functional equation (4.2), we obtain
$$\align & (\log X)\cdot \gamma(s, f\times \widetilde f\,) L(s, f\times \widetilde f\,) + \gamma'(s, f\times \widetilde f\,)  L(s, f\times \widetilde f\,)\tag{5.3}\\
&\hskip 10pt + \gamma(s, f\times \widetilde f\,) L'(s, f\times \widetilde f\,) \; + \; \left(\underset u=1-s\to{\text{Res}}  \; + \; \underset u = -s\to{\text{Res}}  \right) \Lambda(s+u, \, f\times \widetilde f\,) \frac{G(u)}{u^2} X^u\\
&\hskip 30pt = I(X, \, f\times \widetilde f, s) \; + \; \epsilon(f\times \widetilde f\,) I(X^{-1}, \, \widetilde f\times  f, \, 1-s).\endalign$$

Now, we may expand $I(X, f\times f, s)$ into an absolutely convergent series. Let
$L(s, f\times \widetilde f\,) = \sum\limits_{m=1}^\infty \lambda_{f\times \widetilde f} \,(m) m^{-s}.$ Hence
$$\align I(X, f\times \widetilde f, s) & = \sum_{m=1}^\infty \lambda_{f\times \widetilde f}\,(m) m^{-s} \cdot \frac{1}{2\pi i} \int\limits_{3-i\infty}^{3+i\infty} \gamma(s+u, f\times \widetilde f\,) \left( \frac{X}{m}  \right)^u G(u) \;\frac{du}{u^2}\\
& = \gamma(s, f\times \widetilde f\,) \sum_{m=1}^\infty \frac{\lambda_{f\times \widetilde f}\,(m)}{m^s} V_s^*\left( \frac{m}{X}  \right),\endalign$$ 

\pagebreak

where 
$$ V_s^*(y) = \frac{1}{2\pi i} \int\limits_{3-i\infty}^{3+i\infty} y^{-u} G(u) \frac{\gamma(s+u, f\times \widetilde f\,)  }{ \gamma(s, f\times \widetilde f\,) } \; \frac{du}{u^2}.\tag{5.4}$$
 
 Next, we do the same for $I\left(X^{-1}, \widetilde f\times  f, \, 1-s\right)$ and then divide both sides of   (5.3) by $\gamma(s, f\times \widetilde f\,).$ It follows that
 $$\align &\hskip 10pt (\log X)\cdot L(s, f\times \widetilde f\,) \; + \; \frac{\gamma'(s, f\times \widetilde f\,)}{\gamma(s, f\times \widetilde f\,)} L(s, f\times \widetilde f\,) \; + \; L'(s, f\times \widetilde f\,)\tag{5.5}\\
 &\hskip 170pt+ \left(\underset u=1-s\to{\text{Res}}  \; + \; \underset u = -s\to{\text{Res}}  \right) \frac{\Lambda(s+u, \, f\times \widetilde f\,)}{\gamma(s,\, f\times \widetilde f\,)  } \frac{G(u)}{u^2} X^u \\
 &\hskip 63pt = \sum_{m=1}^\infty \frac{\lambda_{f\times \widetilde f}\,(m)   }{m^{1+2int}}
 V_s^*\left( \frac{m}{X} \right) \; + \; \epsilon(f
 \times \widetilde f\,) \sum_{m=1}^\infty \frac{\overline{\lambda}_{f\times \widetilde f}\,(m)   }{m^{-2int}}
 V_s^*\left(mX \right).
 \endalign$$
 
 Recall that $q_\infty = (|t|+1)^{n^2}.$
 \vskip 8pt\noindent
 {\bf Claim:} For $0< \Re(s) \le 1$ and $|\Im(s)| \ge 1$, we have
 $$V_s^*(y) \ll  \cases \left( \frac{y}{\sqrt{q_\infty}}  \right)^{-100}, &
 \text{if}\; y \ge \sqrt{q_\infty},\\
 \log y + \Cal O\left( \frac{y}{\sqrt{q_\infty}}  \right)^{\frac{\Re(s)}{2}}, & \text{if} \; y \le \sqrt{q_\infty}.\endcases$$
   Stirling's asymptotic formula for the Gamma function
   says that for fixed $\sigma$, as $|t| \to \infty$,
$$
\left | \Gamma(\sigma + it)\right | \sim \sqrt{2\pi}\, | t| ^{\sigma-1/2}e^{-\pi |t| /2}.
$$

If we apply the above asymptotic formula to the $n^2$ Gamma functions in the numerator and denominator of $ \frac{\gamma(s+u, f\times \widetilde f\,)  }{ \gamma(s, f\times \widetilde f\,) }$ (see (4.2)) it follows that
$$\frac{\gamma(s+u, f\times \widetilde f\,)  }{ \gamma(s, f\times \widetilde f\,) } \; \ll \; q_\infty^{\frac{\Re(u)}{2}}  e^{\frac{\pi}{2} n^2 |u|}.\tag{5.6}$$

The first bound in the claim follows upon shifting the line of integration in (5.4) to $\Re(u) = 100$ and then using the estimate (5.6). If we shift the line of integration to $\Re(u) = -\Re(s)/2$, we derive the second bound.

\pagebreak

It is clear that $V_{-2it}^*(nX) \ll 1.$ The second bound in lemma 5.1 now follows from  (5.5) and the above claim.

\qed

  \vskip 20pt
 \noindent
  {\bf \S 6 Maass-Selberg Relation:} 
 \vskip 10pt

  The Maass-Selberg relation for $GL(n)$ with $n > 2$ (due to Langlands) is an identity for the inner product of two Arthur truncated Eisenstein series (see \cite{Garrett, 2005}, \cite{Labesse-Waldspurger, 2013}).
 
 \vskip 10pt\noindent
 {\bf Theorem 6.1 (Maass-Selberg relation)} {\it Let $r, s\in \Bbb C$ with $r \ne \bar s$ and $r+\bar s\ne 1.$  Then we have
 $$\align & \Big\langle \widehat{E}_{A}(*, f; \,r), 
 \; \widehat{E}_{A}(*, f;\,s)\Big \rangle \; = \; \langle f,f\rangle  \langle  f,  f\rangle\cdot \frac{ A^{n(r+\bar s-1)} }{ r+\bar s-1 } \; + \; \left| \langle f,  f\rangle  \right|^2 \cdot \overline{c}_s\, \frac{ A^{n(r-\bar s)} }{ r-\bar s} \\
 & \hskip 153pt 
  + \left| \langle f,  f\rangle  \right|^2\cdot   c_r \frac{ A^{n(\bar s - r)} }{ \bar s - r} \;+\;  \langle f,f\rangle  \langle  f, f\rangle\cdot c_r \, \overline{c}_s \frac{ A^{n(1-r-\bar s} }{1- r-\bar s)}, \endalign$$
 where $c_s := \frac{\Lambda(2ns - n, \; f\times \widetilde f\,)}{ \Lambda(1+2ns - n, \; f\times \widetilde f\,)}.$ Here $\langle f,f\rangle$ is the inner product on $SL(n,\Bbb Z)\backslash \frak h^n$, while the inner product of Eisenstein series is on $SL(2n,\Bbb Z)\backslash \frak h^{2n}.$}
 
 \vskip 8pt
 Next, we prove the following corollary.
 
 \vskip 10pt\noindent
 {\bf Corollary 6.2} {\it Define ${c_s}' :=  \frac{d}{ds} c_s.$  Then for $s = \frac12+it$ with $0\ne t\in\Bbb R$, we have
 $$\align \Big\langle \widehat{E}_{A}(*, f;\;s), 
 \; \widehat{E}_{A}(*, f;\;s)\Big \rangle &\; =  \;\left| \langle f,  f\rangle  \right|^2\cdot \overline{c}_s \frac{A^{2nit}}{2it} \; + \;  \left| \langle f,  f\rangle  \right|^2\cdot c_s\frac{A^{-2nit}}{-2it} \;
 \\
 &   \hskip 80pt - \;\langle f,f\rangle  \langle  f,  f\rangle\cdot \frac{c_s'}{c_s} \; +\;\langle f,f\rangle  \langle  f,   f\rangle\cdot 2n\log A.\endalign$$}

 \vskip 8pt
 {\bf Proof of corollary 6.2:}
 Let
 $$s = \frac12+it, \qquad r = s + \epsilon, \qquad (t\in\Bbb R, \;\, \epsilon >0),$$
 for $\epsilon \to 0.$
Then Theorem 6.1 takes the form
$$\align & \Big\langle \widehat{E}_{A}(*, f;\;r), 
 \; \widehat{E}_{A}(*, f;\;s)\Big \rangle \; = \; \langle f,f\rangle  \langle  f,  f\rangle\cdot \frac{ A^{n\epsilon} }{ \epsilon } \; + \; \left| \langle f,  f\rangle  \right|^2\cdot \overline{c}_s\, \frac{ A^{2nit+\epsilon} }{ 2it+\epsilon} \;\\
 &\hskip 160pt + \; \left| \langle f,  f\rangle  \right|^2\cdot  c_r\, \frac{ A^{-n\epsilon-2nit} }{ -\epsilon-2it} \,-\,  \langle f,f\rangle  \langle  f,  f\rangle\cdot c_r \, \overline{c}_s\, \frac{ A^{-n\epsilon} }{\epsilon}.\tag{6.3} \endalign$$
 
 On the line $\Re(s) = \frac12$, we have $|c_s| = 1.$ It follows by the Taylor expansion around $s=0$ that
 $$\align  c_r \cdot \overline{c}_s & = \overline{c}_s \left[ c_s +   {c_s}' \epsilon + o(\epsilon^2)\right]\\
 & = 1 + \overline{c}_s {c_s}' \epsilon + o(\epsilon^2).\endalign$$
 Further,
 $$A^{n\epsilon} = 1 + (\log A)\cdot n\epsilon + o(\epsilon^2).$$
 $$A^{-n\epsilon} = 1 - (\log A)\cdot n\epsilon + o(\epsilon^2).$$
 Corollary 6.2 immediately follows upon letting $\epsilon \to 0$ in (6.3).\qed
 
 \vskip 10pt\noindent
 {\bf Proposition 6.4} Let $\Cal F := SL(2n,\Bbb Z)\big\backslash GL(2n,\Bbb R)\big/(O(n,\Bbb R)\cdot \Bbb R^\times).$ Then we have
 $$\align & \big |L(1+2int,\, f\times \widetilde f\,)\big| \underset 1\le i \le 2n-1\to{ \underset 1\le y_i \le (t^{1+\epsilon})^{n(2n-1)}\to{ \int_{\Cal F}}} \left| \widehat{E}_A^*(z, f;\, 1+it)   \right|^2 \; d^\star z \;\, \ll \;\, A^n (\log A)(\log |t|).\endalign$$

 {\bf Proof:} This follows immediately from corollary 4.19, the Maass-Selberg relation in theorem 6.1 and corollary 6.2, and the upper bounds for Rankin-Selberg L-functions given in lemma 5.1.
 
 \qed

 \vskip 20pt
 \noindent
  {\bf \S 7 Stade's formula for Whittaker functions:} 
 \vskip 10pt
 Let $\nu = (\nu_1, \ldots,\nu_{n-1}) \in \Bbb C^{n-1},$ denote spectral parameters for $GL(n, \Bbb R)$ with $n\ge 2.$  The Whittaker function $W_{ n,\nu}$ on 
  $GL(n,\Bbb R)/(O(n,\Bbb R)\cdot\Bbb R^\times)$ is defined to be
  
  $$W_{n,\nu}(z) :=  \int\limits_{\Bbb R^{\frac{(n-1)n}{2}}} I_\nu(w_nuz) e^{-2\pi i(u_{1,2} +\cdots + u_{n-1,n})} \hskip-5pt\prod_{1 \le i < j\le n} du_{i,j},
  \tag{7.1}$$
    where $w_n = \left(\smallmatrix & &\hskip-2pt 1\\
  & \smallmatrix & &\hskip 2pt.\\
  & . &\\
  \hskip-4pt.&&\endsmallmatrix &\\
  1 & &\endsmallmatrix\right)$ is the long element of the Weyl group of $GL(n, \Bbb R)$ and 
  $$\nu_{j,k} = \sum_{i=0}^{j-1} \frac{n\nu_{n-k+i} - 1}{2}.$$
  The completed Whittaker function is defined to be
  $$W_{n,\nu}^*(z) := W_{n,\nu}(z) \prod_{j=1}^{n-1} \prod_{k=j}^{n-1} \frac{\Gamma\left( \frac12+\nu_{j,k}  \right)}{\pi^{\frac12+\nu_{j,k}}}.$$

 Stade found a beautiful formula which expresses the Whittaker function  
 as an integral of $K$-Bessel functions. For $\nu\in\Bbb C$, the Bessel function $K_\nu$ is defined as
 $$K_\nu(2\pi y) := \frac12\int\limits_0^\infty e^{-\pi y(t+1/t)} \, t^\nu \;\frac{dt}{t}.\tag{7.2}$$
 For $y,a,b > 0,$ define 
 $$K^*_\nu(y; a, b) := 2\left(\frac{a}{b}\right)^\frac{\nu}{2} \cdot K_\nu\left(2\pi y\sqrt{ab}\right).$$
 For any array $\{u_{k,\ell}\}$, define
 $$S(j,i) := 1 + \frac{1}{u_{j-1,i}}\left( 1 + \frac{u_{j-2,i-1}}{u_{j-2,i}}   \left(  1 + \frac{u_{j-3,i-1}}{u_{j-3,i}}\left( 1 + \cdots +\frac{u_{2,i-1}}{u_{2,i}} \left(1 + \frac{u_{1,i-1}}{u_{1,i}}   \right)      \cdots\right)  \right)\right)$$
 for $1 \le j\le i.$
 
 \vskip 10pt\noindent
 {\bf Proposition 7.3 (Stade, 1990)} {\it Let $n \ge 2$ and $\nu = (\nu_1, \ldots,\nu_{n-1})\in\Bbb C^{n-1}$ with $\Re(\nu_i) >\frac{1}{n}$ for $1 \le i \le n-1.$ Set
 $\mu_j := \sum\limits_{k=j}^{n-1} \nu_{j.k}.$  Then
 $$\align & W_{n,\nu}^*(y) =  \prod_{i=1}^{n-1}\; y_{n-i}^{\sum\limits_{j=1}^{n-1}
 b_{i,j}\nu_j -\mu_i     } \hskip-12pt\int\limits_{\left(\Bbb R^+\right)^{\frac{(n-2)(n-1)}{2}}} \prod_{i=1}^{n-1} \Bigg\{K^*_{\mu_i}\left( y_{i}\, ; \, 1+ \sum_{k=i+1}^{n-1} u_{1,k}, \, S(i,i)  \right)\\
 & \hskip 187pt \cdot\left( \prod_{k=i+1}^{n-1} u_{i,k}^{-\nu_{i,k-i+1}} \right) \Bigg\} \; \prod\limits_{1\le i < j \le n-1} \frac{du_{i,j}}{u_{i,j}}.\endalign$$
 }

 \vskip 10pt\noindent
 {\bf Corollary 7.4} {\it Let $n \ge 2,$ and $\nu = (\nu_1, \ldots,\nu_{n-1})\in\Bbb C^{n-1}$ with associated Langlands parameters $\alpha = (\alpha_1, \ldots,\alpha_n) \in\Bbb C^n$ as in (2.3). Assume $\Im(\alpha_i) \le C |t|$ for $t\in \Bbb R$ with $|t|\to\infty$ and $C> 0$ a fixed constant. Further   assume that
 $$\max_{1\le i \le n-1} y_i \; \ge \;\, |t|^{1+\epsilon},  \qquad\quad (\text{for some fixed} \; \epsilon > 0).$$
 Then 
 $$W_{n,\nu}(y) \ll  \Big(\max_{1\le i\le n-1} y_i\Big)^{-N}$$
 for any fixed large $N> 1.$ The implied constant in the $\ll$-symbol above does not depend on $t, y$.}
 
\pagebreak
 
\vskip 10pt
{\bf Proof:}  The proof follows from proposition 7.3 and the bound for the $K$-Bessel function
 $$e^{\frac{\pi}{2} |\Im \nu|}K_\nu(y) \ll y^{-N}, \qquad \; \Big(\text{for} \;\, y \ge (|\Im(\nu) |+1)^{1+\epsilon}, \;\; |\Re(\nu)| \ll 1\Big),$$
 which holds for any fixed $N > 1.$ \qed

\vskip 6pt
Assuming $\Re(\alpha_i) = 0$ for $1 \le i\le n$,  \cite{Brumley-Templier}
 proved corollary 7.4 using a different method.

 \vskip 20pt
 \noindent
  {\bf \S 8 Fourier Expansions of Eisenstein Series:}

 \vskip 10pt\noindent
 {\bf Definition 8.1 (Automorphic function)} {\it Let $k\in\Bbb Z$ with $k\ge 2.$ Let $z = xy \in GL(k,\Bbb R)$ as in (2.1).  A function $$F:SL(k, \Bbb Z)\backslash GL(k,\Bbb R)\big/O(k, \Bbb R)\cdot \Bbb R^\times \; \to \;\Bbb C$$
 is called an automorphic function if it is smooth and $F(z)$  has polynomial growth in $$Y_{_{k_1,\ldots,k_r}} =\prod\limits_{i=1}^r y_{_{k_i}}, \qquad (Y_{_{k_1,\ldots,k_r}}\to\infty),$$   for all $1 \le k_1 < k_2, \cdots < k_r\le n-1$, and where $y_j$ is fixed if $j \ne k_i\, (1\le i\le r).$ }
 
 \vskip 10pt\noindent
 {\bf Definition 8.2} Let $m,n\in\Bbb Z$ with $n\ge2$ and $1\le m <n.$ We define the parabolic subgroup 
 $$\widetilde P_{n,m} := \left( \pmatrix  & * & * & \cdots & *\\
 SL(n-m)& \vdots &\vdots &\cdots&\vdots\\
 & * & * &\cdots & *\\
 & 1 & * &\cdots & *\\
 & &\ddots & &\vdots\\
 & & &1&*\\
 & & & & 1\endpmatrix  \right) \; \subset SL(n).$$
 
 \pagebreak
 
 \vskip 10pt\noindent
 {\bf Theorem 8.3} {\it Suppose $F$ is an automorphic function for $SL(k, \Bbb Z)$ with $k \ge 2.$ Then $F$ has the following Fourier expansion.
 $$\align & F(z) = \int\limits_0^1\cdots \int\limits_0^1 F\left(\left( \smallmatrix 1 & 0 &\cdots & 0 & u_{1,k}\\
  & 1 & \cdots & 0 &\vdots\\
  & & \ddots &\vdots &\vdots\\
  & & &  1 & u_{k-1,k}\\
  &&&&1\endsmallmatrix \right) z \right) \; du_{1,k}\, du_{2,k} \cdots du_{k-1,k}\\
  &
  \\
  &
  \\
  & \hskip 20pt + \sum_{1\le\ell\le k-2} \;\sum_{m_{k-1}\ne 0}\cdots \sum_{m_{k-\ell}\ne0} \;\, \sum_{\gamma_{k-1} \in \widetilde P_{k-1,\ell}(\Bbb Z)\backslash SL(k-1, \Bbb Z)}  \\
  & \hskip 5pt \cdot \int\limits_0^1\cdots \int\limits_0^1
  F\left(\left( \smallmatrix 1&0 &\cdots &0 & u_{1,k-\ell} & \cdots & u_{1,k-1} & u_{1,k}\\
  &1 &\cdots &0 & u_{2,k-\ell} & \cdots & u_{2,k-1} & u_{2,k}\\
  & & \hskip 16pt\ddots& &  &  &\vdots  &\vdots \\
   & & \hskip 10pt& &  &  &\vdots  &\vdots \\
  & & & & & & 1 & u_{k-1,k}  \\
  & & & & & &  & 1 
     \endsmallmatrix\right) \cdot \pmatrix  \gamma_{k-1} &\\
     & 1\endpmatrix z  \right)\\
     &\hskip 60pt\cdot e^{-2\pi i\big(m_{k-1} u_{k-1,k} \; + \;\cdots \; +\; m_{k-\ell} u_{k-\ell,\;k-\ell+1}\big)} \; \prod_{i=1}^{k-1} \;\prod_{  j = \max(k-\ell, \, i+1)         }^k du_{i,j}\\
     &
     \\
     &
     \\
     &\hskip 20pt + \sum_{m_{k-1}\ne0} \cdots \sum_{m_1\ne 0} \; \sum_{\gamma_{k-1} \in U_{k-1}(\Bbb Z)\backslash SL(k-1,\Bbb Z)} \\
     &\hskip 60pt \cdot\int\limits_0^1\cdots \int\limits_0^1   F\left(\left( \smallmatrix 1 & u_{1,2} & \cdots & u_{1,k-1} & u_{1,k}\\
   & 1 & \cdots & u_{2,k-1} & u_{2,k}\\ 
    &  & \ddots & \vdots & \vdots\\
     &  &  &1 & u_{k-1,k}\\ 
      &  &  & & 1\endsmallmatrix\right) \cdot \pmatrix  \gamma_{k-1} &\\
     & 1\endpmatrix z \right)\\
     &\hskip 98pt \cdot e^{-2\pi i\big(m_{k-1} u_{k-1,k} \; + \; \cdots \; + \; m_1 u_{1,2}\big)    } \; \prod_{i=1}^{k-1}\prod_{j=i+1}^k du_{i,j}. \endalign$$  }
 \indent     {\bf Proof:} See \cite{Goldfeld, 2006}, \cite{Ichino, Yamana}.

  \vskip 10pt
  To simplify the later exposition, we define the last sum above as the non degenerate term associated to the automorphic function $F$.   Here is a formal definition of both the degenerate and the non degenerate terms..
 
\pagebreak
   
  \vskip 10pt\noindent
  {\bf Definition 8.4 (degenerate and non degenerate term)} Let $F$ be an automorphic function for $SL(k, \Bbb Z)$ with $k \ge 2.$ We define the non degenerate term associated to $F$ to be
  $$\align & \text{\rm ND}(F) := \hskip 20pt  \sum_{m_{k-1}\ne0} \cdots \sum_{m_1\ne 0} \; \sum_{\gamma \in U_{k-1}(\Bbb Z)\backslash SL(k-1,\Bbb Z)} \\
     &\hskip 90pt \cdot\int\limits_0^1\cdots \int\limits_0^1   F\left( \pmatrix 1 & u_{1,2} & \cdots & u_{1,k-1} & u_{1,k}\\
   & 1 & \cdots & u_{2,k-1} & u_{2,k}\\ 
    &  & \ddots & \vdots & \vdots\\
     &  &  &1 & u_{k-1,k}\\ 
      &  &  & & 1\endpmatrix \cdot \pmatrix  \gamma_{k-1} &\\
     & 1\endpmatrix z \right)\\
     &\hskip 180pt \cdot e^{-2\pi i\big(m_{k-1} u_{k-1,k} \; + \; \cdots \; + \; m_1 u_{1,2}\big)    } \; \prod_{i=1}^{k-1}\prod_{j=i+1}^k du_{i,j}. \endalign$$
     The other terms in the expansion for $F(z)$ given in  theorem 8.3 are defined to be the degenerate term which is denoted $D(F)$.

     \vskip 10pt \noindent
 {\bf Proposition 8.5 (Langlands)} {\it The Eisenstein series $E(z, f; s)$ has constant term $C_{\Cal P}$ (as defined in (4.6)) along the parabolic $\Cal P = P_{n,n}$. Along all other parabolics, the constant term is 0.}
 
 \vskip 10pt
 {\bf Proof:} See theorems 6.2.1, 6.3.1 in \cite{Shahidi, 2010} and proposition 10.10.3 in \cite{Goldfeld, 2006}.  
 \qed
 
\vskip 10pt 
 {\bf Proposition 8.6} {\it Let $z \in GL(2n, \Bbb R)$ with Iwasawa decomposition
 $$z = \pmatrix 1& x_{1,2} & \cdots & x_{1,2n}\\
  & 1 & \cdots & x_{2,2n}\\
   &  & \ddots & \vdots\\
    &  &  & 1\\
   \endpmatrix \cdot \pmatrix Y_1 &\\ & Y_2\endpmatrix ,$$
   where $Y_1, Y_2$ are diagonal matrices in $GL(n,\Bbb R).$
  Let $E(z, f; s)$ be the Eisenstein series defined in (4.5). Then we have
 $$\align & E(z, f; s) = \sum_{\gamma \in \widetilde{P}_{2n-1,\,n-1}(\Bbb Z)\backslash SL(2n-1, \Bbb Z)}
 \Bigg[  \left|\frac{\text{\rm Det}(\frak m_1)}{\text{\rm Det}(\frak m_2)}\right|^{ns} f(\frak m_1)  f^*(\frak m_2)\;\\
 & \hskip 118pt + \; c_s(f)\cdot\left|\frac{\text{\rm Det}(\frak m_1)}{\text{\rm Det}(\frak m_2)}\right|^{n(1-s)}   f(\frak m_1) f^*(\frak m_2)\Bigg] \cdot \Bigg|_{\left( \smallmatrix \gamma &\\ &1\endsmallmatrix\right)} \; + \; \text{\rm ND}(E).\endalign$$
 Here
 $$\align & f^*(\frak m_2) = \sum_{m_{2n-1}\ne 0} \cdots \sum_{m_{n+1}\ne 0}  \frac{ B(m_{2n-1}, \ldots, m_{n+1}) }{\prod\limits_{k=1}^{n-1} |m_{_{2n-k}}|^{^{\frac{k(n-k)}{2}}}  } \cdot W_{n,\nu}\left(MY_2 \right)\\
 &
 \hskip 180pt \cdot e^{2\pi i\big(m_{2n-1} x_{2n-1, 2n} + \quad \cdots \quad + m_{n+1} x_{n+1, n+2}\big)},
 \endalign$$
  where $B$  denotes the Fourier coefficient of $f$,  the matrix $M$ is the diagonal matrix $$M=\text{\rm diag}\big((m_{2n-1}m_{2n-2}\cdots |m_{n+1}|, \;\;\ldots, \;\;m_{2n-1}m_{2n-2}, m_{2n-1},1)\big),$$
 and $W_{n,\nu}$, is the Whittaker function  (see (7.1))  in the Fourier expansion of $f.$ 
 }
 
 \vskip 10pt
{\bf Proof of Proposition 8.6:} There are two cases to consider.

\vskip 8pt
$\underline{\text{\bf CASE 1, $\ell \ne n-1$:}}$
 \hskip 4pt  Let
   $$U_\ell := \pmatrix 1 & u_{2n-\ell, 2n-\ell+1} & \cdots& u_{2n-\ell,\,2n-1} & u_{2n-\ell,\, 2n}\\
  &1 & \cdots & u_{2n-\ell+1,\, 2n-1} & u_{2n-\ell+1,\, 2n}\\
  & &\ddots&\vdots&\vdots\\
  & &  & 1 & u_{2n-1,\, 2n}\\
  & &  &   & 1
  \endpmatrix.$$
 Then we see that  for $L = 2n-\ell-1,$ we have
 
 $$\align & \pmatrix u_{1,\, 2n-\ell} & \cdots & u_{1,\,2n}\\
 u_{2,\, 2n-\ell} & \cdots & u_{2,\,2n}\\
 \vdots & \cdots & \vdots\\
  u_{2n-\ell-1,\, 2n-\ell} & \cdots & u_{2n-\ell-1,\,2n}
  \endpmatrix \cdot U_\ell
  \\
  &
  \\
  & =\hskip-3pt \pmatrix u_{1,\,2n-\ell} & u_{1,\,2n-\ell}\; u_{2n-\ell, \, 2n-\ell+1} + u_{1,\,2n-\ell+1} & \cdots & u_{1,\,2n-\ell}\; u_{2n-\ell,2n} + \cdots + u_{1,2n}\\
  u_{2,\,2n-\ell} & u_{2,\,2n-\ell}\; u_{2n-\ell, \, 2n-\ell+1} + u_{2,\,2n-\ell+1} & \cdots & u_{2,\,2n-\ell}\; u_{2n-\ell,2n} + \cdots + u_{2,2n}\\
 \vdots & \vdots & \cdots & \vdots\\
 u_{L,2n-\ell} & u_{L,2n-\ell}\; u_{2n-\ell, \, 2n-\ell+1} + u_{L,\,2n-\ell+1} & \cdots & u_{L, 2n-\ell}\; u_{2n-\ell,2n} + \cdots + u_{L,2n}\\
  \endpmatrix.\endalign$$
 \vskip 5pt
 
 In the above, we make the change of variables:
 $$\align & u_{1,2n-\ell} \; u_{2n-\ell,  2n-\ell+1} \; + \; u_{1,2n-\ell+1} \; \to \; u_{1, 2n-\ell+1} \;\bigg | \;
 u_{1,2n-\ell}\; u_{2n-\ell,2n} + \cdots + u_{1,2n} \; \to \; u_{1,2n}\\
 & \hskip 105.2pt \vdots \hskip 105.2pt \bigg | \hskip 90pt \vdots\\
 & u_{L,2n-\ell} \, u_{2n-\ell,  2n-\ell+1}  +  u_{L,2n-\ell+1} \; \to \; u_{L, 2n-\ell+1} \;\bigg | \;
 u_{L,2n-\ell}\, u_{2n-\ell,2n} + \cdots + u_{L,2n} \; \to \; u_{L,2n}
  \endalign$$
  
  \pagebreak
 
 It follows from the above calculations that with the notation $E(z) := E(z,f; s)$, we have
   $$\hskip -10pt\overset L=2n-\ell-1\to{\overbrace{\phantom{xxxxxxxxx}}}  \hskip 18pt \overset \ell+1\to{\overbrace{\phantom{xxxxxxxxxxxxxxxx}}} $$
  \vskip -20pt
 $$\int\limits_0^1 \cdots \int\limits_0^1 E
\left(\pmatrix 
 1 & 0 & \cdots & 0 & u_{1,2n-\ell} & \cdots & u_{1,2n}\\
  & 1 & \cdots & 0 & u_{2,2n-\ell} & \cdots & u_{2,2n}\\
   & & \ddots & \vdots & \vdots & \cdots & \vdots\\
    &  & & 1 & u_{L,2n-\ell} & \cdots & u_{L,2n}\\
    &  & &  & \hskip-2pt1 & \cdots & u_{L+1,2n}\\
     &  & &  & \hskip 30pt\ddots & \vdots & \vdots\\
      &  & &  &  &1  & u_{2n-1,2n}\\
     &  & &  &  &  & 1
  \endpmatrix  \pmatrix \gamma &\\
  & 1\endpmatrix \cdot z  \right)$$
  $$\hskip 115pt\cdot {\displaystyle e}^{\displaystyle -2\pi i(m_{2n-1} u_{2n-1,2n} + \cdots + m_{2n-\ell} u_{2n-\ell, 2n-\ell+1})} \; d^\times u$$

\noindent
is equal to (for $ I_{L\times L}$, the $L\times L$ identity matrix)
    
   $$\int\limits_0^1 \cdots \int\limits_0^1 E\left(\pmatrix 
 1 & 0 & \cdots & 0 & u_{1,2n-\ell} & \cdots & u_{1,2n}\\
  & 1 & \cdots & 0 & u_{2,2n-\ell} & \cdots & u_{2,2n}\\
   & & \ddots & \vdots & \vdots & \cdots & \vdots\\
    &  & & 1 & u_{L,2n-\ell} & \cdots & u_{L,2n}\\
    &  & &  & 1 & \cdots & 0\\
     &  & &  &  & \ddots & \vdots\\
     &  & &  &  &  & 1
  \endpmatrix  \pmatrix I_{L\times L} & 0\\
  &
  \\
  0 & U_\ell\endpmatrix   \pmatrix \gamma &\\& 1\endpmatrix\cdot z \right)$$
 $$\hskip 115pt\cdot \, {\displaystyle e}^{\displaystyle -2\pi i(m_{2n-1} u_{2n-1,2n} + \cdots + m_{2n-\ell} u_{2n-\ell, 2n-\ell+1})} \; d^\times u,$$
 which is $0$ by proposition 8.5. 
 \vskip 15pt
$\underline{\text{\bf CASE 2, $\ell = n-1$:}}$
\hskip 4pt In this case, define
 $$\Cal K_n := \int\limits_0^1 \cdots \int\limits_0^1 E\left(  \pmatrix  I_{n\times n}   & \pmatrix
 u_{1,n+1} & \cdots & u_{1,2n}\\
\vdots &  \cdots&\vdots\\
u_{n,n+1} & \cdots & u_{n,2n}
  \endpmatrix\\
  &
  \\
  &
  \\
 & I_{n\times n}
 \endpmatrix 
 \pmatrix I_{n\times n} & 0\\
 &
 \\
  0 & U_{n-1}\endpmatrix   \pmatrix \gamma &\\& 1\endpmatrix\cdot z \right) $$
  $$\hskip 155pt\cdot \, {\displaystyle e}^{\displaystyle -2\pi i(m_{2n-1} u_{2n-1,2n} + \cdots + m_{n+1} u_{n+1,n+2})} \; d^\times u.$$
 
 Let $z' = \pmatrix \gamma &\\&1\endpmatrix \cdot z$. We write  $z'$ in Iwasawa form, so that
  $$z' = \pmatrix X_1' & X_2'\\ & X_3'\endpmatrix \pmatrix Y_1' & \\
 & Y_2'\endpmatrix \subset GL(2n, \Bbb R).$$
 Then
 $$\pmatrix I_{n\times n} &\\
 & U_{n-1}\endpmatrix \cdot z'\; = \;  \pmatrix X_1' & X_2'\\ 
 &
 \\
 & U_{n-1} X_3'\endpmatrix \pmatrix Y_1'& \\
 & Y_2'\endpmatrix.$$
Here $X_3'$ is a unipotent matrix with off diagonal entries $x_{i,j}'$ and $U_{n-1} X_3'$ is equal to
 
 $$\align &   \pmatrix 1 & u_{_{n+1, n+2}} & \cdots& u_{_{n+1,2n-1}} & u_{_{n+1, 2n}}\\
  & 1 & \cdots & u_{_{n+2, 2n-1} }  & u_{_{n+2, 2n}}\\
 & &\ddots & \vdots &\vdots\\
 & & &1 & u_{_{2n-1,2n}}\\
 & & & & 1
 \endpmatrix
 \pmatrix 1 & x_{_{n+1, n+2}}' & \cdots& x_{_{n+1,2n-1}}' & x_{_{n+1, 2n}}'\\
  & 1 & \cdots & x_{_{n+2, 2n-1}}'   & x_{_{n+2, 2n}}'\\
 & &\ddots & \vdots &\vdots\\
 & & &1 & x_{ 2n-1,2n}'\\
 & & & & 1
 \endpmatrix\\
 &\hskip 50pt = \pmatrix u_{n+1,n+2} + x_{n+1,n+2}' &  &\hskip -20pt\cdots & u_{n+1,2n} + \cdots + x_{n+1,2n}'\\
 & 1 & \cdots&u_{n+2,2n} + \cdots + x_{n+2,2n}'\\
 & &\hskip 20pt\ddots &\vdots\\
 & & & 1
 \endpmatrix, \endalign
 $$
 which becomes
 $$ \pmatrix 1 & u_{_{n+1, n+2}} & \cdots& u_{_{n+1,2n-1}} & u_{_{n+1, 2n}}\\
  & 1 & \cdots & u_{_{n+2, 2n-1} }  & u_{_{n+2, 2n}}\\
 & &\ddots & \vdots &\vdots\\
 & & &1 & u_{_{2n-1,2n}}\\
 & & & & 1
 \endpmatrix$$ after making the following  change of variables:
 $$u_{n+1,n+2} + x_{n+1,n+2}' \;\to\; u_{n+1,n+2}, \qquad  \ldots  \qquad ,u_{2n-1,2n} + x_{2n-1,2n}' \;\to\; u_{2n-1,2n}.$$

 \pagebreak
 
 Consequently, 
  $$\align & \Cal K_n = \int\limits_0^1 \cdots  \int\limits_0^1 \left[ \int\limits_0^1 \cdots \int\limits_0^1 E\hskip-5pt\left(  \pmatrix  I_{n\times n}   & \hskip-5pt\smallmatrix
 u_{1,n+1} & \cdots & u_{1,2n}\\
\vdots &  \cdots&\vdots\\
u_{n,n+1} & \cdots & u_{n,2n}
  \endsmallmatrix\\
  &
  \\
  &
  \\
 & I_{n\times n} 
 \endpmatrix \hskip -3pt
 \pmatrix X_1'Y_1' & X_2' Y_2'\\
 &  \\
 &
 \\
  0 & U_{n-1} Y_2'\endpmatrix  \right)\right. \left.  \hskip-3pt\prod_{i=1}^{n} \prod_{j=n+1}^{2n} du_{i,j}
  \vphantom{E\left(  \pmatrix  I_{n\times n}   & \hskip-5pt\pmatrix
 u_{1,n+1} & \cdots & u_{1,2n}\\
\vdots &  \cdots&\vdots\\
u_{n,n+1} & \cdots & u_{n,2n}
  \endpmatrix\\
  &
  \\
  &
  \\
 & I_{n\times n} 
 \endpmatrix \hskip -3pt
 \pmatrix X_1'Y_1' & X_2' Y_2'\\
 &  \\
 &
 \\
  0 & U_{n-1} Y_2'\endpmatrix  \right)}
  \right]\endalign$$
  $$\align &
  \hskip 100pt\cdot \, {\displaystyle e}^{\displaystyle -2\pi i(m_{2n-1} u_{2n-1,2n} + \cdots + m_{n+1} u_{n+1,n+2})} \; \prod_{k=n+1}^{2n-1} \prod_{\ell=k+1}^{2n} du_{k,\ell}\\
  & 
  \hskip 182pt
  \cdot 
  {\displaystyle e}^{\displaystyle 2\pi i(m_{2n-1} x_{2n-1,2n}' + \cdots + m_{n+1} x_{n+1,n+2}')}.  \endalign$$
By proposition 8.5, the inner integral above is
$$  
  \left| \frac{\text{det} (Y_1') }{\text{det} (Y_2')  } \right|^{ns}  f\Big(X_1' Y_1'\Big)  f\Big(U_{n-1} Y_2'\Big) \; + \;  
  c_s(f)\cdot \left| \frac{\text{det} (Y_1') }{\text{det} (Y_2')   } \right|^{n(1-s)}  f\Big(X_1' Y_1'\Big) f\Big(U_{n-1} Y_2'\Big).$$
  
  Further, the outer integral picks up the Fourier coefficient of $f$ and $\widetilde f$, respectively:
  
 $$\align & \int\limits_0^1\cdots \int\limits_0^1 f\Big(U_{n-1} Y_2'  \Big) e^{-2\pi i \big( m_{2n-1} u_{2n-1, 2n} + \;\;\cdots \;\; m_{n+1} u_{n+1, n+2}  \big)}\; \prod_{k=n+1}^{2n-1} \prod_{\ell=k+1}^{2n} du_{k,\ell}
 \\
 &
 \hskip 50pt
 = \frac{B(m_{2n-1}, \ldots, m_{n+1})}{\prod\limits_{k=1}^{n-1} |m_{_{2n-k} }|^{\frac{k(n+k)}{2}}   } \; W_{n,\nu}\left(MY_2'   \right)
 \endalign $$
where $B(m_{2n-1}, \ldots, m_{n+1})$ is the Fourier coefficient of $f$.

\qed

\vskip 10pt 
 {\bf Corollary 8.7} {\it Let $z \in GL(2n, \Bbb R)$  as in proposition 8.6.
  Let $\widehat E_A^*(z, f; s)$ be the smoothed truncated Eisenstein series in definition 4.9, and set
  $$\widetilde \Gamma(2n-1) := \widetilde P_{2n-1,n-1}(\Bbb Z)\big\backslash SL(2n-1,\Bbb Z).$$ 
  Then we have 
  $$\align & \widehat E_A^*(z, f; s) = \underset
 {h\left(\left(\smallmatrix \gamma &\\ & 1\endsmallmatrix\right)  z\right) \;< \;A}\to {\sum_{\gamma \in \widetilde \Gamma(2n-1)}}
 \left| \frac{\text{\rm det}( \frak m_1)}{\text{\rm det}( \frak m_2)}  \right|^{ns} \left[ 1 - \frac{A^{n/2}}{h(y)^{n/2}}   \right] f(\frak m_1)  f^*(\frak m_2) \; \bigg | _{ \left(\smallmatrix \gamma &\\ & 1\endsmallmatrix\right)}
 \\
 &
 \\
 &
+ \frac{\Lambda(2ns-2n,\; f\times \widetilde f\,)}{\Lambda(1+2ns-2n,\; f\times \widetilde f\,)}\hskip-5pt \underset
 {h\left(\left(\smallmatrix \gamma &\\ & 1\endsmallmatrix\right)  z\right) \;< \;A}\to {\sum_{\gamma \in \widetilde \Gamma(2n-1)}}\hskip-4pt \left| \frac{\text{\rm det}( \frak m_1)}{\text{\rm det}( \frak m_2)}  \right|^{n(2-s)}  \left[ 1 - \frac{A^{n/2}}{h(y)^{n/2}}   \right]
 f(\frak m_1) f^*(\frak m_2) \; \bigg | _{ \left(\smallmatrix \gamma &\\ & 1\endsmallmatrix\right)}
 \\
 &
 \\
 &
\hskip 40pt+  {\sum_{\gamma \in \widetilde \Gamma(2n-1)}}  \left| \frac{\text{\rm det}( \frak m_1)}{\text{\rm det}( \frak m_2)}  \right|^{n(1-s)} c_s(f) \cdot  f(\frak m_1) f^*(\frak m_2) \; \bigg | _{ \left(\smallmatrix \gamma &\\ & 1\endsmallmatrix\right)}
 \\
 &
 \\
 &
 \hskip 30pt +  \frac{\Lambda(2ns-2n,\; f\times \widetilde f\,)}{\Lambda(1+2ns-2n,\; f\times \widetilde f\,)}
   \sum_{\gamma \in \widetilde \Gamma(2n-1)}  \left| \frac{\text{\rm det}( \frak m_1)}{\text{\rm det}( \frak m_2)}  \right|^{n(s-1)} c_{2-s}(f) \cdot  f(\frak m_1)  f^*(\frak m_2) \;  \bigg | _{ \left(\smallmatrix \gamma &\\ & 1\endsmallmatrix\right)}
 \\
 &
 \\
 &
 \\
 &
  \hskip 20pt -  \sum_{\gamma \in \widetilde \Gamma(2n-1)} \frac{1}{2\pi i} \int\limits_{2-i\infty}^{2+i\infty}
  \frac{A^{-\frac{n}{2}w}}{w(w+1)} \left| \frac{\text{\rm det}( \frak m_1)}{\text{\rm det}( \frak m_2)}  \right|^{n\left(1-s-\frac{w}{2}\right)} c_{s+\frac{w}{2}}(f) \cdot  f(\frak m_1) f^*(\frak m_2) \; \bigg | _{ \left(\smallmatrix \gamma &\\ & 1\endsmallmatrix\right)}
 \; dw
  \\
  &
  \\
  &
  \hskip 30pt -   \frac{\Lambda(2ns-2n,\; f\times \widetilde f\,)}{\Lambda(1+2ns-2n,\; f\times \widetilde f\,)}\sum_{\gamma \in \widetilde \Gamma(2n-1)} \frac{1}{2\pi i} \int\limits_{2-i\infty}^{2+i\infty}
  \frac{A^{-\frac{n}{2}w}}{w(w+1)}  \left| \frac{\text{\rm det}( \frak m_1)}{\text{\rm det}( \frak m_2)}  \right|^{n\left(-1+s-\frac{w}{2}\right)}
  \\
  &
  \hskip 227pt \cdot  c_{2-s+\frac{w}{2}}(f) \cdot f(\frak m_1)  f^*(\frak m_2) \; \bigg | _{ \left(\smallmatrix \gamma &\\ & 1\endsmallmatrix\right)}
 \; dw
  \\
  &
  \hskip 60pt +\; ND\left( \widehat E_A^*   \right).
\endalign
 $$
  }
  
  {\bf Proof of corollary 8.7:} The proof follows from the well known identity
  $$\frac{1}{2\pi i}\int\limits_{2-i\infty}^{2+i\infty} \frac{x^w}{w(w+1)}\; dw = \cases 1-\frac{1}{x}, & \text{if} \; x\ge 1,\\
  0, &\text{otherwise},\endcases\tag{8.8}$$
 together with  the identity
 $$\align & \widehat E_A^*(z, f; s) = E(z, f; s) - A^{\frac{n}{2}}E(z,  f;  s-1/2) + \frac{\Lambda(2ns-2n,\; f\times \widetilde f\,)  }{ \Lambda(1+2ns-2n, \;f\times \widetilde f\,) } \cdot E(z,  f;\,  2-s)
 \\
 &
\hskip 50pt  - \frac{1}{2\pi i}\int\limits_{1-i\infty}^{1+i\infty} \frac{ A^{-\frac{n}{2} w}  }{  w(w+1)} E\left(z, f; \; s+\frac{w}{2}\right)\; dw 
\tag{8.9}  \\
  &\; - \;  \frac{\Lambda(2ns-2n,\; f\times \widetilde f\,)  }{ \Lambda(1+2ns-2n, \;f\times \widetilde f\,) }
\cdot \frac{1}{2\pi i}\int\limits_{1-i\infty}^{1+i\infty} \frac{ A^{-\frac{n}{2} w}  }{  w(w+1)} E\left(z, f; \; 2-s+\frac{w}{2}\right)\; dw.\endalign$$
 and proposition 8.6.
  \qed

\vskip 20pt
 \noindent
  {\bf \S 9 Coset representatives for $\widetilde{P}_{n, m}\backslash SL(n, \Bbb Z)$:} 
 \vskip 10pt
 
 Recall definition 8.2 which states that
 $$\widetilde P_{n,m}(\Bbb Z) \, := \left( \pmatrix  & * & * & \cdots & *\\
 SL(n-m, \Bbb Z)& \vdots &\vdots &\cdots&\vdots\\
 & * & * &\cdots & *\\
 & 1 & * &\cdots & *\\
 & &\ddots & &\vdots\\
 & & &1&*\\
 & & & & 1\endpmatrix  \right) \; \subset SL(n, \Bbb Z).\tag{9.1}$$
 
 The main goal of this section is the proof of the following proposition.
\vskip 10pt\noindent
{\bf Proposition 9.2} {\it  Let $A, A' \in SL(n,\Bbb Z).$ Then there exists a matrix $X \in \widetilde P_{n,m}(\Bbb Z)$ such that $$A = XA'$$ if and only if for each $1\le k\le m$, the matrices $A, A'$ have the same $k\times k$ minors on the bottom $k$ rows.}
\vskip 8pt
The proof of proposition 9.2 is based on the following lemma.

\vskip 10pt\noindent
{\bf Lemma 9.3} {\it Suppose $n\ge m \ge 1$ and $a_{i,j}, a_{i,j}' \in\Bbb Z$ (for $1 \le i\le m$, $1\le j\le n).$ Define matrices
$$A = \Big ( a_{i,j}\big)_{1 \le i\le m,\; 1\le j\le n}, \qquad\quad A' = \Big ( a_{i,j}'\big)_{1 \le i\le m,\; 1\le j\le n}.$$ 
If $\text{\rm rank}(A) = \text{\rm rank}(A') = m$, then there exists a rational $m\times m$ matrix
$$U = \pmatrix 1 &* &\cdots & *\\
& 1 & \cdots &*\\
&& \ddots &\vdots\\
&&&1\endpmatrix$$
such that 
$$A = UA'$$
if and only if, for each $1 \le k \le m$, the matrices $A, A'$ have the same $k\times k$ minors in the bottom $k$ rows.}

\vskip 10pt {\bf Proof of lemma 9.3:}  It is easy to show  that if $A = UA'$ then the matrices $A, A'$ have the same $k\times k$ minors in the bottom $k$ rows. In the other direction, we use induction. Before proceeding with the induction we establish two claims.

\vskip 8pt\noindent
{\bf Claim 9.4:} Let $k \in \Bbb Z$ with $k\ge 2.$ Let
$$A = \pmatrix \alpha_1&\alpha_2&\cdots & \alpha_k\\
a_{2,1} & a_{2,2} & \cdots & a_{2,k}\\
\vdots &\vdots &\cdots &\vdots\\
a_{k,1} & a_{k,2} & \cdots & a_{k,k}
\endpmatrix, \quad A' = \pmatrix \alpha_1'&\alpha_2'&\cdots & \alpha_k'\\
a_{2,1} & a_{2,2} & \cdots & a_{2,k}\\
\vdots &\vdots &\cdots &\vdots\\
a_{k,1} & a_{k,2} & \cdots & a_{k,k}
\endpmatrix$$
be two matrices in $GL(k,\Bbb Z)$ satisfying $\text{\rm det}(A) =\text{\rm det}(A') \ne 0.$
Then there exists 
$$U = \pmatrix 1 & u_1 & u_2 &\cdots & u_{k-1}\\
& 1 & 0 & \cdots & 0\\
&  &\ddots  & \vdots & \vdots\\
&  &  & 1 & 0\\
&  &  &  & 1\endpmatrix \in SL(k, \Bbb Q)\tag{9.5}$$
such that
$A' = U\cdot A.$

To prove the claim note that $\text{\rm det}(A) =\text{\rm det}(A') \ne 0$ implies  (by expanding in terms of cofactors $C_{i,1}$ along the first row) that
 $\text{\rm det}(A) = \sum_{i=1}^k \alpha_i C_{i,1} = \sum_{i=1}^k \alpha_i' C_{i,1} = \text{\rm det}(A')$ which implies that
 $$\sum_{i=1}^k (\alpha_i-\alpha_i') C_{i,1} = 0.\tag{9.6}$$
 It immediately follows that if we set $$\alpha_i' = \alpha_i + a_{2,i}u_1 + a_{3,i} u_2 + \cdots a_{k,i} u_{k-1}, \qquad (i = 1,2,\ldots, k),$$
 then we can find rational numbers $u_1, \ldots,u_{k-1}$ so that (9.6) is satisfied. Since $\text{\rm det}(A) \ne 0$, one then concludes that the matrix $U = A' A^{-1}$ must be of the form (9.5).
 
 \vskip 10pt\noindent
 {\bf Claim 9.7} Let $n,k\in \Bbb Z$ with $2\le k \le n.$ Let
 $$A = \pmatrix \alpha_1&\alpha_2&\cdots & \alpha_n\\
a_{2,1} & a_{2,2} & \cdots & a_{2,n}\\
\vdots &\vdots &\cdots &\vdots\\
a_{k,1} & a_{k,2} & \cdots & a_{k,n}
\endpmatrix, \quad A' = \pmatrix \alpha_1'&\alpha_2'&\cdots & \alpha_n'\\
a_{2,1} & a_{2,2} & \cdots & a_{2,n}\\
\vdots &\vdots &\cdots &\vdots\\
a_{k,1} & a_{k,2} & \cdots & a_{k,n}
\endpmatrix,$$
be two matrices with integer coefficients. Assume that the $k\times k$ minors of $A$ are the same as the $k\times k$ minors of $A'$, and that at least one of these minors is non-zero. Then there exists 
$$U = \pmatrix 1 & u_1 & u_2 &\cdots & u_{k-1}\\
& 1 & 0 & \cdots & 0\\
&  &\ddots  & \vdots & \vdots\\
&  &  & 1 & 0\\
&  &  &  & 1\endpmatrix \in SL(k, \Bbb Q)\tag{9.8}$$
such that
$A' = U\cdot A.$

  To prove claim 9.7, note that, without loss of generality, we may assume the first minor
  $$\text{\rm det}  \pmatrix \alpha_1&\alpha_2&\cdots & \alpha_k\\
a_{2,1} & a_{2,2} & \cdots & a_{2,k}\\
\vdots &\vdots &\cdots &\vdots\\
a_{k,1} & a_{k,2} & \cdots & a_{k,k}
\endpmatrix   \; =  \; \text{\rm det}  \pmatrix \alpha_1'&\alpha_2'&\cdots & \alpha_k'\\
a_{2,1} & a_{2,2} & \cdots & a_{2,k}\\
\vdots &\vdots &\cdots &\vdots\\
a_{k,1} & a_{k,2} & \cdots & a_{k,k}
\endpmatrix   \;\ne \;0.$$
It immediately follows from claim 9.4 that there exists $U$ of the form (9.8) such that
$$\pmatrix \alpha_1'&\alpha_2'&\cdots & \alpha_k'\\
a_{2,1} & a_{2,2} & \cdots & a_{2,k}\\
\vdots &\vdots &\cdots &\vdots\\
a_{k,1} & a_{k,2} & \cdots & a_{k,k}
\endpmatrix  = U \cdot \pmatrix \alpha_1&\alpha_2&\cdots & \alpha_k\\
a_{2,1} & a_{2,2} & \cdots & a_{2,k}\\
\vdots &\vdots &\cdots &\vdots\\
a_{k,1} & a_{k,2} & \cdots & a_{k,k}
\endpmatrix $$

Similarly, for any $k < \ell \le n$ we must have
$$\pmatrix \alpha_1'&\alpha_2'&\cdots & \alpha_{k-1}' &\alpha_\ell'\\
a_{2,1} & a_{2,2} & \cdots & a_{2,k-1} & a_{2,\ell}\\
\vdots &\vdots &\cdots &\vdots &\vdots\\
a_{k,1} & a_{k,2} & \cdots & a_{k,k-1} & a_{k,\ell}
\endpmatrix  = U \cdot \pmatrix \alpha_1&\alpha_2&\cdots & \alpha_{k-1} &\alpha_\ell\\
a_{2,1} & a_{2,2} & \cdots & a_{2,k-1} &a_{2,\ell}\\
\vdots &\vdots &\cdots &\vdots&\vdots\\
a_{k,1} & a_{k,2} & \cdots & a_{k,k-1} &a_{k,\ell}
\endpmatrix. $$
This is because a solution $U$ must exist even if the determinants of both minors are zero, and $u_1, \ldots,u_{k-1}$ are already determined by the first $k-1$ columns and the previous computation.
This implies that $\alpha_\ell' = \alpha_\ell + a_{2,\ell}u_1 + \cdots a_{k,\ell} u_{k-1}$ for all $a \le \ell \le n.$ This completes the proof of claim 9.7.

Let's assume that 
$$A = \Big ( a_{i,j}\big)_{1 \le i\le m,\; 1\le j\le n}, \qquad\quad A' = \Big ( a_{i,j}'\big)_{1 \le i\le m,\; 1\le j\le n}$$
 have the same $k\times k$ minors in the bottom $k$ rows. We want to prove there exists a unipotent matrix $U$ such that $A = UA'.$ For the $1\times 1$ minors on the bottom row this implies that $a_{m,i} = a_{m,i}'$ for $i = 1,2,\ldots n.$ For the $2\times 2$ minors on the bottom two rows we may use claim 9.7 to conclude that
 $$\pmatrix a_{m-1,1}' & \cdots & a_{m-1,n}'\\
 a_{m,1} & \cdots & a_{m,n}\endpmatrix \; = \; 
 \pmatrix 1 & u_0\\
 & 1\endpmatrix  
 \pmatrix a_{m-1,1} & \cdots & a_{m-1,n}\\
 a_{m,1} & \cdots & a_{m,n}\endpmatrix,$$
 for some rational number $u_0$.
 
 Next, we consider the $3\times 3$ minors on the bottom 3 rows. We want to prove there exists an upper triangular unipotent $3\times 3$ matrix $U$ such that
 $$\pmatrix a_{m-2,1}' & \cdots & a_{m-2,n}'\\
 a_{m-1,1}' & \cdots & a_{m-1,n}'\\
 a_{m,1} & \cdots & a_{m,n}\endpmatrix \; = \; 
 U\cdot  
 \pmatrix a_{m-2,1} & \cdots & a_{m-2,n}\\
 a_{m-1,1} & \cdots & a_{m-1,n}\\
 a_{m,1} & \cdots & a_{m,n}\endpmatrix.$$
 We may set
 $$U = \pmatrix 1 &u_1 & u_2\\
 & 1&u_0\\
 & & 1\endpmatrix \; = \;\pmatrix 1 &0 & 0\\
 & 1&u_0\\
 & & 1\endpmatrix\cdot \pmatrix 1 &u_1 & u_2\\
 & 1&0\\
 & & 1\endpmatrix,$$
 which implies that we need to prove

 $$\align & \pmatrix 1 &0 & 0\\
 & 1&u_0\\
 & & 1\endpmatrix^{-1}\pmatrix a_{m-2,1}' & \cdots & a_{m-2,n}'\\
 a_{m-1,1}' & \cdots & a_{m-1,n}'\\
 a_{m,1} & \cdots & a_{m,n}\endpmatrix \; =  \pmatrix a_{m-2,1}' & \cdots & a_{m-2,n}'\\
 a_{m-1,1} & \cdots & a_{m-1,n}\\
 a_{m,1} & \cdots & a_{m,n}\endpmatrix\;  \\
 &
 \hskip 150pt
 \; = \;\pmatrix 1 &u_1 & u_2\\
 & 1&0\\
 & & 1\endpmatrix\cdot  
 \pmatrix a_{m-2,1} & \cdots & a_{m-2,n}\\
 a_{m-1,1} & \cdots & a_{m-1,n}\\
 a_{m,1} & \cdots & a_{m,n}\endpmatrix.\endalign$$
 But this again follows from claim 9.7. It is clear that we may proceed by induction to complete the proof of lemma 9.3.

   \qed
  
  \vskip 10pt
  {\bf Proof of proposition 9.2:}  Let $A, A' \in SL(n,\Bbb Z).$ Assume there exists $X\in \widetilde P_{n,m}(\Bbb Z)$ such that $A = XA'$. Then it is obvious that for each $1\le k\le m$ the matrices $A,A'$ have the same $k\times k$ minors in the bottom $k$ rows. 
  
  In the other direction, let
  $$A = \pmatrix a_{1,1} &\cdots & a_{1,n}\\
  \vdots & \cdots &\vdots\\
  a_{n-m,1} &\cdots & a_{n-m,n}\\
  a_{n-m+1,1} &\cdots & a_{n-m+1,n}\\
  \vdots & \cdots &\vdots\\
  a_{n,1} & \cdots & a_{n,n}
  \endpmatrix  \; := \; \pmatrix M_{(n-m)\times n}\\
  &\\
  N_{m\times n}\endpmatrix,$$
  
  $$A' = \pmatrix a_{1,1}' &\cdots & a_{1,n}'\\
  \vdots & \cdots &\vdots\\
  a_{n-m,1}' &\cdots & a_{n-m,n}'\\
  a_{n-m+1,1}' &\cdots & a_{n-m+1,n}'\\
  \vdots & \cdots &\vdots\\
  a_{n,1}' & \cdots & a_{n,n}'
  \endpmatrix  \; := \; \pmatrix M_{(n-m)\times n}'\\
  &\\
  N_{m\times n}'\endpmatrix.$$
  By lemma 9.3, there exists a rational  $m\times m$ matrix
  $$U = \pmatrix 1 & * & * &\cdots & *\\
  & 1 & * & \cdots & *\\
  & & \ddots &\vdots &\vdots\\
  & & & 1 & *\\
  & & & & 1
  \endpmatrix$$
  such that
  $$N = UN'.$$
  Let
  $$\big(X_1,\; X_2\big) := M {A'}^{-1}.$$
  Clearly
  $$X := \pmatrix X_1 & X_2\\
  & U\endpmatrix \in \widetilde P_{n,m}(\Bbb Z).$$
  \qed
  
  \vskip 20pt
 \noindent
  {\bf \S 10 Bounds for Eisenstein series:} 
 \vskip 10pt
  
 Fix the parabolic $\Cal P = P_{n,n}$ in $GL(2n, \Bbb R).$ Let $z = \frak n\frak m k \in GL(2n, \Bbb R)$ with
  $$\frak n \in N^{\Cal P}, \qquad \frak m = \pmatrix \frak m_1 &\\& \frak m_2\endpmatrix \in M^{\Cal P}, \qquad k \in K.$$
  \vskip 10pt\noindent
  {\bf Proposition 10.1} {\it Let $f$ be a  Maass form on $GL(n, \Bbb Z)$ as in (4.1) with Whittaker function $W_{n,\nu}$, as in proposition 8.6, and set
   $$M = \text{\rm diag}\big(m_{2n-1}m_{2n-2}\cdots |m_{n+1}|, \;\;\ldots, \;\;m_{2n-1}m_{2n-2}, m_{2n-1},1\big).$$ Then
  $$\sum_{\gamma \in \widetilde P_{2n-1,n-1}(\Bbb Z)\backslash SL(2n-1, \Bbb Z)} f(\frak m_1) W_{n,\nu}\left(MY_2 \right) \; \Bigg |_{ \left(\smallmatrix \gamma &\\ & 1\endsmallmatrix\right)} \; \ll \; 1$$
 for any $z \in GL(2n, \Bbb R)$ for which $|x_{i,j}|\le 1,$ and $y_i\ge 1 $ with $ (1\le i\le j \le 2n)$. The upper bound above does not depend on $M$.
  }
  
  \vskip 10pt
  {\bf Proof of proposition 10.1:}    For $1\le i\le n-1$, set
  $$\mu^i_j = \cases  -\frac{1}{i+1}, & \text{if} \; j = 2n-i-1,\\
 \;\;\; 1, & \text{if} \;  j = 2n-1,\\
\;\;\;  0, & \text{otherwise}.\endcases$$
 Let $\mu^i = (\mu^i_1, \ldots,\mu^i_{2n-1}).$ Then
 $$I_{\mu^i}(z) = \prod_{\ell=1}^i  y_\ell^{\frac{2n}{i+1}(i+1-\ell)}.$$
It follows from lemma 5.7.2 \cite{Goldfeld, 2006} that for
$$\gamma' = \pmatrix \gamma &\\& 1\endpmatrix \; \in \; SL(2n, \Bbb Z),$$
with $\gamma \in SL(2n-1, \Bbb Z)$, that we have
 $$I_{\mu^i}(\gamma' z) = \big|\big| e_{2n-i}\gamma' z \wedge \cdots \wedge e_{2n-1}\gamma' z \wedge e_{2n} \gamma' z\big|\big|^{\frac{2n}{i+1}}\; \cdot \;\big|\big| e_{2n-i} z\wedge \cdots \wedge e_{2n} z\big|\big|^{\frac{-2n}{i+1}}\; I_{\mu^i}(z).\tag{10.2}$$

\vskip 8pt\noindent 
{\bf Remark 10.3} {\it For any $2n \times 2n$ matrix $M$, the wedge product
  $\big |\big | e_{2n-i} M\wedge \cdots \wedge e_{2n} M\big |\big |^2$
 is the sum of the squares of all the $(i+1)\times (i+1)$ minors of the last $i+1$ rows of $M$.}
\vskip 8pt 
  Consequently, for $z\in GL(2n, \Bbb R)$ in upper triangular Iwasawa form, we have
 $$\big|\big| e_{2n-i} z\wedge \cdots \wedge e_{2n} z\big|\big| \; = \; \big[z\big]_{(2n-i, \ldots,2n),\, (2n-i, \ldots,2n)} = \;\prod_{\ell=1}^i y_\ell^{i+1-\ell},\tag{10.4}$$
  where for a matrix $M$, and for vectors of integers $I,J$, we denote the minor of $M$ determined by the rows $I$ and columns $J$ as $[M]_{I,J}$ (see lemma 4.11). 
  
  It immediately follows from (10.2), (10.4)  and remark 10.3 (with $M = \gamma' z$)  that
  $$\align  I_{\mu^i}(\gamma' z) & = \big|\big| e_{2n-i}\gamma' z \wedge \cdots \wedge e_{2n-1}\gamma' z \wedge e_{2n} \gamma' z\big|\big|^{\frac{2n}{i+1}}\\
  & = \left(\big[C\big]_{\alpha_i, \,\beta_{i,1}}^2 \; + \; \big[C\big]_{\alpha_i, \,\beta_{i,2}}^2 \; + \; \cdots\; + \big[C\big]_{\alpha_i, \,\beta_{i,L}}^2  \right)^{\frac{n}{i+1}},\endalign$$
  where
  $$\align & C := \gamma' z,\\
  & \alpha_i := (2n-i, 2n-i+1, \cdots ,2n),\\
 & \beta_{i,1} := (1,2, \ldots, i+1),\\
 & \beta_{i,2} := (1,2, \ldots,i, \, i+2),\\
 & \hskip 24pt \vdots\\
 & \beta_{i,L} := \alpha_i =  (2n-i, 2n-i+1, \ldots, 2n). \endalign$$

  By the Cauchy-Binet formula (4.12) and lemma 4.11 it follows that
  
  $$\align  \big[C\big]_{\alpha_i, \,\beta_{i,1}} & = \big[\gamma'\big]_{\alpha_i,\beta_{i,1}} \;\big[z\big]_{\beta_{i,1},\beta_{i,1}}\\
  \big[C\big]_{\alpha_i, \,\beta_{i,2}} & = \big[\gamma'\big]_{\alpha_i,\beta_{i,1}} \;\big[z\big]_{\beta_{i,1},\beta_{i,2}} \; +  \;\big[\gamma'\big]_{\alpha_i,\beta_{i,2}} \;\big[z\big]_{\beta_{i,2},\beta_{i,2}}\\
  &\\
  & \hskip 5pt \vdots\\
  &\\
   \big[C\big]_{\alpha_i,\alpha_i} & = \big[\gamma'\big]_{\alpha_i,\beta_{i,1}} \big[z\big]_{\beta_{i,1},\alpha_i} \; + \; \big[\gamma'\big]_{\alpha_i,\beta_{i,2}} \big[z\big]_{\beta_{i,2},\alpha_i}\\
   & \hskip 120pt+ \;\cdots \; + \; \big[\gamma'\big]_{\alpha_i,\alpha_i} \big[z\big]_{\alpha_i,\alpha_i}\endalign$$
  
 If $\big[\gamma'\big]_{\alpha_i,\beta_{i,1}} \ne 0$ then
 $$\big[C\big]_{\alpha_i,\beta_{i,1}}^2 \; \ge \; \big[z\big]_{\beta_{i,1},\beta_{i,1}} \; \ge \; 1.$$
 
 If $\big[\gamma'\big]_{\alpha_i,\beta_{i,1}} = 0$ then
 $$\big[C\big]_{\alpha_i,\beta_{i,2}}^2 \; \ge \; \big[z\big]_{\beta_{i,2},\beta_{i,2}} \; \ge 1,$$
 and by repeating the previous argument, we have
 $$I_{\mu^i}\left(\gamma' z\right) \; \ge \; 1.$$
  Hence
  $$\max_{1\le \ell \le i} y_\ell\left(\gamma' \right) \ge 1.$$
  We adopt the notation $$y_i(\gamma') :=  y_i \big |_{\gamma'}, \qquad (\text{for}\; 1 \le i\le n-1).$$ 
  Now it follows from corollary 7.4 that 
  $$\align  W_{n,\nu}(MY_2)\Big |_{\left(\smallmatrix \gamma &\\&1\endsmallmatrix\right)} & \ll \left(\max_{1\le j<n}  y_j(\gamma')\right)^{-Nn^3} \ll \left| I_{\mu^i} \left( \gamma' z\right) \right|^{-N}
  \\
  &
 \ll \left( \big[C\big]_{\alpha_i, \beta_{i,1}}^2 \; + \;  \big[C\big]_{\alpha_i, \beta_{i,2}}^2 \; + \;\; \cdots \;\; + \big[C\big]_{\alpha_i, \beta_{i,L}}^2
\right)^{-N},   \endalign$$ 
  for arbitrary large fixed $N > 1,$ and every $1 \le i < n.$

  Since $f$ is a cusp form, it is absolutely bounded, and it immediately follows from the above calculations and proposition 9.2  that  $$\align & \sum_{\gamma \in \widetilde P_{2n-1,n-1}(\Bbb Z)\backslash SL(2n-1, \Bbb Z)} f(\frak m_1) \; W_{n,\nu}\left(MY_2 \right) \; \Bigg |_{ \left(\smallmatrix \gamma &\\ & 1\endsmallmatrix\right)} \\
  &
  \\
  & \hskip 44pt \ll\; \sum_{\gamma \in \widetilde P_{2n-1,n-1}(\Bbb Z)\backslash SL(2n-1, \Bbb Z)} W_{n,\nu}\left(MY_2 \right) \; \Bigg |_{ \left(\smallmatrix \gamma &\\ & 1\endsmallmatrix\right)}\\
  &
  \\
  &
   \hskip 44pt \ll\; \sum_{\gamma \in \widetilde P_{2n-1,n-1}(\Bbb Z)\backslash SL(2n-1, \Bbb Z)} \left( \big[C\big]_{\alpha_i, \beta_{i,1}}^2 \; + \;  \big[C\big]_{\alpha_i, \beta_{i,2}}^2 \; + \;\; \cdots \;\; + \big[C\big]_{\alpha_i, \beta_{i,L}}^2
\right)^{-N}\\
&
\\
&
 \hskip 44pt \ll\; 1. \endalign$$

  \qed
  
   \vskip 10pt\noindent
  {\bf Proposition 10.5} {\it Let $z \in GL(k,\Bbb R)/O(k,\Bbb R)\cdot\Bbb R^\times$ with Iwasawa decomposition $z = xy$ as in (2.1) where
  $-\frac12 \le x_{i.j} < \frac12, \; y_i \ge 1, \; (1\le i<j\le k).$
  Let $$M = \text{\rm diag}\big(m_{2n-1}m_{2n-2}\cdots |m_{n-1}|, \;\;\ldots, \;\;m_{2n-1}m_{2n-2}, m_{2n-1},1\big),$$ and for spectral parameters $\nu = (\nu_1, \ldots,\nu_{k-1}) \in \Bbb C^{k-1},$
define  $$\align & F(z) := \sum_{\gamma\in U_{k-1}(\Bbb Z)\backslash SL(k,\Bbb Z)}\; \sum_{m_1=1}^\infty \cdots \sum_{m_{k-1}\ne 0}   \frac{A(m_1, \ldots ,m_{k-1})  }{ \prod\limits_{\ell=1}^{k-1}|m_\ell|^{\frac{\ell(k-\ell)}{2}} }\; W_{k,\nu}\left(M y\right)\\
  & \hskip 190pt \cdot
  e^{2\pi i(m_1x_{1,2}+m_2x_{2,3} + \cdots + m_{k-1} x_{k-1,k})}\; \bigg |_{\left(\smallmatrix \gamma &\\&1\endsmallmatrix    \right)    },\endalign$$
  where the coefficients $A(m_1, \ldots,m_{k-1})\in\Bbb C$ satisfy
  $$\sum_{\left |m_1^{k-1}\cdots m_{k-1}\right |\le N} |A(m_1, \ldots,m_{k-1})|^2 \; \ll \; N.$$
   Assume the Langlands parameters $(\alpha_1, \ldots,\alpha_k)$ attached to $W_{k,\nu}$ (see definition (2.3))  satisfy 
  $\; |\Im(\alpha_i)| \ll |t|$ as $|t|\to\infty,$ for all $1\le i\le k.$ Then for
  $$\max_{1 \le i \le k-1} y_i\; \ge \; \left( |t|^{1+\epsilon}  \right)^{\frac{k(k-1)}{2}}$$
  we have
  $$F(z) \ll  \left(\max_{1\le i\le k-1} y_i\right)^{-N}$$
 for any fixed constant $N>1.$ }
 
 \vskip 10pt
 The proof of proposition 10.5 requires the following lemma.
 
 \vskip 10pt
 {\bf Lemma 10.6} {\it Let $\gamma' = \pmatrix \gamma &\\&1\endpmatrix \in SL(k, \Bbb Z)$ for $k\ge 2.$   Assume that $\min\limits_{1\le i<k} y_i \ge 1.$ Then
 $$\max\limits_{1\le i<k} y_i(\gamma') \; \ge \; \left(\max\limits_{1\le i<k} y_i\right)^{\frac{2}{(k-1)k}}.$$}
 
 \vskip 10pt 
 {\bf Proof of lemma 10.6:} Set $\mu = (0,\ldots,0,1)$  and let $z = xy \in GL(k,\Bbb R)/(O(k,\Bbb R)\cdot\Bbb R^\times)$. It follows from (2.2) that
 $$I_\mu(z) = \prod_{\ell=1}^{k-1} y_\ell^{k-\ell} \; \ge \; \max_{1\le i\le k-1} y_i.$$
 On the other hand,
 $$I_\mu(\gamma' z) = \prod_{\ell=1}^{k-1} y_\ell(\gamma')^{k-\ell} \; \le \; \left(\max_{1\le i<k} y_i(\gamma')\right)^{\frac{(k-1)k}{2}}.$$

\noindent
{\bf Claim:}  $I_\mu(\gamma' z) = I_\mu(z).$ 

\vskip 5pt
The claim is a consequence of the fact that
$$I_\mu(\gamma' z) = \prod_{i=0}^{k-2} \big|\big|e_{k-i} \gamma' z \wedge \cdots \wedge e_k \gamma' z\big|\big|^{-k\mu_{k-i-1}} \cdot
\big|\big|e_{k-i}  z \wedge \cdots \wedge e_k  z\big|\big|^{k\mu_{k-i-1}}\, I_\mu(z),$$
and that $\mu_{k-1} =1$ only if $i=0,$ (the other ${\mu_i}'$s are zero).
Hence, only the last row of $z$ contributes which establishes the claim. 
 
 It follows that
 $$\max_{1\le i<k} y_i(\gamma') \ge \left(\max_{1\le i<k} y_i
 \right)^{\frac{2}{(k-1)k}}. \qed$$
 
 When $k = 3,$ the above lemma 10.6 and proposition 10.5 are due to \cite{Brumley-Templier} with a better power.
 
 \vskip 10pt
 {\bf Proof of proposition 10.5:}  
 For every $1\le i\le k-2$, we choose $\mu^i:=(\mu_1^i, \ldots,\mu_{k-1}^i)$ where
 $$\mu_j^i = \cases \hskip-5pt -\frac{1}{i+1}, & \text{if} \; j = k-i-1,\\
 \;1, & \text{if} \; j = k-1,\\
 \; 0, &\text{otherwise}.\endcases$$
  With this choice, we see that
  $$I_{\mu^i}(z)  = \prod_{\ell=1}^i \;y_\ell^{\frac{i-\ell+1}{i+1} k}.\tag{*} $$
  
 It follows from lemma 5.7.2 in \cite{Goldfeld, 2006} that for
  $$\gamma' = \pmatrix\gamma & \\& 1\endpmatrix \in SL(k, \Bbb Z)$$
we have
$$I_{\mu^i}(\gamma' z) = \big|\big| e_{k-i}\gamma' z \wedge \;\cdots\; \wedge e_k \gamma' z\big|\big|^{\frac{k}{i+1}} \cdot \big|\big| e_{k-i} z \wedge \;\cdots\; \wedge e_k z\big|\big|^{-\frac{k}{i+1}} \cdot I_{\mu^i}( z).$$
Since
$$\big|\big| e_{k-i} z \wedge \;\cdots\; \wedge e_k z\big|\big| = \big[z\big]_{(k-i,\ldots,k)(k-i,\ldots,k)} = \prod_{\ell=1}^i y_\ell^{i-\ell+1},$$
if we combine the above with (*), we obtain
$$I_{\mu^i}(\gamma' z) = \big|\big| e_{k-i}\gamma' z \wedge \;\cdots\; \wedge e_k \gamma' z\big|\big|^{\frac{k}{i+1}} .$$
Since
 $$I_{\mu^i}(\gamma'  z) \ll  \left(\max_{1\le i \le k-1} y_i(\gamma')\right)^{2k^2}$$
and
$$\max_{1 \le i \le k-1} y_i \ge \left(t^{1+\epsilon}\right)^{\frac{k(k-1)}{2}},$$
 it follows from corollary 7.4 and the above that
 $$\align  W_{k,\nu}\big(M\cdot y(\gamma)\big) \ll \max_{1\le i\le k-1}\big(m_i \,y_i(\gamma')\big)^{-N_1} & \ll I_{\mu^i}(\gamma' z)^{\frac{-N_1}{4k^2}}\left( \max_{1\le i\le k-1} y_i  \right)^{\frac{-N_1}{k(k-1)}}\\
 & = \big|\big|e_{k-i}\gamma' z\wedge \; \cdots\; e_k\gamma' z\big|\big|^{-N'} \cdot\left( \max_{1\le i\le k-1} y_i  \right)^{-N}\endalign$$
 where $N' = \frac{N_1}{4k(i+1)}$ and $N = \frac{N_1}{k(k-1)}$ are very large.
 \vskip6pt
 Hence, for $k \ge 4$, we obtain
 $$\align  F(z) & \ll \sum_{i=1}^{k-2} \sum_{\gamma \in U_{k-1}(\Bbb Z)\backslash SL(k-1,\Bbb Z)}
 \;\sum_{m_1\ge 1} \cdots\sum_{m_{\ell-1}\ne0} \frac{|A(m_1,\ldots,m_{k-1})|}
 {\prod\limits_{\ell=1}^{k-1} |m_\ell|^{\frac{\ell(k-\ell)}{2}}     }\\
 &\hskip 180pt\cdot  \big|\big|e_{k-i}\gamma' z\wedge\;\cdots\; \wedge e_k\gamma' z\big|\big|^{-N'}
 \left(\max y_i\right)^{-N}\\
 & \ll \left(\max_{1\le i \le k-1} y_i\right)^{-N},\endalign$$
 by theorem 11.3.2 in \cite{Goldfeld, 2006} and similar arguments as in the proof of proposition 10.1.
 
 For $k = 3$ and
 $$\gamma' = \pmatrix a& b &\\
 c& d\\
 & & 1\endpmatrix$$
 we have
 $$y_1(\gamma') = |cz_2+d| y_1, \qquad y_2(\gamma') = \frac{y_2}{|cz_2+d|^2}.$$
 Hence $$y_1(\gamma')^2 \,y_2(\gamma') = y_1^2y_2.$$
 It follows that
 $$\align \max_{1\le j\le 2}  m_j  y_j(\gamma') \; & \ge \; \Big(m_1^2\, m_2 \, y_1(\gamma')^2 \, y_2(\gamma')\Big)^\frac13 \; = \;\big(m_1^2m_2 \, y_1^2 \, y_2\big)^\frac13\\
 & \ge \left(m_1^2m_2\right)^\frac13\endalign$$
 and
 $$\align W\left(M \, y(\gamma')\right) & \ll \Big(\max_{1\le j\le 2} m_j\, y_j(\gamma')\Big)^{-1} \cdot \Big(\max_{1\le j\le 2} y_j(\gamma')\Big)^{-N_1}\\
 & \ll \left(m_1^2 m_2   \right)^{-\frac13} \cdot \big| \big| e_2 \gamma' z\wedge e_3\gamma' z\big| \big |^{-N'} \, \Big(\max_{1\le j\le 2}  y_j \Big)^{-N},\endalign$$
 for any large $N, N'.$
 Hence
  $$\align F(z) & \ll \sum_{\gamma\in U_2(\Bbb Z)\backslash SL(2,\Bbb Z)} \;\sum_{m_1\ge 1}\;\sum_{m_2\ne 0} \frac{|A(m_1,m_2)|}{|m_1|^\frac53 \, |m_2|^\frac43} \; \big|\big | e_2 \gamma' z  \wedge    e_3\gamma' z\big| \big|^{-N'} \cdot \Big(\max_{1\le j\le 2} y_j\Big)^{-N}\\
  & \ll \Big(\max_{1\le j\le 2} y_j\Big)^{-N}.\endalign$$
 \qed
 
 \vskip 10pt\noindent
{\bf Proposition 10.7} {\it Let us rewrite the identity for $\widehat E_A^*$ given in corollary 8.7 as $$\widehat E_A^*(z, f; s) = A_1+A_2+A_3+A_4 -A_5-A_6 +ND\left(\widehat E_A^\star\right),$$
where $$A_1 =  \underset
 h\left(\left(\smallmatrix \gamma &\\ & 1\endsmallmatrix\right) z\right) < A\to {\sum_{\gamma \in \widetilde \Gamma(2n-1)}}
 \left| \frac{\text{\rm det}( \frak m_1)}{\text{\rm det}( \frak m_2)}  \right|^{ns} \left[ 1 - \frac{A^{n/2}}{h(y)^{n/2}}   \right] f(\frak m_1)  f^*(\frak m_2) \; \bigg | _{ \left(\smallmatrix \gamma &\\ & 1\endsmallmatrix\right)}
 $$ and $A_2,A_3,A_4,A_5,A_6$ are defined similarly. Recall that $h(z) =  \left|\frac{\text{\rm det}( \frak m_1)}{\text{\rm det}( \frak m_2)}\right|$. The following bounds hold.
 \pagebreak
 \vskip 8pt
{\bf (i)} If $h(z) \ge A$ and $\min\limits_{1\le i <2n} y_i \ge \frac{\sqrt{3}}{2},$ then we have $A_1, A_2 \ll A^n \ll A^{\frac{n}{2}} h(z)^{\frac{n}{2}}\left(  \frac{h(z)}{A} \right)^{-\frac12}$,\break \phantom{xx} \hskip 220pt  $A_i \ll 1\; \text{for} \;3\le i \le 6$.
 \vskip 5pt
{\bf (ii)} If $h(z) \le A $ and $\min\limits_{1\le i <2n}  y_i \ge \frac{\sqrt{3}}{2}$, then
$$A_1 \; \ll \; A^{\frac{n}{2}}\cdot\left(y_n^n y_{n+1}^{n-1} y_{n+2}^{n-2}\cdots y_{2n-1}   \right)^{\frac{n}{2}} := A', \qquad A_2\ll A', \qquad A_3,A_4,A_5,A_6 \ll1.$$
\vskip 5pt
{\bf (iii)} If $\max\limits_{1\le i < 2n} y_i \ge\left(|t|^{1+\epsilon}\right)^{n(2n-1)}$ and $\min\limits_{1\le i <2n} y_i \ge \frac{\sqrt{3}}{2},$ and $g \in \Bbb C_0^\infty\left( [1,2]  \right)$, then
$$\int\limits_0^\infty g\left( \frac{A}{\beta}  \right) ND\left(\widehat{E}_A^*(z,f; 1+it)  \right) \; \frac{dA}{A} \ll |t|^{-N}$$
for any fixed $N > 1.$ }
 
 \vskip 10pt
 {\bf Proof of proposition 10.7:} 
 \vskip 8pt
\hskip 15pt $\underline{\text{\bf Part (i)}}$ We first prove part (i) for $A_1.$  It follows from proposition 10.1 that
 $$A_1 \ll A^n \sum_{\gamma\in \widetilde\Gamma(2n-1)} \left| f(\frak m_1)  f^*(\frak m_2)    \right| \;\bigg |_{\left(\smallmatrix \gamma &\\& 1\endsmallmatrix\right)} \; \ll \; A^n \ll A^{\frac{n}{2}} h(z)^{\frac{n}{2}} \left(\frac{h(z)}{A}   \right)^{-\frac12}.$$
   
   Next, we prove the bounds for $A_2.$ Recall that $s = 1+it.$ Then
   $$\left|\frac{\Lambda(2ns-2n, \; f\times \widetilde f\,)}{   \Lambda(1+2ns-2n, \; f\times \widetilde f\,)  }\right|\; = \; 1,$$
   the same argument as the above shows $A_2 \ll \;A^n \ll A^{\frac{n}{2}} h(z)^{\frac{n}{2}} \left(\frac{h(z)}{A}   \right)^{-\frac12}.$

To obtain the bound for $A_3$ we note that
$$|c_s(f)| = \left|\frac{\Lambda(2ns-n, \; f\times \widetilde f\,)}{   \Lambda(1+2ns-n, \; f\times \widetilde f\,)  }\right| \ll 1,$$
from which it follows that
$$A_3 \ll \sum_{\gamma \in \widetilde \Gamma(2n-1)} \big| f(\frak m_1) f^*(\frak m_2)   \big| \;\bigg |_{\left(\smallmatrix \gamma & \\&1\endsmallmatrix   \right)} \ll 1.$$
by proposition 10.1.
To obtain the bound for $A_4$, note that
 $$\left|\frac{\Lambda(2ns-2n, \; f\times \widetilde f\,)}{   \Lambda(1+2ns-2n, \; f\times \widetilde f\,)  }\right|\; = \; 1, \qquad |c_{2-s}(f)| \ll 1,$$
 so by the same argument as for $A_3$, we obtain $A_4 \ll 1.$
 
 To obtain the bound for $A_5$ we write $A_5 = A_{5,1} + A_{5,2}$ where
 (for $\gamma' =  \left(\smallmatrix \gamma &\\ & 1\endsmallmatrix\right))$
 $$A_{5,1} = \underset  h(\gamma' z) < \frac{1}{A}  \to{\sum_{\gamma \in \widetilde \Gamma(2n-1)}} \frac{1}{2\pi i} \int\limits_{2-i\infty}^{2+i\infty}
  \frac{A^{-\frac{n}{2}w}}{w(w+1)} \left| \frac{\text{\rm det}( \frak m_1)}{\text{\rm det}( \frak m_2)}  \right|^{n\left(1-s-\frac{w}{2}\right)} c_{s+\frac{w}{2}}(f) \cdot   f(\frak m_1) f^*(\frak m_2) \; \bigg | _{\gamma'}
 \; dw
$$
 and
 $$A_{5,2} = \underset  h(\gamma' z) \ge \frac{1}{A}  \to{\sum_{\gamma \in \widetilde \Gamma(2n-1)}} \frac{1}{2\pi i} \int\limits_{2-i\infty}^{2+i\infty}
  \frac{A^{-\frac{n}{2}w}}{w(w+1)} \left| \frac{\text{\rm det}( \frak m_1)}{\text{\rm det}( \frak m_2)}  \right|^{n\left(1-s-\frac{w}{2}\right)} c_{s+\frac{w}{2}}(f) \cdot  f(\frak m_1) f^*(\frak m_2) \; \bigg | _{\gamma'}
 \; dw.
$$
It is clear that $A_{5,2} \ll 1$ by proposition 10.1. For $A_{5,1}$ shift the line of integration to $\Re(w) = -\frac12,$ and write $A_{5,1} = A_{5,1,1} + A_{5,1,2},$ where $A_{5,1,1}$ is the residue at $w=0$ given by
$$A_{5,1,1} = \underset  h(\gamma' z) < \frac{1}{A}  \to{\sum_{\gamma \in \widetilde \Gamma(2n-1)}} h(z)^{n(1-s)} c_s(f) \cdot  f(\frak m_1) f^*(\frak m_2) \; \bigg | _{ \gamma'} \; \ll \; 1,
$$
by proposition 10.1. Further
$$A_{5,1,2} = \underset  h(\gamma' z) < \frac{1}{A}  \to{\sum_{\gamma \in \widetilde \Gamma(2n-1)}} \frac{1}{2\pi i} \int\limits_{-\frac12-i\infty}^{-\frac12+i\infty}\hskip -5pt
  \frac{A^{-\frac{n}{2}w}}{w(w+1)} \left| \frac{\text{\rm det}( \frak m_1)}{\text{\rm det}( \frak m_2)}  \right|^{n\left(1-s-\frac{w}{2}\right)} \hskip - 5pt c_{s+\frac{w}{2}}(f) \cdot  f(\frak m_1) f^*(\frak m_2) \; \bigg | _{ \gamma'}
 \; dw,
$$
and again $A_{5,1,2} \ll 1$ by proposition 10.1.

The case of $A_6$ is similar to $A_5.$ The only difference is that
 $$\left|\frac{\Lambda(2ns-2n, \; f\times \widetilde f\,)}{   \Lambda(1+2ns-2n, \; f\times \widetilde f\,)  }\right|\; = \; 1.$$
The same arguments as for $A_5$ then give $A_6 \ll 1.$

\vskip 8pt
\hskip 15pt $\underline{\text{\bf Part (ii)}}$ Next we prove the bounds in part (ii) of proposition 10.7. As in the proof of part (i), let
 $$\mu' = \Big(0,\ldots,0,\underset n^{th}\; \text{position} \to{ \underbrace{-\frac{1}{2n}, \;}} 0,\ldots,\underset (2n-1)^{th}\; \text{position} \to{ \underbrace{\frac12}}\Big).$$
Then
$I_{\mu'}(z) = y_1^{n-1} y_2^{n-2} \cdots y_{n-1} \gg 1$
and
$I_{\mu'}(\gamma' z) = Z\cdot I_{\mu'}(z).$
The proof of proposition 4.10 shows that $Z \ge 1.$ Hence,
$I_{\mu'}(\gamma' z) \ge I_{\mu'}(z).$

Now
$$\left(\max_{1\le i < n-1} y_i(\gamma' z)\right)^{n^3} \ge I_{\mu'}(\gamma' z),$$
and
$$A_1 \ll A^{\frac{n}{2}}\sum_{\gamma \in \widetilde\Gamma(2n-1)} h(z)^{\frac{n}{2}} \big|  f^*(\frak m_2)\big | \;\; \bigg | _{ \left(\smallmatrix \gamma &\\ & 1\endsmallmatrix\right)}.
$$
By proposition 4.10
$$h(\gamma' z)^{\frac{n}{2}} \le h(z)^{\frac{n}{2}} = \left(y_n^n y_{n+1}^{n-1}\cdots y_{2n-1}\right)^{\frac{n}{2}} \left(y_1y_2^2\cdots y_{n-1}^{n-1}   \right)^{\frac{n}{2}},$$
and since
$$f^*\left(\frak m_2(\gamma')\right) \ll \left( \max_{1\le i < n} y_i\left(\gamma' z\right)  \right)^{-N-n^6}$$
we have
$$\align  A_1 & \ll A'\left(y_1 y_2^2 \cdots y_{n-1}^{n-1}   \right)^{\frac{n}{2}} \cdot \sum_{\gamma \in \widetilde\Gamma(2n-1)} \left( \max_{1\le i < n} y_i\left(\gamma' z\right)  \right)^{-N-n^6}
\\
&
\ll A'\left(y_1 y_2^2 \cdots y_{n-1}^{n-1}   \right)^{\frac{n}{2}} I_{\mu'}(z)^{-n^3}\cdot \sum_{\gamma \in \widetilde\Gamma(2n-1)} \left( \max_{1\le i < n} y_i\left(\gamma' z\right)  \right)^{-N}
\\
&
\ll A'
\endalign$$
The same argument also shows that $A_2 \ll A'.$ The proof of the bounds $A_3,A_4,A_5,A_6\ll 1$, are the same as in the proofs in part (i).

\vskip 8pt
\hskip 15pt $\underline{\text{\bf Part (iii)}}$ Finally we prove the third part of proposition 10.7. We have
$$\align & \int\limits_0^\infty g\left(\frac{A}{\beta}\right) \widehat E_A^*(z, f; s) \; \frac{dA}{A}\\
& \hskip 30pt= E_0(z, f; s) \int\limits_0^\infty g(A) \frac{dA}{A} - \frac{1}{2\pi i}\int\limits_{2-i\infty}^{2+i\infty} \frac{\beta^{-\frac{n}{2} w}}{w(w+1)} E\left(z, f; \;s+\frac{w}{2}\right) \widetilde g\left(-\frac{n}{2} w   \right) \;dw
\\
& \hskip 30pt
 - \frac{\Lambda(2ns-2n, \; f\times \widetilde f\,)}{   \Lambda(1+2ns-2n, \; f\times \widetilde f\,)  } \cdot \frac{1}{2\pi i} \int\limits_{2-i\infty}^{2+i\infty} \frac{\beta^{-\frac{n}{2} w}}{w(w+1)} E\left(z,  f; \;2-s+\frac{w}{2}\right) \widetilde g\left(-\frac{n}{2} w   \right) \;dw,
\endalign$$
where $$E_0(z, f; s) := E(z, f; s) - A^{\frac{n}{2}} E\left(z, f; s-\frac12\right) + \frac{\Lambda(2ns-2n, f\times \widetilde f\,)}{\Lambda(1+2ns-2n, f\times \widetilde f\,)} E(z,  f; 2-s).$$

Now, we are assuming that $\max\limits_{1\le i < 2n} y_i \ge\left(|t|^{1+\epsilon}\right)^{n(2n-1)}$.
By proposition 10.5, for $\Re(s')$ fixed with $s'\in\Bbb C, \; |\Im(s')| \le |t|^{1+\epsilon}$, we have
$ND\left(E(z, f; s')\right) \ll |t|^{-N}.$
The conclusion follows.

\qed

\vskip 10pt\noindent
{\bf Corollary 10.8} {\it Let  $g \in  \Bbb C_c^\infty\left([1,2]\right)$.

\vskip 5pt {\bf (i)} If $h(z) \ge \beta $ and $\min\limits_{1\le i<2n} y_i \ge \frac{\sqrt{3}}{2}$, 
 then
$$\int\limits_0^\infty g\left(\frac{A}{\beta}\right) \widehat E_A^*(z, f;  s) \; \frac{dA}{A} \ll \beta^{\frac{n}{2}} h(z)^{\frac{n}{2}} \left( \frac{h(z}{\beta}  \right)^{-\frac12}.$$

\vskip 5pt {\bf (ii)} If $h(z) \le \beta $ and $\min\limits_{1\le i<2n} y_i \ge \frac{\sqrt{3}}{2}$, 
and $\max\limits_{1\le i < 2n} y_i \ge \left(t^{1+\epsilon}\right)^{n(2n-1)}$,  then
$$\int\limits_0^\infty g\left(\frac{A}{\beta}\right) \widehat E_A^*(z, f; s) \; \frac{dA}{A} \; \ll \; \beta' = \beta^{\frac{n}{2}}\left(y_n^n y_{n+1}^{n-1} \cdots y_{2n-1}\right)^{\frac{n}{2}}.$$}

{\bf Proof:} This follows directly from proposition 10.7. \qed

   \vskip 20pt
 \noindent
  {\bf \S 11 Upper bound for the integral $\Cal I$:} 
  
 \vskip 10pt
 
   The key idea for proving our main theorem 1.3 is in obtaining suitable upper and lower bounds
 for the integral $\Cal I$ which we now define.
 
 \vskip 10pt\noindent
 {\bf Definition 11.1 (The integral $\Cal I$)}
{\it   Let $t\in \Bbb R$ with $|t| \gg 1$ and set   $$\beta = |t|^{n^{10}}, \qquad \delta = \beta^{-1} .$$ Let $f$ be a cusp form on $GL(n)$ as in (4.1). Set $$\frak h^{2n} := GL(2n,\Bbb R)/(O(2n,\Bbb R)\cdot\Bbb R^\times).$$
   We now define the integral $\,\Cal I = \Cal I_{g,\psi}(t)$ for test functions
   $g \in \Bbb C_0^\infty([1,2])$ (with $\widetilde g(n/2) = 1$) and $\psi:[0,\infty]\to\Bbb C,$ where $g, \psi$ are non-negative. In addition, we require that for some fixed $a >0$, the Mellin transform $\widetilde{\psi}(\omega)$ is holomorphic in $-a\le \Re(w) \le a$, and satisfies $\widetilde{\psi}(\omega) \ll e^{-n|w|}$ in this strip. Let
     $$\align & \Cal I  \; := \; \big| L(1+2int, f\times \widetilde f\,)\big|^2 \hskip-10pt\int\limits_{P_{2n-1,1}(\Bbb Z)\backslash\frak h^{2n}} \left| \; \int\limits_0^\infty \widehat E_A^*(z, f; 1+it)\, g\left(  \frac{A}{\beta} \right) \frac{dA}{A} \; \right|^2 \hskip-2pt\cdot \psi\left( \frac{\text{\rm det}(z)}{\delta}  \right)  \, d^*z.  \endalign$$
     }
   
   \vskip 10pt\noindent
   {\bf Theorem 11.2 (Upper bound for $\Cal I$)} {\it  For $|t| \gg 1$ we have the upper bound
   $$\boxed{ \Cal I_{g,\psi}(t) \; \ll \; \big| L(1+2int, f\times \widetilde f\,)\big| \cdot \delta^{-\frac12} \beta^{\frac12+n} \big(\log |t|\big)^2\cdot \Big[1 +  \big| L(1+2int, f\times \widetilde f\,)\big|  \Big], }$$
   where the $\ll$-constant depends at most on $n$ and $f$.}
   
   \pagebreak
   
   \vskip 10pt {\bf Proof of theorem 11.2:} We begin with some standard computations involving Mellin inversion and Rankin-Selberg unfolding. 
   $$\align  \Cal I_{g,\psi}(t) & =  \big| L(1+2int, f\times \widetilde f\,)\big|^2 \int\limits_{P_{2n-1,1}(\Bbb Z)\backslash\frak h^{2n}} \left| \;\int\limits_0^\infty \widehat E_A^*(z, f; 1+it)\, g\left(  \frac{A}{\beta} \right) \frac{dA}{A} \;\right|^2\\
   &
   \hskip 180pt \cdot \left(\; \int\limits_{2-i\infty}^{2+i\infty} \widetilde\psi(-w) \left( \frac{\text{det}(z)}{\delta}   \right)^w \,dw \right)   d^*z \\
   &
   \\
   &
   = \;\frac{\big| L(1+2int, f\times \widetilde f\,)\big|^2}{2\pi i}\, \int\limits_{2-i\infty}^{2+i\infty} \widetilde\psi(-w)\, \delta^{-w}
   \\
   &
   \hskip 50pt
   \cdot
   \int\limits_{SL(2n,\Bbb Z)\backslash \frak h^{2n}} \left| \int_0^\infty g\left( \frac{A}{\beta}  \right)  \widehat E_A^*(z, f; 1+it) \;\frac{dA}{A}\right|^2\cdot   E_{P_{2n-1,1}}(z, w)\;  d^*z \; dw.
   \endalign$$

   Next, we shift the line of integration in the $w$-integral to $\Re(w) = \frac12.$ In doing so, we cross a pole of the maximal parabolic Eisenstein series $E_{P_{2n-1,1}}(z, w)$ at $w = 1.$ It follows that
  $$\align  \Cal I_{g,\psi}(t) & =   c\cdot \frac{\big| L(1+2int, f\times \widetilde f\,)\big|^2  \; \widetilde \psi(-1)}{ \delta} \hskip-8pt  \int\limits_{SL(2n,\Bbb Z)\backslash \frak h^{2n}} \left|  \int_0^\infty g\left( \frac{A}{\beta}  \right)  \widehat E_A^*(z, f; 1+it) \;\frac{dA}{A}\right|^2 d^*z\\
  &
  \\
  & \hskip 37pt + \frac{\big| L(1+2int, f\times \widetilde f\,)\big|^2}{2\pi i}\, \int\limits_{\frac12-i\infty}^{\frac12+i\infty} \widetilde\psi(-w)\, \delta^{-w}
   \\
   &
   \hskip 54pt
   \cdot
   \int\limits_{SL(2n,\Bbb Z)\backslash \frak h^{2n}} \left| \int_0^\infty g\left( \frac{A}{\beta}  \right)  \widehat E_A^*(z, f; 1+it) \;\frac{dA}{A}\right|^2\cdot   E_{P_{2n-1,1}}(z, w) \;d^*z\; dw\\
   &
   \\
   & := \Cal I_{g,\psi}^{(1)}(t)  \; + \; \Cal I_{g,\psi}^{(2)}(t).
  \endalign$$
   Here $c\ne 0$ is a constant.

   Next, we will obtain upper bounds for $\Cal I^{(1)}$ and  $ \Cal I^{(2)}$. To bound $\Cal I^{(1)}$ we  break it into 3 pieces: 
  $$c^{-1}\cdot\Cal I^{(1)} = \Cal I^{(1)}_1 + \Cal I^{(1)}_2 + \Cal I^{(1)}_3,$$
 where 
  $$\Cal I^{(1)}_1 =   \big| L(1+2int, f\times \widetilde f\,)\big|^2  \;\frac{ \widetilde \psi(-1)}{ \delta} \hskip-5pt \underset h(z) \ge \beta \to{ \int\limits_{SL(2n,\Bbb Z)\backslash \frak h^{2n}}} \left|  \int_0^\infty g\left( \frac{A}{\beta}  \right)  \widehat E_A^*(z, f; 1+it) \;\frac{dA}{A}\right|^2 d^*z,$$
  
  $$\Cal I^{(1)}_2 =   \big| L(1+2int, f\times \widetilde f\,)\big|^2  \; \frac{\widetilde \psi(-1)}{ \delta} \hskip-18pt \underset h(z) \le \beta \to{\underset \max\limits_{1\le i<2n} y_i  \;\ge\; \left( t^{1+\epsilon}  \right)^{n(2n-1)}\to{ \int\limits_{SL(2n,\Bbb Z)\backslash \frak h^{2n}}} }\left|  \int_0^\infty g\left( \frac{A}{\beta}  \right)  \widehat E_A^*(z, f; 1+it) \;\frac{dA}{A}\right|^2 d^*z,$$ 
  
$$\Cal I^{(1)}_3 =   \big| L(1+2int, f\times \widetilde f\,)\big|^2  \; \frac{\widetilde \psi(-1)}{ \delta} \hskip-18pt  \underset h(z) \le \beta \to{\underset \max\limits_{1\le i<2n} y_i \;\le\; \left( t^{1+\epsilon}  \right)^{n(2n-1)}  \to{ \int\limits_{SL(2n,\Bbb Z)\backslash \frak h^{2n}}}} \left|  \int_0^\infty g\left( \frac{A}{\beta}  \right)  \widehat E_A^*(z, f; 1+it) \;\frac{dA}{A}\right|^2 d^*z.$$   

By corollary 10.8 (i) and the fact that $h(z)^n\prod\limits_{\ell=1}^{2n-1} y_\ell^{-\ell(2n-\ell)} \ll 1$,  it follows that
$$\Cal I^{(1)}_1 \ll \delta^{-1} \beta^n \cdot \big|L(1+2int, \; f\times \widetilde f\,)\big|^2.$$
Similarly, by corollary 10.8 (ii)
$$\Cal I^{(1)}_2 \ll \delta^{-1} \beta^n (\log |t|)^2 \cdot\big|L(1+2int, \; f\times \widetilde f\,)\big|^2,$$
and by proposition 6.4, it follows that
$$\Cal I^{(1)}_3 \ll \delta^{-1} \beta^n  (\log |t|)^2 \cdot \big|L(1+2int, \; f\times \widetilde f\,)\big|.$$
\vskip 8pt
It remains to bound $\Cal I^{(2)}$ which we also split into 3 parts:
$$2\pi i\cdot \Cal I^{(2)} = \Cal I^{(2)}_1 + \Cal I^{(2)}_2 + \Cal I^{(2)}_3.$$
Here
$$\align &\Cal I^{(2)}_1 = \big| L(1+2int, f\times \widetilde f\,)\big|^2 \int\limits_{\frac12-i\infty}^{\frac12+i\infty} \widetilde\psi(-w) \, \delta^{-w} \hskip-10pt \underset h(z) \ge \beta\to{ \int\limits_{SL(2n,\Bbb Z)\backslash \frak h^{2n}}} \left|  \int_0^\infty g\left( \frac{A}{\beta}  \right)  \widehat E_A^*(z, f; 1+it) \;\frac{dA}{A}\right|^2\\
&\hskip 285pt\cdot E_{P_{2n-1,1}}(z, w) \; d^*z\; dw,\endalign$$

$$\align &\Cal I^{(2)}_2 = \big| L(1+2int, f\times \widetilde f\,)\big|^2 \int\limits_{\frac12-i\infty}^{\frac12+i\infty} \widetilde\psi(-w) \, \delta^{-w} \hskip-28pt \underset h(z) \le \beta \to{\underset \max\limits_{1\le i<2n} y_i  \;\ge\; \left( t^{1+\epsilon}  \right)^{n(2n-1)}\to{ \int\limits_{SL(2n,\Bbb Z)\backslash \frak h^{2n}}} } \hskip-20pt\left|  \int_0^\infty g\left( \frac{A}{\beta}  \right)  \widehat E_A^*(z, f; 1+it) \;\frac{dA}{A}\right|^2\\
&\hskip 282pt\cdot E_{P_{2n-1,1}}(z, w) \;  d^*z \;dw, \endalign$$

$$\align &\Cal I^{(2)}_3 = \big| L(1+2int, f\times \widetilde f\,)\big|^2 \int\limits_{\frac12-i\infty}^{\frac12+i\infty} \widetilde\psi(-w) \, \delta^{-w}  \hskip-33pt  \underset h(z) \le \beta \to{\underset \max\limits_{1\le i<2n} y_i \;\le\; \left( t^{1+\epsilon}  \right)^{n(2n-1)}  \to{ \int\limits_{SL(2n,\Bbb Z)\backslash \frak h^{2n}}}} \hskip-15pt\left|  \int_0^\infty g\left( \frac{A}{\beta}  \right)  \widehat E_A^*(z, f; 1+it) \;\frac{dA}{A}\right|^2\\
&\hskip 283pt\cdot E_{P_{2n-1,1}}(z, w) \;  d^*z \;dw. \endalign$$
By corollary 10.8 (i), proposition 3.4, and the fact that
$$\align & h(z)^n  \sum_{1\le k \le 2n}\Bigg[ \Big(y_1 y_2^2 \; \cdots \; y_{2n-k}^{2n-k}    \Big)^{\frac12}\Big(y_{2n-k+1}^{k-1} y_{2n-k+2}^{k-2} \;\;\cdots\;\; y_{2n-1}   \Big)^{ \frac12  }\\
 & + \Big(y_1 y_2^2 \; \cdots \; y_{2n-k}^{2n-k}    \Big)^{\frac{k-1}{2n}} \Big(y_{2n-k+1}^{k-1} y_{2n-k+2}^{k-2} \;\;\cdots\;\; y_{2n-1}   \Big)^{ \frac{2n-k+1}{2n}  }\Bigg]\cdot \prod_{\ell=1}^{2n-1} y_\ell^{-\ell(2n-\ell)} \;\ll \;h(z)^{\frac12},\endalign$$
 it follows that
$$\Cal I^{(2)}_1 \ll \delta^{-\frac12} \beta^{\frac12+n} (\log |t|) \cdot\big|  L(1+2int, f\times \widetilde f\,)\big|^2.$$
By corollary 10.8 (ii) and proposition 3.4 it follows that
$$\Cal I^{(2)}_2 \ll \delta^{-\frac12} \beta^{\frac12+n} (\log |t|)\cdot  \big|  L(1+2int, f\times \widetilde f\,)\big|^2.$$
 Next, by proposition 3.4, for $\frac{\sqrt{3}}{2} \le y_i\le |t|^{(1+\epsilon)n (2n-1)}$ $(1\le i < 2n)$, we obtain the bound
  $$\align  E_{P_{2n-1,1}}(z, w) & \ll \sum_{1\le k\le 2n} \bigg[\left(|t|^{1+\epsilon}\right)^{n(2n-1) \left(  \frac{(2n-k)(2n-k+1)}{2}  \;+\; \frac{(k-1)k}{2}\right)} \;\\
  & 
  \hskip 95pt + \;  \left(|t|^{1+\epsilon}\right)^{n(2n-1) \left(  \frac{(2n-k)(2n-k+1)(k-1)}{4n}  \;+\; \frac{(k-1)k(2n-k+1)}{4n}\right)} \bigg]\\
  &
  \\
  & \ll \; |t|^{100n^4} (\log |t|) \; \ll \; \beta^{\frac12}.
  \endalign$$
  Consequently, by proposition 6.4 it follows that
  $$\Cal I^{(2)}_3 \ll \delta^{-\frac12} \beta^{n+\frac12} (\log |t|)^2 \cdot \big|  L(1+2int, f\times \widetilde f\,)\big|.$$
  
  Finally, combining all the above bounds, we get
  $$ \Cal I_{g,\psi}(t) \ll \big|  L(1+2int, f\times \widetilde f\,)\big| \cdot \delta^{-\frac12} \beta^{n+\frac12} (\log | t|)^2 \cdot \Big[ 
  1 +  \big|  L(1+2int, f\times \widetilde f\,)\big| \Big].$$
  \qed

   \vskip 20pt
 \noindent
  {\bf \S 12 Lower bound for the integral $\Cal I$:} 
  
 \vskip 10pt
  The main aim of this section is to prove the following  lower bound for the integral $\Cal I$.
 \vskip 10pt\noindent
 {\bf Theorem 12.1(Lower bound for $\Cal I$)} {\it Assume that the cusp form $f$ for $SL(n, \Bbb Z)$ has Langlands parameters $(i\alpha_1, \ldots,i\alpha_n)$. Assume further that   $\widetilde\psi(w)$ vanishes (to order $n^4$)  at $w = i(\alpha_j-\alpha_k) \ne 0$ for $1\le j, \; k\le n$ and that $\widetilde\psi(0) = 1.$  Then, under the same assumptions as in definition 11.1 and theorem 11.2, we have
 $$\boxed{\Cal I_{g,\psi}(t) \; \gg \; \beta^n \delta^{-1}\big/(\log |t|),} $$
 where the $\gg$-constant depends at most on $n$ and $f$.}
 
 \vskip 10pt\noindent
 {\bf Remark:}  To show that it is possible for $\psi$ to vanish to high order at finitely many pure imaginary points and also satisfy the conditions specified in definition 11.1 (i.e., positivity and exponential decay of the Mellin transform in a strip) can be seen as follows. Let $\Psi$ be a function satisfying the conditions in definition 11.1.
   For $\lambda > 0$, define $\psi(y) := \Psi(\lambda y) + \Psi(y).$ Clearly $\psi$  satisfies the conditions specified in definition 11.1. Then
 $$\widetilde \psi(w) = \big(\lambda^{-w} + 1\big)\widetilde\Psi(w).$$
If we choose $\lambda = e^{\pi/\alpha}$ then it is clear that $\widetilde\psi(i\alpha) = 0.$ By iterating this procedure we can construct a test function $\psi$ having the properties required in theorem 12.1. For example, we may initially choose a test function of type $\Psi(y) = y^R e^{-y^{1/2n}}$ (for some large positive $R$) and then apply the above procedures. We may, therefore, take
$$\widetilde \psi(w) := \frac{\Gamma\big(2n(R+w)\big)}{2^{4(n^2-L)} \,\Gamma(2nR)}\underset \alpha_j\ne\alpha_k\to{ \prod_{1\le j, k\le n}}\left( e^{\frac{\pi w}{\alpha_j-\alpha_k}} +1  \right)^4. $$
where $L = \#\{j,k \mid \alpha_j = \alpha_k\}.$
  From now on we shall assume that $\widetilde\psi$ is of this form with $R \gg 1$ sufficiently large and independent of $n$ and $f$.
 
 \vskip 10pt
 We defer the proof of theorem 12.1 until later. A key ingredient of the proof is the following orthogonality condition stating that the degenerate part of the Fourier expansion of $F$ is orthogonal to the non-degenerate part. 
 
 \vskip 10pt\noindent
 {\bf Proposition 12.2} {\it  Suppose $F$ is an automorphic function for $SL(k,\Bbb Z)$ with $k\ge 2$ as in theorem 8.3. Define $D(F)$ and $ND(F)$ as in definition 8.4. Then
 $$\int\limits_{P_{k-1,1}(\Bbb Z)\backslash \frak h^k}\overline{D}(F) \cdot ND(F) \psi\left( \frac{\text{\rm det}(z)}{\delta}   \right) \, d^\star z = 0.$$
  }
  
  {\bf Proof of proposition 12.2:} We have
  $$\align & \overline{D}(F) \cdot ND(F) = \sum_{m_1\ne 0}\cdots \sum_{m_{k-1}\ne0} \; \sum_{\gamma \in U_{k-1}(\Bbb Z)\backslash SL(k-1,\Bbb Z)} \\
  &
 \hskip 20pt \cdot \overline{D}(F)\left( \pmatrix \gamma &\\& 1\endpmatrix z  \right) \int\limits_0^1 \cdots \int\limits_0^1  F\left(\pmatrix 1 & u_{1,2} & &\cdots & u_{1,k}\\
  & 1 & u_{2,3} & \cdots & u_{2k}\\
  & & \ddots &\ddots&\vdots\\
  & & & 1 & u_{k-1,k}\\
  & & && 1 \endpmatrix  \pmatrix \gamma &\\&1\endpmatrix z\right)\\
  & \hskip 230pt \cdot e^{-2\pi i(m_1u_{1,2} + \cdots + m_{k-1}u_{k-1,k})} \; d^\times u.\endalign$$
  Since 
  $$\bigcup_{\gamma \in U_{k-1}(\Bbb Z)\backslash SL(k-1,\Bbb Z)}  \pmatrix \gamma & \\ & 1 \endpmatrix \; P_{k-1,1} \backslash \frak h^k \; \equiv \; U_k(\Bbb Z)\backslash \frak h^k,$$
  we have
  $$\align & \int\limits_{ P_{k-1,1} \backslash \frak h^k}  \overline{D}(F) \cdot ND(F) \psi\left( \frac{\text{\rm det}(z)}{\delta}   \right) \, d^*z = \sum_{m_1\ne 0} \cdots \sum_{m_{k-1} \ne 0} \;\int\limits_{ U_k(\Bbb Z)\backslash \frak h^k } \overline{D}(F) \int\limits_0^1 \cdots  \int\limits_0^1 
  \\
  &
 \hskip 50pt \cdot F\left(\pmatrix 1 & u_{1,2} & &\cdots & u_{1,k}\\
  & 1 & u_{2,3} & \cdots & u_{2,k}\\
  & & \ddots &\ddots&\vdots\\
  & & & 1 & u_{k-1,k}\\
  & & && 1 \endpmatrix  z\right) \cdot e^{-2\pi i\big(m_1u_{1,2} + \cdots + m_{k-1}u_{k-1,k}\big)}\\
  &
  \hskip 260pt\cdot \psi\left( \frac{\text{\rm det}(z)}{\delta}   \right)  \, d^\times u \cdot  d^\star z
  \\
  &
  \\
  & = \int\limits_{y_1=0}^\infty \cdots \int\limits_{y_{k-1}=0}^\infty \;\int\limits_{x_{1,2}=0}^1\cdots \int\limits_{x_{k-1,k}=0}^1
   \overline{D}(F)   \cdot \int\limits_{u_{1,2}=0}^1 \cdots \int \limits_{u_{k-1,k}=0}^1\\
   &\hskip 50pt \cdot F\left(\pmatrix 1 & u_{1,2} & &\cdots & u_{1,k}\\
  & 1 & u_{2,3} & \cdots & u_{2k}\\
  & & \ddots &\ddots&\vdots\\
  & & & 1 & u_{k-1,k}\\
  & & && 1 \endpmatrix  z\right) \cdot e^{-2\pi i\big(m_1u_{1,2} + \cdots + m_{k-1}u_{k-1,k}\big)}\\
  & \hskip 260pt \cdot \psi\left( \frac{\text{\rm det}(z)}{\delta}   \right) d^\times u \;d^\star z.
  \endalign$$
  Consequently, it is enough to prove that
  $$\Cal J := \int\limits_{x_{1,2}=0}^1\cdots \int\limits_{x_{k-1,k}=0}^1
   \overline{D}(F)\; e^{2\pi i\big(m_1 x_{1,2} + \cdots + m_{k-1}x_{k-1,k}\big)} \;d^\times  x \; = 0.$$
   We write
   $$D(F) := D_1(F) + \cdots + D_{k-1}(F)$$
   and
   $$\Cal J := \Cal J_1 + \cdots + \Cal J_{k-1},$$
  where
  $$\align  \Cal J_1 & = \int\limits_{x_{1,2}=0}^1\cdots \int\limits_{x_{k-1,k}=0}^1
\; \int\limits_{u_{1,k}=0}^1\cdots \int\limits_{u_{k-1,k}=0}^1 F\left( \left(\smallmatrix  1 & 0 & \cdots & 0 & u_{1,k}\\
& 1 &\cdots & 0&u_{2,k}\\
& &\ddots &\vdots&\vdots \\
& & & 1&u_{k-1,k}\\
& & & & 1  \endsmallmatrix\right) 
\left(\smallmatrix 1 & x_{1,2} &\cdots & x_{1,k-1} & 0\\
 & 1 & \cdots & x_{2,k-1} & 0\\
 & & \ddots & \vdots & \vdots\\
 & & & 1 & 0\\
 & & & & 1 \endsmallmatrix\right)y \right)\\
 & \hskip 205pt\cdot e^{2\pi i\big(m_1x_{1,2} + \cdots + m_{k-1}x_{k-1,k}\big)} \; d^\times x \; d^\times u
 \\
 & = 0,\endalign$$
after computing the $x_{k-1,k}$-integral. 

For $1 \le \ell \le k-2,$ we define
  $$\align & \Cal J_{\ell +1} := \int\limits_{x_{1,2}=0}^1\cdots \int\limits_{x_{k-1,k}=0}^1  F\left( \left(\smallmatrix  1 & 0 &\cdots & 0 & u_{1,k-\ell} & \cdots & u_{1,k-1} & u_{1,k}\\ 
  & 1 &\cdots & 0 & u_{2,k-\ell} & \cdots & u_{2,k-1}& u_{2,k} \\
  & &  & \ddots &\vdots &\cdots & \vdots&\vdots
  \\
  & & & & & &1 & u_{k-1,k}\\
  &&&&&&&1 \endsmallmatrix\right) \pmatrix \gamma_{k-1} &\\&1\endpmatrix z \right)
  \\
  &\hskip 60pt\cdot \int\limits_{u_{1,k-l}=0}^1\cdots \int\limits_{u_{k-1,k}=0}^1 e^{2\pi i\big(-m_{k-1}' u_{k-1,k} - \cdots - m_{k-\ell}' u_{k-\ell,k-\ell+1}\big)} \; d^\times u\\
  &
 \hskip 200pt \cdot \; e^{2\pi i\big(m_{1} x_{1,2} + \cdots + m_{k-1} x_{k-1,k}\big)} \; d^\times x.\endalign$$
  Here
  $$\align
  &\pmatrix \gamma_{k-1} &\\&1\endpmatrix z \; \equiv \; \pmatrix    
  1 & x_{1,2}(\gamma_{k-1}) &x_{1,3}(\gamma_{k-1}) & \cdots &  x_{1,k}(\gamma_{k-1})\\
  &  1 & x_{2,3}(\gamma_{k-1}) & \cdots & x_{2,k}(\gamma_{k-1})\\
  & & \ddots & \cdots&\vdots\\
  & & & 1&x_{k-1,k}(\gamma_{k-1})\\
  & & & & 1\endpmatrix y(\gamma_{k-1}) \\
  & \hskip 290pt\Big(\hskip-8pt\mod{\Bbb Z_k \,O(k,\Bbb R)}\Big).
  \endalign$$
 
  Now $x_{1,k}, x_{2,k}, \ldots, x_{k-1,k}$ do not appear in $y(\gamma_{k-1}).$ Note also that
  $$x_{k-1,k}(\gamma_{k-1}) = a_{k-1,1} x_{1,k} + a_{k-1,2}x_{2,k} + \;\cdots\; + a_{k-1,k-1} x_{k-1,k},$$
  $$\align & \left(\smallmatrix  1 & 0 &\cdots & 0 & u_{1,k-\ell} & \cdots & u_{1,k-1} & u_{1,k}\\ 
  & 1 &\cdots & 0 & u_{2,k-\ell} & \cdots & u_{2,k-1}& u_{2,k} \\
  & &  & \ddots &\vdots &\cdots & \vdots&\vdots
  \\
  & & & & & &1 & u_{k-1,k}\\
  &&&&&&&1 \endsmallmatrix\right) 
  \left( \smallmatrix    
  1 & x_{1,2}(\gamma_{k-1}) & \;\cdots\; &x_{1,k-\ell}(\gamma_{k-1}) & \;\cdots\; &  x_{1,k}(\gamma_{k-1})\\
  &  1 &\cdots & x_{2,k-\ell}(\gamma_{k-1}) & \cdots & x_{2,k}(\gamma_{k-1})\\
  & & \ddots& & \vdots&\vdots\\
   & &&&&\\
  & & & & 1&x_{k-1,k}(\gamma_{k-1})\\
  & & && & 1\endsmallmatrix\right)
  \\
  & \hskip 40pt =  \left( \smallmatrix    
  1 & x_{1,2}(\gamma_{k-1}) &\; \cdots\; &x_{1,k-\ell}(\gamma_{k-1})+u_{1,k-\ell} & \;\cdots\; &  x_{1,k}(\gamma_{k-1})+u_{1,k-\ell}x_{k-\ell,k}(\gamma_{k-1}) + \; \cdots \; + u_{1,k}\\
  &  1 & \cdots & x_{2,k-\ell}(\gamma_{k-1})+u_{2,k-\ell} & \cdots & x_{2,k}(\gamma_{k-1})+u_{2,k-\ell}x_{k-\ell,k}(\gamma_{k-1}) + \; \cdots \; + u_{2,k}\\
  & & \ddots &\vdots & \vdots&\vdots\\
  & &&&&\\
  & & & & 1&x_{k-1,k}(\gamma_{k-1})+u_{k-1,k}\\
  & & & & &1\endsmallmatrix\right).\endalign$$
 Now $x_{1,k}, x_{2,k}, \ldots, x_{k-1,k}$ only appear in $x_{1,k}(\gamma_{k-1}), \ldots ,x_{k-1,k}(\gamma_{k-1})$. If we make the change  variables
  $$\align & x_{1,k}(\gamma_{k-1})\; +\; u_{1,k-\ell} x_{k-\ell,k}(\gamma_{k-1}) \;+ \;\;\cdots\;\; +\; u_{1,k} \;\; \longrightarrow \;\; u_{1,k}
  \\
  &
   x_{2,k}(\gamma_{k-1})\; +\; u_{2,k-\ell} x_{k-\ell,k}(\gamma_{k-1}) \;+ \;\;\cdots\;\; +\; u_{2,k} \;\; \longrightarrow \;\; u_{2,k}
   \\
   &
   \hskip 20pt\vdots
   \\
   &
   x_{k-1,k}(\gamma_{k-1})\; + \; u_{k-1,k} \;\; \longrightarrow \;\; u_{k-1,k},\endalign
  $$
  and compute the $x_{1,k}, x_{2,k}, \ldots, x_{k-1,k}$ integral, it follows that   $$a_{k-1,1} = \; \cdots \; = a_{k-1,k-2} = 0, \qquad a_{k-1,k-1} = 1,$$
  which implies that
  $$ \gamma_{k-1} \; \in \; \widetilde P_{k-1,1}.$$
  
  To complete the proof, we use  induction. Assume 
  $$\gamma_{k-1} \in \widetilde P_{k-1,\ell} \backslash \widetilde P_{k-1,i-1}$$ for $2\le i \le \ell.$ We will show that $\gamma_{k-1} \in \widetilde P_{k-1,i}.$
  
  Since
  $$\widetilde P_{k-1,i-1} = \left\{ \pmatrix \gamma_{k-i} & *\\
  0 & U_{i-1}\endpmatrix  \;\;\bigg | \;\; \gamma_{k-i} \in SL(k-i,\Bbb Z)\right\},$$
  we can assume $\gamma_{k-1}$ is of the form $\pmatrix \gamma_{k-i} & 0\\0 & U_{i-1}\endpmatrix,$ by multiplying a suitable matrix in $\widetilde P_{k-1,\ell}$ on the left.
  
  Let
  
  $$\gamma_{k-i} = \pmatrix a_{1,1} & \cdots & a_{1,k-i}\\
  \vdots &\cdots & \vdots\\
  a_{k-i,1} & \cdots & a_{k-i,k-i}\endpmatrix .$$
  Then
  $$\align &  \pmatrix \gamma_{k-i}&0\\0&1\endpmatrix  \pmatrix 1 & x_{1,2} & x_{1,3} & \cdots & x_{_{1,k-i+1}}\\
   & 1& x_{2,3} & \cdots & x_{_{2,k-i+1}}\\
   & & \ddots &\cdots &\vdots\\
   & & & 1 & x_{_{k-i,k-i+1}}\\
   &&&& 1\endpmatrix y
   \\
   &
   \\
   &
   \hskip 8pt \equiv \pmatrix 1 & x_{1,2}(\gamma_{_{k-i}}) & x_{1,3}(\gamma_{_{k-i}}) & \cdots & & x_{_{1,k-i+1}}(\gamma_{_{k-i}})\\
   & 1 & x_{2,3}(\gamma_{_{k-i}}) & \cdots & & x_{_{2,k-i+1}}(\gamma_{_{k-i}})\\
   & & \ddots & \cdots & & \vdots\\
   & & &  &1 & x_{_{k-i,k-i+1}}(\gamma_{_{k-i}})\\
   &&&&& 1\endpmatrix y(\gamma_{_{k-i}})\quad \Big(\text{mod} Z_k\, O(k,\Bbb R)   \Big).\endalign$$
  Here
  $$x_{k-i,k-i+1}(\gamma_{k-i}) =  a_{k-i,1} x_{1, k-i+1} \; + \; 
  a_{k-i,2} x_{1, k-i+1} \; + \; \cdots \; + a_{k-i,k-i} x_{k-i,k-i+1}$$
  does not appear in $ y(\gamma_{k-i})
$.

Next, change variables
$$\align & u_{1,k-i+1} + \; \cdots \; + x_{1,k-i+1} (\gamma_{k-i})
\; \longrightarrow \; u_{1,k-i+1}\\
 & \phantom{xxx} \vdots \\
 & u_{k-i, k-i+1} \; + \; x_{k-i, k-i+1}(\gamma_{k-i}) \; \longrightarrow \; u_{k-i, k-i+1}.\endalign$$
 Computing the
 $x_{_{1,k-i+1}},  \quad\ldots\quad x_{_{k-i,k-i+1}} 
$
 integrals,  it follows that
$$a_{k-i,1} = \cdots = a_{k-i,k-i-1} = 0, \qquad a_{k-i, k-i} = 1,
$$
which implies that
$$\gamma_{k-i} \in \widetilde P_{k-1, i}.$$
Hence
$$\align & \Cal J_{\ell +1} := \int\limits_{x_{1,2}=0}^1\cdots \int\limits_{x_{k-1,k}=0}^1  F\left( \left(\smallmatrix  1 & 0 &\cdots & 0 & u_{1,k-\ell} & \cdots & u_{1,k-1} & u_{1,k}\\ 
  & 1 &\cdots & 0 & u_{2,k-\ell} & \cdots & u_{2,k-1}& u_{2,k} \\
  & &  & \ddots &\vdots &\cdots & \vdots&\vdots
  \\
  & & & & & &1 & u_{k-1,k}\\
  &&&&&&&1 \endsmallmatrix\right)  z \right)
  \\
  &\hskip 60pt\cdot \int\limits_{u_{1,k-l}=0}^1\cdots \int\limits_{u_{k-1,k}=0}^1 e^{2\pi i\big(-m_{k-1}' u_{k-1,k} - \cdots - m_{k-\ell}' u_{k-\ell,k-\ell+1}\big)} \; d^\times u\\
  &
 \hskip 200pt \cdot \; e^{2\pi i\big(m_{1} x_{1,2} + \cdots + m_{k-1} x_{k-1,k}\big)} \; d^\times x.\endalign$$

Here $1 \le \ell \le k-2$ and
$$\align  & \matrix \hskip-105pt \underbrace{{\scriptscriptstyle k-\ell-1 \;\text{column}}} \\
 \left(\smallmatrix 1 & 0 & \cdots & 0 & u_{1,k-\ell} & \cdots & u_{1,k-1} &u_{1,k}\\
& 1 & \cdots & 0 & u_{2,k-\ell} & \cdots & u_{2,k-1} & u_{2,k}\\
& &\ddots & \vdots &\vdots & \cdots & \vdots&\vdots\\
& & & 1 & u_{_{k-\ell-1,k-\ell}} & \cdots & u_{_{k-\ell-1,k-1}} & u_{_{k-\ell-1,k}}\\
&&&&&\ddots&\vdots&\vdots\\
&&&&&&1&u_{k-1,k}\\
&&&&&&&1
\endsmallmatrix\right) 
\endmatrix
\left(\smallmatrix 1& x_{1,2} &\cdots & x_{1,k-\ell} & \cdots & x_{1,k}\\
& 1 & \cdots &x_{2,k-\ell} & \cdots & x_{2,k}\\
& & \ddots & \vdots &  &\vdots\\
& & &\;\;1 \;\; x_{_{k-\ell-1,k-\ell}} & \cdots & x_{_{k-\ell-1,k}}\\
&&&\qquad \ddots&&\vdots\\
& &  & & 1& x_{k-1,k}\\
& & & & & 1 \endsmallmatrix\right)\\
&
\\
&\hskip 3pt = \left(\smallmatrix 1 & x_{1,2} & \;\; \cdots \;\; & x_{_{1,k-\ell-1}} & x_{_{1,k-\ell}}+u_{_{1,k-\ell}} & \;\;\; \cdots \;\;\; & u_{1,k} +\; \cdots\; + x_{1,k}\\ 
& \;\;\;\ddots & & \vdots & \vdots & &\vdots\\
& & \;\;  \; &  & && \\
& & & & &  &\\
& &  &1 &  x_{_{k-\ell,k-\ell}}+u_{_{k-\ell-1,k-\ell}} &\cdots& x_{_{k-\ell-1,k}}+u_{_{k-\ell-1,k-\ell}}x_{_{k-\ell,k}} + \; \cdots \; + u_{_{k-\ell-1,k}}
\\
& & & &\ddots & &\vdots\\
&&&&&&\\
& & & & && 1\endsmallmatrix\right).\endalign$$
The change of variables
$$\align &  x_{k-\ell-1,k-\ell} +u_{k-\ell-1,k-\ell} \; \longrightarrow \; u_{k-\ell-1,k-\ell}\\
& x_{k-\ell-1,k-\ell+1} + u_{k-\ell-1,k-\ell}\, x_{k-\ell,k-\ell+1} + u_{k-\ell-1,k-\ell+1} \; \longrightarrow \; u_{k-\ell-1,k-\ell+1}\\
&\phantom{xxxx}\vdots \\
&x_{k-\ell-1,k} + u_{k-\ell-1,k-\ell} \,x_{k-\ell,k} +\;\;\cdots\;\; + u_{k-\ell-1,k} \; \longrightarrow \; u_{k-\ell-1,k} \endalign$$
and the computation of the $x_{k-\ell-1,k-\ell}$-integral shows that $\Cal J_{\ell+1} = 0.$

\qed

  \vskip 10pt
  \noindent
  {\bf Corollary 12.3} {\it We have for any $N > 1$ that
  $$\align & \Cal I_{g,\psi}(t) \; \gg \; \big |L(1+2int, \; f\times \widetilde f\,)\big |^2\; \int\limits_0^\infty \cdots \int\limits_0^\infty \sum_{N \le p \le 2N}\\
  & \hskip 60pt \cdot \Bigg|\int\limits_0^1\cdots \int\limits_0^1 \;\int \limits_0^\infty \widehat{E}_A^*
  \left(  \pmatrix 1 & u_{1,2} & \cdots & u_{1,2n}\\
  & 1 & \cdots & u_{2,2n}\\
  & &\ddots &\vdots\\
  &&&1\endpmatrix y, \; f; \; 1+it \right) g\left( \frac{A}{\beta}   \right)  \;\frac{dA}{A}
  \\
  &
  \hskip 110pt\cdot e^{2\pi i\big(-pu_{2n-1,2n}\; - \; u_{2n-2,2n-1} - \;\; \cdots \;\; - u_{1,2}\big)}\, d^*u\Bigg|^2 \; \psi\left(  \frac{\text{\rm Det}(y)}{\delta} \right)\; d^*y,
  \endalign$$
  where the sum above goes over primes $N\le p \le 2N$. }
  
  \vskip 10pt
  {\bf Proof:} This follows  from proposition 12.2 and definition 11.1 after noting that only the absolute value squared of the non-degenarate terms (see definition 8.4) contribute. The lower rank average in the Fourier expansion is gone because of the unfolding in the beginning of the proof of proposition 12.2. For the lower bound we only consider the contributions of the $(1,\ldots, 1,p)$ Fourier coefficients with $N\le p \le 2N$. Another key step in the proof of theorem 12.1 is the choice of $N$. \qed
  
  \vskip 10pt
   The following lemma is due to \cite{Stade, 2002}.
   \vskip 10pt
  \noindent
  {\bf Lemma 12.4} {\it Fix an integer $m\ge 2$, a vector $\nu \in\Bbb C^{m-1},$ and the associated Whittaker function $W_{m,\nu}$.  Let $(\beta_1, \ldots,\beta_m)$
 be the Langlands parameters associated to $(m,\nu)$ as defined in (2.3). Then
 $$\align & \int\limits_0^\infty\cdots \int\limits_0^\infty\; \Big | W_{m,\nu}(y)  \Big |^2 \left( \text{\rm Det} \,y)  \right)^w \; d^*y \; = \; \frac{\pi^{-\frac{(m-1)m}{2}w} \cdot  2^{\frac{(m-1)m(m+1)}{6}}    }{\Gamma\left(\frac{mw}{2}\right)} \\
 &
 \hskip 80pt\cdot \Bigg| \prod_{j=1}^{m-1} \prod_{j\le k\le m-1}
  \cdot \pi^{-\frac12-\frac12\left(\beta_{m-k} - \beta_{m-k+j}    \right)}\;
 \Gamma\left( \frac{ 1+\beta_{m-k}-\beta_{m-k+j} }{2}  \right)\Bigg |^{-2}\\
 &
 \hskip 270pt
 \cdot\prod_{j=1}^m\prod_{k=1}^m \Gamma\left(  \frac{w+\beta_j+\overline{\beta_k}}{2}  \right).\endalign$$
 }
 
 \vskip 10pt
 Next, we use corollary 12.3,  lemma 12.4 and results from sieve theory to get a lower bound for $\Cal I_{g,\psi}.$ 
 \vskip 10pt
 {\bf Proof of theorem 12.1:} By corollary 12.3, we need to consider the contribution from different parts of (8.9).
  The main contribution to the lower bound will come from the sum (over primes $p$) of the $(1,\ldots, 1, p)$ Fourier coefficients of $E(z, f; \frac12+it)$.  The Fourier coefficient for an individual prime $p$  (denoted  $a_1(1,\ldots,1,p)$) is given by (see \cite{Go}, Proposition 10.9.3, \;\cite{Shahidi}, 2010)

   $$\align a_1(1,\ldots, 1,p) \; & := \;  \int\limits_0^1\cdots \int\limits_0^1 \; E
  \left(  \pmatrix 1 & u_{1,2} & \cdots & u_{1,2n}\\
  & 1 & \cdots & u_{2,2n}\\
  & &\ddots &\vdots\\
  &&&1\endpmatrix y, \; f; \; \frac12+it \right)\\
  & \hskip 120pt\cdot e^{2\pi i\big(-pu_{2n-1,2n}\; - \; u_{2n-2,2n-1} - \;\; \cdots \;\; - u_{1,2}\big)}\, d^*u\\
  &
  \\
  & =\; \frac{c_n\, \overline{\lambda(p)\;\eta_{nit}(p)}}{p^{\frac{2n-1}{2}}} \;W_{2n,\nu}(My) \cdot\frac{1}{L(1+2int, \;f\times\widetilde  f\,)},
   \endalign$$
   where $c_n$ is an absolute constant depending only on $n$, 
 $$\eta_s(n) := \underset a,b \ge 1\to{\sum_{ab = n}} \left(\frac{a}{b}\right)^s, \qquad\qquad M ={\scriptstyle \left( \smallmatrix \displaystyle p & & &\\
 & \displaystyle1 & &\\
 & &\displaystyle\ddots &
\\
&&& \displaystyle1\endsmallmatrix\right)},$$
 and $\lambda(p)$ is the $p$-th Hecke eigenvalue of $f$. 
 The Langlands parameters of $E(z, f; \frac12+it)$ are
 $$\Big(i(nt+\alpha_1), \;\;\ldots\;\;,i(nt+\alpha_n), \; i(-nt+\alpha_1), \;\;\ldots \;\;,i(-nt+\alpha_n)\Big),$$
 where $(i\alpha_1, \;\;\ldots\;\;,i\alpha_n)$ are the Langlands parameters of $f$. For simplicity, we assume $\Re(\alpha_j) = 0$ for $1 \le j\le n.$
 
 Let
 $$\align I_1^*(p) &  := \int\limits_0^\infty \cdots  \int\limits_0^\infty\; \left| \,\int\limits_0^\infty A^{\frac{n}{2}} \cdot a_1(1,\ldots,1,p) \; g\left( \frac{A}{\beta}  \right) \;\frac{dA}{A}\right|^2 \psi\left( \frac{\text{\rm Det (z)}}{\delta}  \right) \; d^*y\\
 &
 \\
 & = \beta^n |\lambda(p)|^2\, |\eta_{nit}(p)|^2 \cdot\frac{|c_n|^2}{\left|L(1+2int, \;f\times \widetilde f\,)\right|^2} \cdot I_{1,1}^*(p),
 \endalign$$
 where
 
 $$\align & I_{1,1}^*(p)  = \frac{1}{2\pi i} \int\limits_{1-i\infty}^{1+i\infty}  \widetilde\psi(-w)\; (p\delta)^{-w} \int\limits_0^\infty \cdots \int\limits_0^\infty \left| W_{2n,\nu}(y)  \right|^2 \; \left(  \text{\rm Det}\, z \right)^w \; d^*y\; dw
 \tag 12.5 \\
 &\hskip3pt
 = \frac{\pi^{2n^2} 2^{\frac{(2n-1)2n(2n+1)}{6}}}{2\pi i}\prod_{j=1}^n \prod_{k=1}^n  \left| \Gamma\left( \frac{1+i(\alpha_k-\alpha_j)}{2}  \right)  
 \Gamma\left( \frac{1+ 2int +i(\alpha_k-\alpha_j}{2}  \right) \right|^{-2}
 \\
 &\hskip1pt
 \cdot\int\limits_{1-i\infty}^{1+i\infty} \frac{ \widetilde\psi(-w) (p\delta)^{-w}}
 {\pi^{\frac{2n(2n-1)w  }{ 2 }} }\;
 \frac{ \prod\limits_{j=1}^n \prod\limits_{k=1}^n  
  \Gamma\hskip-1pt\left( \frac{ w+i(\alpha_k-\alpha_j )}{2}  \right)^2  
  \Gamma\hskip-1pt\left( \frac{ w+2int+i(\alpha_k-\alpha_j }{2}  \right) 
   \Gamma\hskip-1pt\left( \frac{ w-2int+i(\alpha_k-\alpha_j }{2}  \right)          }{  \Gamma(nw)  }dw.\endalign$$

   Under the assumptions on $\widetilde\psi$, we can assume the above integrand is analytic at the points $w = i(\alpha_j-\alpha_k) \ne 0$ for $1\le j, \; k\le n.$ Shifting the line of integration in the $w$-integral above to $\Re(w) = -\frac{1}{\log\log R}$,  we pick up residues from poles at
   $$w = 0, \;\quad 2int+i(\alpha_j-\alpha_k), \;\quad -2int+i(\alpha_j-\alpha_k), \;\quad (1\le j, k \le n).$$
   
     At $w = 0$ there will be a pole of order $2L-1$ where
     $$L = \#\{(j,k) \mid \alpha_j = \alpha_k, \; 1 \le j,k \le n\}.$$
   The residue at $w = 0$ is
   
   $$\align & \underset w=0\to{\text{Res}}\;\Bigg[\frac{ \widetilde\psi(-w) (p\delta)^{-w}}
 {\pi^{\frac{2n(2n-1)w  }{ 2 }} }\; \frac{\Gamma\left(  \frac{w}{2} \right)^{2L}
}{\Gamma(nw)} \underset \alpha_k\ne\alpha_j\to{\prod_{1 \le j < k \le n} } \Gamma\hskip-1pt\left( \frac{ w+i(\alpha_k-\alpha_j )}{2}  \right)^2 \Gamma\hskip-1pt\left( \frac{ w-i(\alpha_k-\alpha_j )}{2}  \right)^2 \tag 12.6\\
   & 
\hskip 90pt \cdot \prod_{j=1}^n\prod_{k=1}^N   \Gamma\hskip-1pt\left( \frac{ w+i(2nt+\alpha_k-\alpha_j) }{2}  \right) 
    \Gamma\hskip-1pt\left( \frac{ w-i(2nt-\alpha_k+\alpha_j) }{2}  \right) \Bigg].
  \endalign$$
  
   Let $z = x+iy.$  Now, by Stirling's asymptotic expansion \cite{Whittaker-Watson, \S 12.33, 1927} we have for  $|z|\to\infty$ and $|\arg z| < \pi$ the asymptotic expansion
     $$\Gamma(z) = \sqrt{\frac{2\pi}{z}}\, \left( \frac{z}{e}  \right)^{z} \left( 1 + \frac{1}{12z} + \frac{1}{288 z^2} - \frac{139}{51840 z^3} \; + \; \Cal O\left( \frac{1}{|z|^4}  \right)  \right).$$
   It follows that for $-\frac12 \le x \le \frac12$ and $|y|\to \infty$ that
    $$\align & \Gamma\left(\frac{x+iy}{2}\right)  \Gamma\left(\frac{x-iy}{2}\right) = 2\pi e^{-x} \left(\frac{x^2+y^2}{4}\right)^{\frac{x-1}{2}} e^{-|y| \Big(\arctan\left(\frac{|y|}{x}\right)+\delta(x)\Big)}\\
   &
   \hskip 150pt\cdot \left | 1 + \frac{1}{12z} + \frac{1}{288 z^2} - \frac{139}{51840 z^3} \; + \; \Cal O\left( \frac{1}{|z|^4}  \right)   \right|^2,\endalign$$  
   where $$\delta(x) = \cases 0 &\text{if}\; x>0,\\ \pi &\text{if}\; x < 0.\endcases$$
 Such asymptotic expansions can also be derived for derivatives of the Gamma function. Consequently, for any fixed integer $\ell \ge 0$, if we take the $\ell^{th}$ derivative in $x$, we get for $-\frac12 \le x \le \frac12$ and $y\to\infty$, the asymptotic formula:   
  
   $$\align
   &
   \left(\frac{\partial}{\partial x }\right)^\ell \;\Gamma\left(\frac{x+iy}{2}\right)  \Gamma\left(\frac{x-iy}{2}\right) = \Bigg[\left(\frac{\partial}{\partial x }\right)^\ell 2\pi e^{-x} \left(\frac{x^2+y^2}{4}\right)^{\frac{x-1}{2}} e^{-|y|\left(\arctan\left(\frac{|y|}{x}\right)+\delta(x)\right)}\Bigg]  \tag{12.7}\\
   &\hskip 310pt \cdot \left(1+\Cal O\left(\frac{1}{|y|}\right)   \right),
  \\
 &\\
  & \left(\frac{\partial}{\partial x }\right)^\ell \;\Gamma\left(\frac{x+iy}{2}\right)  \Gamma\left(\frac{x-iy}{2}\right) \Bigg|_{x=0} \hskip-3pt=\;   \big(\log |y|/2  \big )^\ell \;\frac{4\pi e^{-\frac{\pi |y|}{2}}}{|y|}\left(1+\Cal O\left(\frac{1}{|y|}\right)   \right) \\
  &
  \\
  &
  \hskip 161pt
  =   \big(\log |y|/2   \big )^\ell \; \Gamma\left( \frac{iy}{2}  \right)  \Gamma\left( \frac{-iy}{2}  \right)\left(1+\Cal O\left(\frac{1}{|y|}\right)   \right). \endalign$$ 
  \vskip 5pt
  It follows from (12.7) that if we choose 
  $$N\delta = \prod_{j=1}^n \prod_{k=1}^n \exp\left(\log\left (\frac{|2nt+\alpha_k -\alpha_j|}{2}\right)\right)   \sim \left(n |t|\right)^{n^2}, $$   then, for any integer $\ell \ge 1$, we have
  $$\align
  & \left(\frac{\partial}{\partial w }\right)^\ell \Bigg[ (N\delta)^{-w}  \prod\limits_{j=1}^n \prod\limits_{k=1}^n  \Gamma\hskip-1pt\left( \frac{ w+i(2nt+\alpha_k-\alpha_j )}{2}  \right)  \Gamma\hskip-1pt\left( \frac{ w-i(2nt-\alpha_k+\alpha_j )}{2}  \right)\Bigg]_{w=0}\\
   &
  \\
   & \hskip 30pt =\sum_{m=0}^\ell \left(\matrix \ell \\
   m\endmatrix \right) \big(-\log(N\delta)\big)^m\\
   &\hskip 38pt\cdot \left(\frac{\partial}{\partial w }\right)^{\ell-m} \Bigg[\prod\limits_{j=1}^n \prod\limits_{k=1}^n  \Gamma\hskip-1pt\left( \frac{ w+i(2nt+\alpha_k-\alpha_j )}{2}  \right)  \Gamma\hskip-1pt\left( \frac{ w-i(2nt-\alpha_k+\alpha_j )}{2}  \right)\Bigg]_{w=0}\\
   &
   \\
   & \hskip 30pt =\sum_{m=0}^\ell \left(\matrix \ell \\
  m\endmatrix \right) \big(-\log(N\delta)\big)^m \cdot \left(\log(N\delta)^{\ell-m} \; + \; \Cal O\left( |t|^{-1} \right)\right)\\
   &\hskip 85pt\cdot \prod\limits_{j=1}^n \prod\limits_{k=1}^n  \Gamma\hskip-1pt\left( \frac{ w+i(2nt+\alpha_k-\alpha_j )}{2}  \right)  \Gamma\hskip-1pt\left( \frac{ w-i(2nt-\alpha_k+\alpha_j )}{2}  \right)\Bigg|_{w=0}   \\
   &
   \\
   &\hskip 50pt \ll \;\frac{(\log |t|)^\ell}{|t|} \cdot \prod_{j=1}^n\prod_{k=1}^n \left|\Gamma\left(\frac{i(2nt+\alpha_k-\alpha_j) }{2}    \right)\right|^2.
    \endalign $$
  \vskip 5pt  
  It now follows from the above computation that for $p\delta = N\delta\lambda_p$ with $1 \le \lambda_p \le 2$ and any integer $\ell \ge 0$ 
  
   $$\align
  &\left(\frac{\partial}{\partial w }\right)^\ell \Bigg[ (p\delta)^{-w}  \prod\limits_{j=1}^n \prod\limits_{k=1}^n  \Gamma\hskip-1pt\left( \frac{ w+i(2nt+\alpha_k-\alpha_j )}{2}  \right)  \Gamma\hskip-1pt\left( \frac{ w-i(2nt-\alpha_k+\alpha_j )}{2}  \right)\Bigg]_{w=0}\tag{12.8}\\
  &
  \\
  &
  \hskip 50pt = \left((-1)^\ell|\log\lambda_p|^\ell + \Cal O\left( \frac{(\log |t|)^\ell}{|t|}  \right)\right)\cdot  \prod_{j=1}^n\prod_{k=1}^n \left|\Gamma\left(\frac{i(2nt+\alpha_k-\alpha_j) }{2}    \right)\right|^2
  \\
  &
  \\
  &
  \hskip 50pt \le \left(|\log 2|^\ell + \Cal O\left( \frac{(\log |t|)^\ell}{|t|}  \right)\right)\cdot  \prod_{j=1}^n\prod_{k=1}^n \left|\Gamma\left(\frac{i(2nt+\alpha_k-\alpha_j )}{2}    \right)\right|^2.
\endalign$$

  \vskip 30pt
  To complete the determination of the residue at $w = 0$ in (12.6), we need to compute
  $$\align &  \left(\frac{\partial}{\partial w }\right)^{\ell} \left[\frac{ \widetilde\psi(-w)}
 {\pi^{\frac{2n(2n-1)w  }{ 2 }} }\; \underset \alpha_k\ne\alpha_j\to{\prod_{1 \le j < k \le n} } \Gamma\hskip-1pt\left( \frac{ w+i(\alpha_k-\alpha_j )}{2}  \right)^2  \Gamma\hskip-1pt\left( \frac{ w-i(\alpha_k-\alpha_j )}{2}  \right)^2\right] 
  \endalign$$
for all $0 \le \ell \le 2L-1.$

Recall that 
$$\widetilde \psi(-w) = \frac{\Gamma\big(2n(R-w)\big)}{2^{4(n^2-L)} \,\Gamma(2nR)}\underset \alpha_j\ne\alpha_k\to{ \prod_{1\le j, k\le n}}\left( e^{\frac{\pi w}{\alpha_j-\alpha_k}} +1  \right)^4. $$
If we assume that $n$ and $\alpha_j \;(j=1,2\ldots,n)$ are fixed and $R$ is sufficiently large then by Stirling's asymptotic formula we see that
$$ \left(\frac{\partial}{\partial w }\right)^{\ell} \widetilde \psi(-w)\; \bigg |_{w=0}  \;\sim \;  \big(-2n\log(2nR)\big)^\ell                 $$
as $R\to \infty.$ 
Consequently

$$\align &  \left(\frac{\partial}{\partial w }\right)^{\ell} \left[\frac{ \widetilde\psi(-w)}
 {\pi^{\frac{2n(2n-1)w  }{ 2 }} }\; \underset \alpha_k\ne\alpha_j\to{\prod_{1 \le j < k \le n} } \Gamma\hskip-1pt\left( \frac{ w+i(\alpha_k-\alpha_j )}{2}  \right)^2  \Gamma\hskip-1pt\left( \frac{ w-i(\alpha_k-\alpha_j )}{2}  \right)^2\right]_{w=0}\tag{12.9}\\
 &
 \\
 & \hskip 60pt \sim  \;   \big(-2n\log(2nR)\big)^\ell \cdot   \underset \alpha_k\ne\alpha_j\to{\prod_{1 \le j < k \le n} } \left|\Gamma\hskip-1pt\left( \frac{ i(\alpha_k-\alpha_j )}{2}  \right)\right|^4    \endalign$$
as $R\to\infty$, for all $0 \le \ell \le 2L-1.$

\vskip 5pt

The estimations (12.8) and (12.9) show that the main contribution to the residue at $w = 0$ in (12.6) comes from the $(2L-2)^{th}$ derivative of $\widetilde \psi(w)$ at $w=0.$
  It now follows from (12.8) and(12.9) that the residue at $w = 0$ in (12.6) is asymptotic to 
   $$\align
   & \sim\; n\cdot 2^{2L}\frac{ \big(2n\log(2nR)\big)^{2L-2} }{(2L-2)!}\underset \alpha_k\ne\alpha_j\to{\prod_{1 \le j < k \le n} } \left|\Gamma\hskip-1pt\left( \frac{ i(\alpha_k-\alpha_j )}{2}  \right)\right|^4 \; \prod_{j=1}^n\prod_{k=1}^n \left|\Gamma\left(\frac{i(2nt+\alpha_k-\alpha_j }{2}    \right)\right|^2,
   \endalign$$
for $R$ sufficiently large and $|t|\to\infty.$

 \vskip 8pt
   The other residues at $\pm 2int + i(\alpha_j-\alpha_k)$ will be bounded by
   $$\align & \ll \; \left|\widetilde\psi^{(j)}\Big(-\big(\pm 2int + i(\alpha_j-\alpha_k\big)  \Big)\right|\cdot |t|^{-n^2}\\
   & \ll \; |t|^{-B}\endalign$$  
   for arbitrary $B > 1,$  because of the rapid decay of $\widetilde\psi$ and its derivatives.

     \vskip 8pt
   On the line $\Re(w) = -\frac{1}{\log\log R}$  we estimate (with Stirling's asymptotic formula)  the integral
   $$
   \align
   & \int\limits_{- \frac{1}{\log\log R} -i\infty}^{ -\frac{1}{\log\log R} +i\infty} \frac{ \widetilde\psi(-w) (p\delta)^{-w}}
 {\pi^{\frac{2n(2n-1)w  }{ 2 }} }\; \frac{\Gamma\left(\frac{w}{2}  \right)^{2L}}{\Gamma(nw)} \; \underset \alpha_j\ne\alpha_k\to{ \prod\limits_{j=1}^n \prod\limits_{k=1}^n }   \Gamma\hskip-1pt\left( \frac{ w+i(\alpha_k-\alpha_j )}{2}  \right)^2 \tag{12.10}\\
 &\hskip 90pt \cdot
 \prod\limits_{j=1}^n \prod\limits_{k=1}^n  
 \hskip-2pt 
  \Gamma\hskip-1pt\left( \frac{ w+2int+i(\alpha_k-\alpha_j }{2}  \right) 
   \Gamma\hskip-1pt\left( \frac{ w-2int+i(\alpha_k-\alpha_j }{2}  \right)        dw\\
   &
   \\
   &\hskip 50pt\ll \; \left( \log\log R\right)^{2L-1} \cdot \prod\limits_{j=1}^n \prod\limits_{k=1}^n  
 \hskip-2pt 
  \left|\Gamma\hskip-1pt\left( \frac{ 2int+i(\alpha_k-\alpha_j }{2}  \right) \right|^2.
   \endalign
   $$
  Here we have used the fact that if we shift the line of integration in the above integral then the growth in $|t|$ does not change because of  $(p\delta)^{-w} \approx (t^{n^2})^{-w}$  which cancels the polynomial term (coming from Stirling's formula)  in the product of Gamma functions.
   
   \vskip 5pt

      It follows  from (12.5), (12.9), (12.10) that
   $$I_{1,1}^*(p) \; \gg \; (\log R)^{2L-2}\cdot |t|^{-n^2}, \qquad (|t|\to\infty),\tag{12.11}$$
    where the $\gg$-constant depends at most on $n$ and the $\alpha_i \;(i = 1,2,\ldots ,n).$

   \vskip 10pt
  \noindent
  {\bf Lemma 12.12} {\it Let $f$ be a Hecke Maass cusp form for $SL(n, \Bbb Z)$ which is tempered at every rational prime. Then there exist at least $\frac{1}{10n^2}\frac{N}{\log N}$ primes $p$ in the interval $[N, 2N]$ such that
  $$|\lambda(p)| \; \ge \frac{1}{100},$$
  where $\lambda(p)$ is the $p$-th Hecke eigenvalue of $f$.}
  
  \vskip 10pt
  {\bf Proof:} By \cite{Liu, Wang, Ye, 2005} it is known that
  $$\sum_{N\le m\le 2N} \Lambda(m) \cdot |\lambda(m)|^2 \; \sim\; N.\tag{12.13}$$
   
   Define a prime to be good if $|\lambda(p)| \ge \frac{1}{100}$, and otherwise define the prime to be bad.
    Now, suppose the conclusion of lemma 12.12 is wrong. We will show  that this leads to a contradiction.  If lemma 12.12 is false, the left side of (12.13) can be asymptotically estimated as follows:
    $$\align \sum_{N\le m\le 2N} \Lambda(m) \cdot |\lambda(m)|^2\; & \sim \; \underset p\;\text{good}\to{\sum_{N\le p\le 2N}}\Lambda(p) \cdot |\lambda(p)|^2 \; + \; \underset p\;\text{bad}\to{\sum_{N\le p\le 2N}}\Lambda(p) \cdot |\lambda(p)|^2
    \\
    &
    \le \; n^2 \frac{1}{10n^2}N  \;+ \;  \sum_{N\le p\le 2N} \Lambda(p) \frac{1}{100^2}\tag{12.14}
    \\
    & \le \frac{N}{10} \; + \; \frac{N}{100^2}.\endalign$$
    Here we used the Ramanujan bound
    $|\lambda(p)| \le n,$  i.e., the fact that $f$ is tempered at $p$. Since (12.14) contradicts (12.13) this proves the lemma.  \qed
    
    \vskip 10pt\noindent
    {\bf Lemma 12.15} {\it Let $N \ge t^2\ge 1.$ Then there exist at least $\left(1 - \frac{1}{20n^2}\right) \frac{N}{\log N}$ primes $p$ in the interval $[N, 2N]$ such that
    $$\big |\eta_{nit}(p)\big | \; \ge \; \frac{1}{2000^2 n^4}.$$}
    
    {\bf Proof:} We have
    $$\big|\eta_{nit}(p)\big| \; = \; \left| p^{2nit} + 1\right| \; \ge \; \big|\cos(2nt) \,(\log p) + 1\big|.$$
   Let
   $$\Delta = \frac{1}{2000^2 n^4}.$$
   Assume for some integer $m$ and a prime $N \le p \le 2N$ that
   $$\big| 2nt(\log p) - (2m+1)\pi\big| \; \ge 10\sqrt{\Delta}.\tag{12.16}$$ 
   Then it follows that
   $$\big|\cos(2nt) (\log p) + 1\big| \; \ge \; \Delta.$$
   Define $S_m$ to be the set of primes $p$ in the interval $[N,2N]$ which don't satisfy (12.16). For $p \in S_m$ we have the inequalities:
   $$e^{  \frac{(2m+1)\pi}{ 2n|t| }   }\left( 1 - \frac{ \sqrt{\Delta}  }{2nt}   \right) \; \le \; p \; \le \; 
   e^{  \frac{(2m+1)\pi}{ 2nt }   }\left( 1+ \frac{ \sqrt{\Delta}  }{2n|t|}   \right).$$
   Sieve theory tells us \cite{Bombieri-Davenport, 1969} that for 
   $$M = e^{  \frac{(2m+1)\pi}{ 2nt }   } \frac{\sqrt{\Delta}}{n|t|}$$
   we have the following bound for the cardinality of $S_m$:
   $$\#S_m \le \frac{3M}{\log M}$$
   
   Now
   $$\frac{N}{2} \le  e^{  \frac{(2m+1)\pi}{ 2n|t|}   } \le 3N.$$
  This implies that
  $$\frac{2n|t|(\log N)}{2\pi} - \frac12 \; \le \; m \; \le \; \frac{2n|t|(\log N)}{2\pi} - \frac12 + \frac{n|t|(\log 2)}{\pi}.$$
  Hence
  $$\#\left(\bigcup_m S_m   \right) \; \le \; \frac{ 2N\sqrt{\Delta} }{ n|t| }\cdot \frac{n|t|(\log 9)}{\log N - \log |t|} \; \le \; \frac{100 N\sqrt{\Delta}}{\log N} \; = \; \frac{1}{20n^2}\frac{N}{\log N}.$$
  
  \qed
  
  \vskip 10pt\noindent
  {\bf Lemma 12.17} {\it Let $f$ be a Hecke Maass cusp form for $SL(n, \Bbb Z)$ which is tempered at every rational prime. Then there exist at least $\frac{1}{20n^2} \frac{N}{log N}$ primes $p$ in the interval $[N, 2N]$ such that
  $$\big|\lambda(p) \eta_{nit}\big| \gg_{n,f} \; 1.$$}
  
  {\bf Proof:} This follows immediately from lemmas 12.12 and 12.15 since the densities $\frac{1}{10n^2}\frac{N}{\log N}$ and $\left(1 - \frac{1}{20n^2}\right) \frac{N}{\log N}$ imply an overlap of positive density. \qed
  
  \vskip 10pt
  Now, it follows from  (12.5), (12.11), lemma 12.17, and the previous choice we made, by setting $N\delta   \sim \left(n |t|\right)^{n^2},$ that 
  $$\boxed{\big| L(1+2nit, \; f\times \widetilde f\,)\big|^2 \sum_{N\le p \le 2N} I_1^*(p) \; \gg \; \beta^n \frac{N}{\log N} |t|^{-n^2} \; \gg \; \frac{\beta^n \delta^{-1}}{\log |t|}}\tag{12.18}$$
  where the $\gg$-constant depends at most on $R,n,f$ and is independent of $|t|\to\infty.$
  From now on we assume
  $$L(1+2int, \; f\times \widetilde f\,) \;\ll\; \frac{1}{(\log |t|)^3}.$$
  Otherwise, we already proved our main theorem.
  
  Consider the contribution of the $(1,1,\ldots,1,p)$ coefficient of $E(z, f; 1+it)$ which we denote as
  $$a(p, 1+it) := \frac{c_n\, \overline{\lambda(p) \eta_{\frac{n}{2}+nit}(p)}    }{  p^{\frac{2n-1}{2}}    } \,W_{2n,\nu'}\left(My\right) \cdot \frac{1}{L(n+1+2nit, \; f\times \widetilde f)}.$$ 
  
 Let $s = 1+it.$  The Langlands parameters of $W_{2n, \nu'}$ are
  $$\left(ns-\frac{n}{2}+i\alpha_1, \quad\cdots\quad ,ns-\frac{n}{2}+i\alpha_n, \; -ns+\frac{n}{2}+i\alpha_1, \quad \cdots\quad ,-ns+\frac{n}{2}+i\alpha_n\right).$$
  Define 
  $$\align I_2^*(p) & := \int\limits_0^\infty\cdots \int\limits_0^\infty \left| \int\limits_0^\infty a(p, 1+it) g\left(\frac{A}{\beta}\right) \; \frac{dA}{A}\right|^2 \psi\left(  \frac{\text{Det} \; z}{\delta}  \right) \; d^\times y\\
  & = |\lambda(p)|^2 \left|\eta_{\frac{n}{2}+nit}(p)\right|^2\cdot \frac{|c_n|^2}{|L(n+1+2nit, \; f\times \widetilde f\,)|^2}\; I_{2,1}^*(p)\endalign$$
  where for $c > n$ and $\alpha_{k,j} := \alpha_k-\alpha_j$,
   we set
   $$\align I_{2,1}^*(p)   & :=   \frac{1}{2\pi i} \int\limits_{c-i\infty}^{c+i\infty}  \widetilde\psi(-s') \big(p\delta\big)^{-s'}
    \int\limits_0^\infty \cdots  \int\limits_0^\infty \left|W_{2n,\nu'}(y)\right|^2 \left( \text{Det} z  \right)^{s'} \; d^\times y \; ds'
  \\
   &
   \\
   & :=  \frac{1}{2\pi i}
    \int\limits_{c-i\infty}^{c+i\infty}  \widetilde\psi(-s') \big(p\delta\big)^{-s'}
   \frac{\pi^{n^3+2n^2-n(2n-1)s'} 2^{  \frac{ (2n-1)2n(2n+1)   }{6}  }     }{  \Gamma(ns')} \prod_{k=1}^n \prod_{j=1}^n
   \\
   &
   \hskip -25pt
  \cdot  \frac{  \Gamma\left( \frac{s'+n(s+\bar s)-n+i \alpha_{k,j}}{2}  \right)   \Gamma\left( \frac{s'+n(s-\bar s)+i \alpha_{k,j}}{2}  \right)      \Gamma\left( \frac{s'+n(\bar s-s)+i \alpha_{k,j}}{2}  \right) 
  \Gamma\left( \frac{s'-n(s+\bar s)+n+i \alpha_{k,j}}{2}  \right)} 
 {  \left|\Gamma\left( \frac{1+i \alpha_{k,j}}{2}  \right) \Gamma\left( \frac{1+2ns-n+i \alpha_{k,j}}{2}  \right)\right|^{2}   }\; ds'
 \\
 &\,
 \ll \; |t|^{(-n+c-1)n^2} \big(p\delta\big)^{-c}.
    \endalign$$
    It follows that
    $$\big|L(1+2nit, \; f\times \widetilde f\,)\big|^2 \sum_{N\le p\le 2N} I_2^*(p)
    \; \ll \frac{\delta^{-1-n}}{(\log |t|)^2} \; \ll \; \frac{\beta^n\delta^{-1}}{(\log |t|)^2},$$
    upon recalling that $\beta = \delta^{-1}$.
    
  Similarly one shows that the contribution from
  $$\frac{ \Lambda(2ns-2n, \; f\times \widetilde f\,)   }{  \Lambda(1+2ns-2n, \; f\times  \widetilde f\,)  }\; E(z,  f; 2-s)$$
  is at most $ \beta^n\delta^{-1}\big/ (\log |t|)^2$.
 \vskip 8pt 
 For $c > 0$, let
  $$E_4(z,s) := \frac{1}{2\pi i}\int\limits_{c-i\infty}^{c+i\infty}  \frac{ A^{-\frac{n}{2} w}  }{w(w+1)} E\left(z, f; \,s+\frac{w}{2}\right)\; dw.$$

  Define $a^*(p, 1+it)$ to be the $(1,1,\ldots,1,p)$ coefficient of $E_4(z,s),$ which is given by
  $$a^*(p, 1+it) = \frac{1}{2\pi i}\int\limits_{c-i\infty}^{c+i\infty}  \frac{ A^{-\frac{n}{2} w}  }{w(w+1)}
  \frac{c_n\, \overline{\lambda(p) \, \eta_{ns+\frac{nw}{2}-\frac{n}{2}} (p)}  }{  p^{ \frac{2n-1}{2}   }  }  \frac{W_{2n,\nu''}(My) \; dw}{L\left(1+2n(s+\frac{w}{2})-n, \; f\times \widetilde f\,\right)}$$
  and the associated Langlands parameters are
  $$\align & \bigg(n\left(s+\frac{w}{2}\right) - \frac{n}{2}+i\alpha_1, \quad\ldots\quad
   ,n\left(s+\frac{w}{2}\right)-\frac{n}{2}+i\alpha_n, \; -n\left(s+\frac{w}{2}\right)+\frac{n}{2}+i\alpha_1,
  \\
  &
 \hskip 228pt  \quad \ldots \quad,-n\left(s+\frac{w}{2}\right)+\frac{n}{2}+i\alpha_n\bigg).\endalign    $$
For $c > 1,$ let
 $$\align I_4^*(p)   & := \int\limits_0^\infty\cdots  \int\limits_0^\infty \left| \int\limits_0^\infty a^*(p, 1+it) g\left(\frac{A}{\beta}\right) \; \frac{dA}{A}     \right|^2 \; \psi\left(\frac{ \text{Det}\; z }{ \delta } \right) \; d^\times y\\
 &
 \\
 &\hskip-15pt \ll    \int\limits_0^\infty\cdots  \int\limits_0^\infty 
 \;\left|\; \int\limits_{c-i\infty}^{c+i\infty} \hskip-1pt \frac{ \beta^{-\frac{n}{2} w}  }{w(w+1)}  \widetilde g(-w) 
 \frac{ \overline{\lambda(p) \,\eta_{n(s+\frac{w}{2})-\frac{n}{2}}(p)} }{ p^{\frac{2n-1}{2}} } 
\frac{ W_{2n,\nu"}(My) \; dw }{L\left(1+2n\left(s+\frac{w}{2}\right)-n, \; f\times \widetilde  f\,\right)}\right|^2\\
&
\hskip 295pt \psi\left( \frac{\text{Det} \;z}{\delta}  \right) \; d^\times y\\
&
\\
&\hskip-10pt
\ll \int\limits_0^\infty\cdots \int\limits_0^\infty \int\limits_{c-i\infty}^{c+i\infty} \, \left| \;\beta^{-\frac{nc}{2}} \widetilde g(-w) \;  \frac{ \overline{\lambda(p) \,\eta_{n(s+\frac{w}{2})-\frac{n}{2}}(p)} }{ p^{\frac{2n-1}{2}} } \, W_{2n,\nu"}(My)\; \right|^2 dw \cdot \psi\left( \frac{\text{Det} \;z}{\delta}  \right) \, d^\times y
\\
&
\\
&
\ll  \int\limits_{c-i\infty}^{c+i\infty} \beta^{-nc} \left|\widetilde g(-w)\right |^2 \; |\lambda(p)|^2 \; \left|  \eta_{n(s+\frac{w}{2})-\frac{n}{2}}(p)  \right|^2 \; I_{4,1}^*(p)\; dw,\endalign$$
where, for $c' > n+nc$ and $\alpha_{k,j} = \alpha_k-\alpha_j$, we have
$$\align I_{4,1}^*(p) & := \int\limits_0^\infty \cdots \int\limits_0^\infty \Big | W_{2n,\nu"}(y)\Big |^2 \; \psi\left( \frac{\text{Det} \;z}{p \,\delta}  \right) \, d^\times y\\
&
= \frac{1}{2\pi i}\int\limits_{c'-i\infty}^{c'+i\infty} \widetilde\psi(-s') \;(p\delta)^{-s'}  
\int\limits_0^\infty \cdots \int\limits_0^\infty\; \Big|W_{2n,\nu"}(y)\Big |^2  
\left(\text{Det}\; z\right)^{s'} \; d^\times y \; ds'\\
& = \frac{1}{2\pi i} \int\limits_{c'-i\infty}^{c'+i\infty} \widetilde\psi(-s') \;(p\delta)^{-s'} \; \frac{ \pi^{n^3+2n^2+n^3c-n(2n-1)s'}  2^{ \frac{(2n-1)2n(2n+1) }{6}  }}{ \Gamma(ns') }\\
&
\hskip -18pt
\cdot  \prod_{k=1}^n \prod_{j=1}^n  \Gamma\left( \frac{s'+n\left(s +\frac{w}{2}  +\bar s +\frac{\bar w}{2}-1   \right)+i\alpha_{k,j}}{2}  \right) \Gamma\left( \frac{s'+n\left(s +\frac{w}{2}  -\bar s -\frac{\bar w}{2}   \right)+i\alpha_{k,j}}{2}  \right)  
   \\
   &\hskip 50pt
  \cdot  \frac{       \Gamma\left( \frac{s'+n\left(\bar s +\frac{\bar w}{2}  - s -\frac{w}{2}   \right)+i\alpha_{k,j}}{2}  \right) 
   \Gamma\left( \frac{s'-n\left(s +\frac{w}{2}  +\bar s +\frac{\bar w}{2}-1   \right)+i\alpha_{k,j}}{2}  \right) } 
 {  \left|\Gamma\left( \frac{1+ i\alpha_{k,j}}{2}  \right) \Gamma\left( \frac{1+2n(s+\frac{w}{2})-n+ i\alpha_{k,j}}{2}  \right)\right|^{2}   }\; ds'
 \\
 &
 \\
 & \;\ll |t|^{-n^3-n^3 c + n^2c'-n^2} (p\delta)^{-c'}.
\endalign
$$
 It follows that
    $$\big|L(1+2nit, \; f\times \widetilde f\,)\big|^2 \sum_{N\le p\le 2N} I_4^*(p) \;
    \ll  \frac{ \delta^{-n-nc}\, \beta^{-nc} \,\delta^{-1} }{(\log |t|)^2   }          \; \ll \; \frac{\beta^n\delta^{-1}}{(\log |t|)^2},$$
    upon recalling that $\beta = \delta^{-1}.$
  Similarly one shows that the contribution from
  $$\frac{ \Lambda(2ns-2n, \; f\times \widetilde f\,)   }{  \Lambda(1+2ns-2n, \; f\times \widetilde f\,)  }\;\int\limits_{1-i\infty}^{1+i\infty} \frac{A^{-\frac{n}{2}w}}{w(w+1)} \;E\left(z, f; \, 2-s+\frac{w}{2}\right) \; dw$$
  is at most $ \beta^n\delta^{-1}\big/ (\log |t|)^2$.  
  \vskip 8pt
  
  Combining all of the above, we finish the proof of theorem 12.1, where we see that the main contribution to the lower bound comes from (12.18).
  
  \qed
  
  \vskip 25pt\noindent
 {\bf Acknowledgement:}
\vskip 8pt
We would like to thank Lei Zhang for helping us determine the correct scaling factor in the constant term formula for Eisenstein series. We are particularly grateful  to Steven D. Miller for carefully reading the manuscript and correcting  several technical errors.  We also thank E. Lapid  for many helpful comments and encouragement. Especially, the second named author would like to thank P. Garrett for explaining his online notes patiently and answering many of her questions. The first named author is supported by  Simons Collaboration Grant 317607 while the second named author is supported by NSF grant DMS64380.
   \vskip 25pt\noindent
 {\bf Bibliography}
 \vskip 10pt
 
 \noindent
  Arthur, James:, {\it A trace formula for reductive groups, II,}  Applications of a truncation operator. Compositio Math. 40 (1980), no. 1, 87Ð121.
 
 \vskip 5pt\noindent  
  Bombieri E.;  Davenport H,; {\it On the large sieve method, Number theory and Analysis,} Plenum, New York, (1969),  9-22.

  \vskip 5pt\noindent 
  Brumley, F.; {\it Effective multiplicity one on $GL(N)$ and narrow zero-free regions for Rankin-Selberg L-functions,} Amer. J. Math. 128 (2006), no. 6, 1455-1474.
  
   \vskip 5pt\noindent 
  Brumley, F.; {\it Lower bounds on Rankin-Selberg L-functions,}  Appendix to Lapid's paper, {\it On the Harish-Chandra Schwartz space of $G(F)\backslash G(\Bbb A),$} Tata Inst. Fundam. Res. Stud. Math., 22, Automorphic representations and L-functions,  Tata Inst. Fund. Res., Mumbai, (2013),  335--377. 
  \vskip 5pt\noindent 
  Brumley, F.; Templier, N,; {\it Large values of cusp forms on $GL(n)$,} \hfil\break
  www.math.cornell.edu/$\sim$templier

   \vskip 5pt\noindent
   Garrett, P.; {\it Truncation and Maass-Selberg relations,} (2005)
   
   \noindent
 http://www.math.umn.edu/$\sim$garrett/m/v/maass$\underline{\phantom{x}}$selberg.pdf  
    \vskip 5pt\noindent
  Gelbart, S.; Lapid, E.; {\it Lower bounds for L-functions at the edge of the critical strip,} Amer. J. Math. 128 (2006), no. 3, 619-638.
  
  \vskip 5pt\noindent
  Goldfeld, D.; {\it Automorphic forms and L-functions for the group $GL(n, \Bbb R)$,}   Cambridge Studies in Advanced Mathematics, 99, Cambridge University Press, Cambridge, (2006).
  
   \vskip 5pt\noindent
   Ichino, A.; Yamana, S.; {\it Periods of automorphic forms: the case of $(GL_{n+1}\times GL_n, GL_n)$,} Compos. Math., to appear.
   
    \vskip 5pt\noindent
   Ichino, A.; Yamana, S.; {\it Periods of automorphic forms: the case of $(U_{n+1}\times U_n, U_n)$,} www.math.kyoto-u.ac.jp/$\sim$ichino
  
  \vskip 5pt\noindent
Iwaniec, H.; Kowalski, E.;  {\it Analytic Number Theory,} volume 53 of American Mathematical Society Colloquium Publications, American Mathematical Society, Providence, RI, (2004).

 \vskip 5pt\noindent
 Iwaniec, H.; Sarnak, P.; {\it Perspectives on the analytic theory of L-functions,} GAFA 2000 (Tel Aviv, 1999). Geom. Funct. Anal. 2000, Special Volume, Part II, 705-741.
 
  \vskip 5pt\noindent
Jacquet, H.;
 {\it Automorphic forms on GL(2), Part II,} in Lecture Notes in Mathematics, vol.
278, Berlin-Heidelberg-NewYork, (1972).

 \vskip 5pt\noindent
Jacquet, H.; Lai, K. F.; {\it  A relative trace formula,} Compositio Math. 54 (1985), no. 2, 243-310.

\vskip 5pt\noindent
Jacquet, H.; Shalika, J.A.; {\it A non-vanishing theorem for zeta functions of $GL_n$}, Invent. Math. 38 (1976/77), 1-16.

\vskip 5pt\noindent
Langlands, R.P.; {\it  Euler products,}  A James K. Whittemore Lecture in Mathematics given at Yale University, 1967, Yale Mathematical Monographs, 1, Yale University Press, New Haven, Conn.-London, (1971).

\vskip 5pt\noindent
Labesse, J.P.;  Waldspurger, J.;
{\it La formule des traces tordue d'apr\'es le Friday Morning Seminar,}  [The twisted trace formula, from the Friday Morning Seminar],  CRM Monograph Series, 31, American Mathematical Society, Providence, RI, (2013).

\vskip 5pt\noindent
Lapid, E.M.; {\it On the fine spectral expansion of JacquetÕs relative trace formula,} J. Inst. Math. Jussieu 5 (2006),
no. 2, 263-308.

 \vskip 5pt\noindent
Liu, J.; Wang, Y.; Ye, Y.; {\it A proof of Selberg's orthogonality for automorphic L-functions,}  Manuscripta Math. 118 (2005), no. 2, 135-149.

 \vskip 5pt\noindent
 Moglin, C.; Waldspurger, J.L.; {\it P\^oles des fonctions $L$ de paires pour $GL(N)$,} appendix to Le spectre r\'esiduel de $GL(n),$ Ann. Sci. ENS (4\` eme s\'erie) 22 (1989) 605-674.

 \vskip 5pt\noindent 
  Moreno, C.; {\it Analytic proof of the strong multiplicity one theorem,} Amer. J. Math. 107 (1985), no. 1, 163-206.
  
  \vskip 5pt\noindent
 Sarnak, P.;  {\it  Nonvanishing of L-functions on $\Re(s)=1$,} Contributions to automorphic forms, geometry, and number theory, 719Ð732, Johns Hopkins Univ. Press, Baltimore, MD, (2004).
 
 \vskip 5pt\noindent
 Shahidi, F.; {\it On certain L-functions,} Amer. J. Math. 103 (1981), 297-355.
 
 \vskip 5pt\noindent
Shahidi, F.; {\it  Eisenstein series and Automorphic L-functions,}  AMS Colloquium Publications, 58, American Mathematical Society, Providence, RI, (2010).

 \vskip 5pt\noindent
  Stade, E.; {\it On explicit integral formulas for GL(n,R)-Whittaker functions,}  Duke Math. J. 60 (1990), no. 2, 313-362.
  
   \vskip 5pt\noindent
  Stade, E.; {\it Archimedean L-factors on $GL(n)\times GL(n)$ and generalized Barnes integrals,} Israel Journal of math., 127 (2002), 201-219.
  
  \vskip 5pt\noindent
  Whittaker, E. T.; Watson, G. N.; {\it A course of modern analysis, An introduction to the general theory of infinite processes and of analytic functions; with an account of the principal transcendental functions.} Reprint of the fourth (1927) edition, Cambridge Mathematical Library, Cambridge University Press, Cambridge, (1996).

  \enddocument